\documentclass[a4paper,10pt]{article}
\def\1{{\bf 1}}
\def\R{{\bf R}}
\def\N{{\bf N}}
\def\Д{{Доказательство.}}
\def\s{{\bf S}}
\def\eps{\varepsilon}
\def\e{\varepsilon}
\def\Z{{\bf Z}}
\def\T{{T}}

\def\F{{\cal F}}
\def\J{{\cal J}}
\def\jj{{\bf J}}
\def\P{{\cal P}}
\def\A{{\cal A}}
\def\B{{\cal B}}
\def\M{{\cal M}}

\usepackage[T1]{fontenc} 
\usepackage[cp1251]{inputenc} 
\usepackage[russian]{babel}

\begin{document}
\large


\thispagestyle{empty}
\begin{center}
Московский государственный университет  \\
имени М.В.Ломоносова
\end{center}
\vspace{20mm}

\vspace{20mm}
\begin{center}
{\Large \bf Рыжиков Валерий Валентинович}
\end{center}
\vspace{20mm}
                         {\Large
\begin{center}
{\bf \Large Марковские сплетающие операторы, джойнинги}
\end{center}
\begin{center}
{\bf \Large и асимптотические свойства динамических систем}
\end{center}
\
                        }
\begin{center}
01.01.01 -- математический анализ\\
\vspace{10mm}
Диссертация \\ на соискание ученой степени\\ доктора физико-математических
наук \\
\end{center}
\vspace{30mm}
\begin{center}
Москва 2004
\end{center}
\newpage

\bf ВВЕДЕНИЕ. \rm
\\{0.1.    Проблема Рохлина о кратном перемешивании} \hfill {3}
 \\{ 0.2.    Теория джойнингов и ее приложения} \hfill {5}
 \\{ 0.3.    Теория марковских сплетающих операторов}\hfill {8}
 \\{ 0.4.    Структура и основные результаты диссертации} \hfill{12}

\bf Глава 1. 
   МАРКОВСКИЕ СПЛЕТАЮЩИЕ ОПЕРАТОРЫ  И 

ТЕНЗОРНАЯ ПРОСТОТА   ДИНАМИЧЕСКИХ СИСТЕМ \rm 
\\
1.1. Несколько методологических принципов теории сплетений  \hfill{22}
 \\{ 1.2.  Дополнительная симметрия} \hfill{27}
 \\{ 1.3.  Индуцированные джойнинги} \hfill{31}
 \\{ 1.4.  Примеры тензорно простых систем} \hfill{33}
 \\{ 1.5.  Связь типов тензорной  простоты} {38}

\bf Глава 2.  СИСТЕМЫ С МИНИМАЛЬНЫМ, ПРОСТЫМ 

 И КВАЗИПРОСТЫМ  ЦЕНТРАЛИЗАТОРОМ \rm 
 \\{ 2.1.  Простые  системы с несчетным централизатором}\hfill {43}
\\  2.2.  Наследственная независимость и квазипростота действий\hfill{45}
 \\ 2.3.  Минимальные самоприсоединения  и  кратная         возвращаемость\hfill{49}
 \\{ 2.4.  Четная и нечетная тензорная простота}\hfill {55}

\bf Глава 3. ДЖОЙНИНГИ И ТЕНЗОРНАЯ ПРОСТОТА 

 НЕКОТОРЫХ ПОТОКОВ.  \rm                   
 \\{ 3.1.  Гладкие джойнинги и  тензорная простота потоков}\hfill{58}
 \\{ 3.2.  Тензорная простота  $\omega$-простых   потоков}\hfill {59}
\\ 3.3.  Перемешивающие  потоки положительного          локального  ранга\hfill{73}

\bf Глава 4. ДЖОЙНИНГИ ДЕЙСТВИЙ КОНЕЧНОГО 

 И ПОЛОЖИТЕЛЬНОГО ЛОКАЛЬНОГО РАНГА \rm
\\
 4.1.  D-свойство перемешивающих автоморфизмов
        конечного  ранга  \hfill{68}
\\ 4.2.  D-свойство  перемешивающих $\Z^n$-действий
 и локальный ранг\hfill{71}
\\ 4.3.  Тензорная простота перемешивающих систем                                   
 с D-свойством\hfill{73}
 \\{ 4.4.  Кратное перемешивание и локальный ранг}\hfill{77}
 \\{ 4.5.  Ранги и джойнинги  автоморфизма $T\times T$}\hfill{78}

\bf Глава 5. НЕКОТОРЫЕ СПЕКТРАЛЬНЫЕ,  АЛГЕБРАИЧЕСКИЕ  И 

 АСИМПТОТИЧЕСКИЕ СВОЙСТВА  ДИНАМИЧЕСКИХ СИСТЕМ \rm
 \\{ 5.1.  Проблема Рохлина об однородном спектре}\hfill{82}
 \\  5.2.  Перемешивающие автоморфизмы с однородным
          \\{ \ \ \ \ \ \ непростым спектром}\hfill{86}
\\ 5.3.  Изоморфизм  декартовых степеней преобразований
 \\{ \ \ \ \ \ \ \ и $ \kappa$-перемешивание}\hfill{89}
\\ 5.4.  Асимметрия прошлого и будущего  системы
          \\{ \ \ \ \ \ \ \ и кратная возвращаемость}\hfill{92}
\\ 5.5.  Расширения, сохраняющие кратное перемешивание
          \\{ \ \ \ \ \ \ \ и тензорную простоту    }\hfill{95}

{\bf ЛИТЕРАТУРА}\hfill {99}                                            
\newpage

\begin{center}    { \Large  \bf ВВЕДЕНИЕ }
\end{center}

\medskip
\medskip
\begin{center}
{\bf  \  0.1.  Проблема Рохлина о кратном перемешивании }
\end{center}
\medskip
Основным объектом исследования в диссертации является обратимое
сохраняющее меру $\mu$  преобразование $T$ пространства
Лебега $(X,\B,\mu)$, которое
называют автоморфизмом.
 Динамической системой  называется
четверка $(T, X,\B,\mu)$
или, в более общей ситуации,  сохраняющее меру  действие некоторой группы.
Среди свойств $T$, которые представляют интерес для эргодической теории,
особую роль играют асимптотические свойства (свойства систем для больших
значений времени).
Рассмотрим пример такого свойства,
который является ключевым для нашей работы.

{\bf Кратное перемешивание. }
Говорят,  что  автоморфизм $T$  перемешивает с кратностью $k$,
если для любых множеств $A_0,\dots,A_k\in \B $
и любых последовательностей $n_1, \dots ,n_k\to \infty $ выполнено:
$$
\mu(A_0\cap T^{n_1}A_1 \dots \cap T^{n_1+\dots +n_k}A_n)
\to  \mu (A_0)\mu (A_1)\dots \mu (A_k).
$$

В.А. Рохлин  в работе \cite{Rok} ввел понятие кратного перемешивания
и доказал, что  эргодический
эндоморфизм компактной
коммутативной группы обладает кратным премешиванием всех порядков.
Проблема эквивалентности свойств перемешивания разных порядков,
получившая название проблемы Рохлина о кратном перемешивании, стала
популярна  после выхода книги  Халмоша \cite{Hal}.
Напомним  историю результатов.

В.П. Леонов \cite{Leo}
показал, что перемешивающие гауссовские системы обладают
перемешиванием всех кратностей.
\\ Ф. Ледраппье \cite{Led} обнаружил контрпример к
проблеме о кратном
перемешивании  для действий группы ${\bf Z}^2$.
Он построил перемешивающее
действие группы $\Z^2$, которое не обладает перемешиванием
кратности 4. Это действие образовано коммутирующими сдвигами
(автоморфизмами)
подгруппы  $H$  группы $2^{\Z^2}$, где
$H$ состоит из всех
последовательностей $\{h(z)\}$,  $z\in \Z^2, \ h(z)\in \Z_2$,
таких, что  условие
$$
       h(z_1+1,z_2)+ h(z_1,z_2+1)+ h(z_1-1,z_2)+
       h(z_1,z_2-1)= h(z_1,z_2)
$$
выполнено для всех $z=(z_1,z_2)$.
Идея Ледраппье позволяет варьировать результаты:
для каждого  $k$ найдется коммутативное действие,
обладающее перемешиванием кратности
$k$, но не обладающее перемешиванием кратности    $(k+1)$.
Однако проблема Рохлина о кратном премешивании,
поставленная для
$\Z$-действий, остается открытой более полувека.

Упомянем результаты, дающие  положительный ответ  для некоторых классов
динамических систем.
Я.Г. Синай высказал гипотезу о том, что орициклический поток
является перемешивающим всех степеней, которую подтвердил  Б. Маркус,
 доказавший более общее  утверждение:
свойством кратного перемешивания обладают унипотентные потоки.
 Ряд обобщений теоремы Маркуса был получен позднее
в \cite{flow}, где автор
применил метод джойнингов, а также  Ш. Мозесом \cite{Moz} и
А.Н. Старковым \cite{Star},\cite{Stara}.  Так, например, в
 \cite{Star} доказано
свойство кратного перемешивания для однородных перемешивающих потоков.
Ряд общих результатов и наблюдений  о кратном перемешивании
получены авторами \cite{UMN89}, \cite{FTh}, \cite{R85},  \cite{Sch}.
Проблема Рохлина о кратном перемешивании допускает модификации.
Например, влечет ли слабое перемешивание  за собой
слабое перемешивание  всех порядков?
\cite{UMN89}.

Один из наиболее общих результатов принадлежит Б. Осту
\cite{Hos}:
перемешивающие автоморфизмы с сингулярным спектром не допускают
нетривиальных самоприсоединений с парной независимостью и по этой
причине обладают перемешиванием бесконечной кратности.
Вывод из теоремы Оста: контрпримеры к проблеме Рохлина
следует искать в классе систем с быстрым перемешиванием кратности 1.

В  \cite{Kal} С. Каликов   установил свойство перемешивания
кратности 2
для перемешивающих  автоморфизмов ранга 1.
   Результат Каликова был
несколько неожиданным, так как здесь свойство кратного перемешивания
получено для систем со слабыми статистическими свойствами.
Автор в \cite{R0} привел
обобщение теоремы Каликова для всех кратностей,
основанное на технике джойнингов,
показав, что перемешивающие автоморфизмы ранга 1 не допускают
самоприсоединения с парной независимостью. В диссертации
 эквивалентное свойство, сформулированное в терминах сплетающих операторов,
 называется тензорной простотой (определение
приводится ниже). Но интерес к этому свойству связан не только  с тем, что
тензорная простота перемешивающей системы влечет за собой
кратное перемешивание.

 Тензорную простоту
можно рассматривать как аналог свойства
взаимной сингулярности спектральной меры автоморфизма
и ее сверточного квадрата (первое указание на это появилось в работе
 Оста \cite{Hos}).
 Упомянутое  спектральное   свойство
в относительном варианте было обнаружено А.М. Степиным \cite{St66}
для  групповых действий  при решениии
проблемы Колмогорова о групповом свойстве спектров динамических систем
(этой проблеме посвящены  также работы В.И.Оселедца \cite{Os} и
А.М.Степина \cite{St86}).
Таким образом, свойство тензорной простоты, появившееся внутри теории
джойнингов \cite{JR}, оказалось замечательным образом
 связанным с проблемами Колмогорова и Рохлина.

Предположим, что контрпример к проблеме Рохлина найден. Тогда
можно задать меру $\nu$ на кубе $X^{n+1}$, определяя значения
$
\nu(A_0\times A_1 \dots \times A_n)$
как предел выражений вида
$\mu(A_0\cap T^{n_1}A_1 \dots \cap T^{n_1+\dots +n_k}A_n).
$
Такая мера является самоприсоединением: она
 инвариантна относительно прямого произведения
$T_{(0)}\times T_{(1)} \dots \times T_{(n)}$,
а ее проекции на двумерные грани куба $X^{n+1}$ стандартны, т.е.
совпадают с мерой $\mu\otimes\mu$. Говорим, что такие самоприсоединения
обладают  попарной независимостью.
Причем мера $\nu$ нетривиальна, т.е.
$\nu\neq\mu\otimes\dots\otimes\mu$.

Хотя наш пример для действий группы $\Z$ является гипотетическим,
для упомянутого действия $\Z^2$ из работы  Ледраппье мы получим
нетривиальное самоприсоединение.

Если же будет доказано, что  рассматриваемая
система не допускает таких нетривиальных джойнингов, мы установим
кратное перемешивание (или слабое кратное перемешивание,
когда система обладала только слабым перемешиванием). В этом и состоит
подход в изучении проблемы Рохлина, использующий джойнинги.
\newpage
\begin{center}
{\bf   0.2. Теория джойнингов и ее  приложения}
\end{center}

 Понятие джойнинга  возникло в работе Фюрстенберга
\cite{Fur}, где он ввел понятие дизъюнктности двух систем
($T$ и $S$ дизъюнктны, если $\mu\otimes\mu$ -- их единственный
джойнинг) и доказал дизъюнктность К-автоморфизма
с автоморфизмом нулевой энтропии. Впрочем, этот факт вытекает из
теоремы
Пинскера \cite{Pi}: К-фактор и фактор с нулевой энтропией независимы.

Толчком к развитию теории джойнингов и их приложений
послужила  статья
  Д. Рудольфа \cite{Rud}, в которой построен
 автоморфизм со свойством  минимальных самоприсоединений. Приведем
 определение этого свойства.

Пусть $T:X\to X$ -- сохраняющий меру автоморфизм
вероятностного пространства $(X,\B,\mu),$ $\mu(X)=1.$
Мера  $\nu$, инвариантная относительно преобразования
  $T_{(1)}\times\dots\times T_{(n)}$, действующего в кубе
 $X_{(1)}\times\dots\times X_{(n)}$,
называется самоприсоединением   порядка $n$, если в дополнение к
сказанному выполнено условие: проекции меры  $\nu$  на
сомножитель $X_{(i)}$ совпадают с мерой  $\mu$.

Так, например,  мера  $\mu\otimes\mu$ и  образы $\Delta_{T^n}$  меры
$\mu$ при отображениях
$\varphi_n: X\to X\times X$, где $\varphi_n(x)=(x,T^n(x))$,
являются очевидными самоприсоединениями второго порядка.

Мера
$\Delta=\Delta_{Id}$ называется диагональной. Ее можно задать по-другому:
$\Delta(A\times B)=\mu(A\cap B)$ для всех $A,B\in \B$.

Автоморфизм $T$ пространства $(X,\B,\mu)$
 называется автоморфизмом с минимальными самоприсоединениями
порядка
$n$ (пишем  $T\in MSJ(n)$), если
любой эргодический  джойнинг  $n$ копий $T$,
исключая меру $\mu^{\otimes n}=\mu_{(1)}\otimes\dots\otimes\mu_{(n)}$,
обладает следующим свойством:
одна из его проекций  на двумерную грань  в
$X\times\dots\times X$ является мерой  $\Delta_{T^i}$.
(Неформально говоря,  такой автоморфизм $T$ имеет только
очевидные джойнинги.)

Используя автоморфизм $T\in MSJ$ как элемент конструкций,
можно построить разнообразные контрпримеры (примеры действий
с необычными свойствами).
Так, например, в  \cite{Rud}  приведены примеры
неизоморфных автоморфизмов $U$, $V$
таких, что  автоморфизм $U^n$ изоморфен автоморфизму $V^n$ для всех $n>1$,
даны примеры  автоморфизмов с несчетным семейством неизоморфных
квадратных корней,
построен  автоморфизм $U$ без корней такой, что $U^2$ имеет корни всех степеней.

В работе Рудольфа имеется ряд других примеров, из которых
мы приведем следующий пример
неизоморфных автоморфизмов $U$, $V$, которые слабо изоморфны.
Согласно определению  Синая \cite{Sin},
две системы $U$, $V$ слабо изоморфны,
если $U$ содержит $V$-фактор, а    система   $V$ имеет $U$-фактор,
т.е. факторсистему, изоморфную системе $V$.
Пусть  $T$ обладает свойством  минимальных самоприсоединений
всех порядков (свойством MSJ).
Рассмотрим $U=T\times T \times T \times \dots$.
Так как  автоморфизм $U$ коммутирует с  инволюцией $S$:
$$S(x_1,x_2,x_3,x_4\dots)= (x_2,x_1,x_3,x_4\dots),
$$
подалгебра измеримых множеств, инвариантных относительно $S$, будет
инвариантна относительно автоморфизма
$U$. Пусть  $V$ -- факторсистема, соответствующая инвариантной подалгебре.
Имеем представление
$V=T^{\odot 2} \times T \times T \times T \times \dots$, где
$T^{\odot 2}$ -- действие  $T\times T$ на подалгебре
множеств, инвариантных относительно отображения $(x_1,x_2) \to (x_2,x_1)$.
Из теории минимальных самоприсоединений вытекает, что
автоморфизмы  $U$ и $V$  не изоморфны. Слабый изоморфизм
$U$ и $V$ очевиден.

Отметим, что автоморфизм
$T\in MSJ(2)$ является слабо перемешивающим, имеет тривиальный
централизатор
(коммутирует только с $T^p$) и не имеет собственных факторов.
Если некоторый эргодический автоморфизм обладает набором
различных факторов, каждый из которых изоморфен некоторому
действию класса $MSJ(2)$,  то все эти факторы попарно независимы.

Другие примеры автоморфизмов и потоков с аналогичными свойствами
появились в работах дель Юнко, Парк, Рае, Ратнер, Свансон
\cite{Jun},\cite{JRS},\cite{JuP},\cite{Rat}. Этими авторами
установлено, что некоторые перекладывания отрезков,
автоморфизм Чакона,
специальные потоки над автоморфизмом Чакона,
некоторые орициклические потоки обладают свойством
минимальных самоприсоединений.
Упомянутые примеры автоморфизмов со свойством MSJ  принадлежат
классу автоморфизмов ранга 1 (в терминологии Катка и Степина \cite{KS}
-- допускают циклическую аппроксимацию).
  Приходько расширил множество примеров автоморфизмов
с минимальными самоприсоединениями:
вероятностными методами строятся автоморфизмы  бесконечного ранга
со свойством MSJ (см. \cite{Pri}).

Понятие простой системы обобщает
свойство минимальных самоприсоединений (MSJ).
Примеры простых систем и фрагменты теории
имеются в статьях    Вича \cite{Veech}, Вейса, Глазнера  \cite{GW},
Глазнера, Оста, Рудольфа \cite{GHR},
Рудольфа, дель Юнко \cite{JR},
Тувено \cite{Tho} и статьях других авторов.

Отметим, что  2-простая система является групповым расширением
каждого из ее факторов при условии нетривиальности фактора
 \cite{Veech}.
Недавно А. дель Юнко  построил пример простой системы,
не имеющей минимального нетривиального фактора  (это
контрастирует с отсутствием нетривиальных факторов у автоморфизма
класса MSJ).
Из результата Глазнера, Оста и Рудольфа \cite{GHR}
следует, что 2-простой слабо перемешивающий
автоморфизм, не являющийся простым порядка 3, должен обладать
перемешиванием. Эти авторы доказали, что
3-простота влечет $n$-простоту для $\Z$-действий.
Однако для действий некоммутативных групп имеются контрпримеры:
2-простота не совпадает с 3-простотой, а  последняя, вообще говоря,
не влечет за собой простоту всех порядков \cite{MZ96}.
Ситуация меняется при рассмотрении свойства простоты порядка 4.
Это свойство влечет за собой  простоту всех порядков
для любого группового  действия (см. работу Кинга \cite{King}).

Понятие свойства минимальных самоприсоединений   и
простоты обобщаются в разных направлениях (см.
 \cite{KiT}, \cite{Tho},\cite{2}), одним  из обобщений является 
квазипростота.
Отметим, что слабо перемешивающий
 автоморфизм, входящий в поток со свойством Ратнер (например,
в орициклический поток),
является квазипростым \cite{Rat},\cite{Tho}.
К.Парк \cite{Par} показала, что стандартное
$SL(2,Z)$-действие
автоморфизмов двумерного тора  ${\bf T}^2$   является квазипростым
действием порядка 2.
Эргодический  джойнинг этого действия сосредоточен на
подмножестве тора ${\bf T}^2\times {\bf T}^2$, которое
задается уравнением $mx=ny\ mod(1)$, где $x,y\in {\bf T}^2$,
$m,n$ - фиксированные натуральные числа.
Рассмотренное  действие не является квазипростым порядка 3,
так как равномерно распределенная мера на многообразии
$\{(x,y,z)\ :\ x+y+z=0\}$ является нетривиальным джойнингом
порядка 3.
Из  результатов цитированной работы Кинга  следует, что
квазипростота порядка 4 влечет за собой квазипростоту
всех порядков.

Результаты Парк  обобщены в статье Приходько \cite{Pr}, где
дано полное описание джойнингов группы автоморфизмов n-мерного
тора.  Эти работы вместе с \cite{PR},\cite{PR98} и \cite{MZ96} относятся
к теории джойнингов некоммутативных действий, которая контрастирует
с  теорией коммутативных действий.
\medskip


{\bf Джойнинги и спектр. }
В работе \cite{JL}  Леманчик и дель Юнко предложили
новый метод  построения
контрпримеров, использующий  вместо свойства минимальных
самоприсоединений свойство взаимной сингулярности
сверточных  степеней спектральной меры.
 Типичность последнего свойства,
доказанная   Степиным \cite{St86},  дает дополнительные
возможности (счетное пересечение типичных множеств непусто).
С позиций теории марковских сплетающих операторов эти подходы
при всей их оригинальности выглядят родственными: если спектральные
эффекты можно сформулировать в терминах сплетающих операторов,
то свойства джойнингов -- в терминах марковских сплетений.
Любопытно, что автоморфизм Чакона одновременно обслуживает
оба подхода: он обладает
свойством минимальных самоприсоединений (как мы упомянули выше)
и свойством взаимной сингулярности  сверточных  степеней спектра,
что недавно показали Приходько и автор \cite{PR}.

Связь между джойнингами и спектром тем сильнее,
чем  сингулярней  спектр.
 А системы с лебеговским спектром
могут быть (стохастически)  дизъюнктными, т.е.
не иметь марковских
сплетений, за исключением единственного тривиального сплетения, которое
функциям с нулевым     средним сопостовляет нулевую функцию.
Классический пример:  спектрально изоморфные  геодезический
и орициклический
потоки дизъюнктны, так как  геодезический поток является К-системой
(см.
\cite{An1}), орициклический является системой с нулевой энтропией
 (\cite{Gur}), а таковые дизъюнктны
(\cite{Pi},\cite{Fur}).
 То, что спектр этих потоков счетнократный лебеговский, было установлено
 в \cite{GeF} и    \cite{Para}.
Другой пример: спектрально изоморфные автоморфизмы
$T$ и $T^{-1}$
из \cite{Rud}
(стохастически)    дизъюнктны, т.е.
не имеют марковских
сплетений за исключением тривиального сплетения.
\newpage

\begin{center}
{\bf  0.3. Теория марковских сплетающих операторов }
\end{center}

Систематическое использование языка марковских сплетений было предпринято
автором, начиная с работы  \cite{preprint}.
Изложение метода сплетений также имеется в работах  Дж. Гудзона
\cite{G}
и  М. Леманчика, Ф. Парро, Ж.-П. Тувено \cite{LPT}.
 Связь между
марковскими операторами и полиморфизмами на декартовых произведениях
пространств с мерой описана в работе А.М. Вершика \cite{Ver}.
В его работе изучаются свойства самих полиморфизмов,
а в теории джойнингов --
свойства полиморфизмов, коммутирующих с динамической системой, или
в более общей ситуации -- полиморфизмов,  сплетающих
системы.
Полиморфизмом называется мера на
 $(X,{\cal B})\times (Y,{\cal B})$, где $X=Y$,
проекции которой на сомножители совпадают с $\mu$.
Полиморфизму $\nu$ соответствует оператор $P$, который
задается формулой:
$$
  Pf(y)= \int_{X}f(x) d\nu_y(x),
$$
где $\{\nu_y \ : \ y\in Y\}$ -- разложение меры $\nu$ в систему условных
мер  $\nu_y$.

Говорят, что мера $\eta$ есть джойнинг автоморфизмов
$T$ и $S$, если
 $\eta$ инвариантна относительно $T\times S$, а
проекции этой меры на сомножители в произведении
$({X},{\cal B})\times ({X},\cal{B})$  совпадают с мерой $\mu$.
Отвечающий джойнингу $\eta$ бистохастический оператор $P$
сплетает $T$ и $S$:
$PT=SP$.

Пусть $(T,X,{\cal B},\mu)$ -- динамическая система,
где $T$ обозначает  обратимое, сохраняющее меру $\mu$ преобразование
множества $X$ (фазового пространства), ${\cal B}$ - алгебра
$\mu$-измеримых множеств.
Преобразование $T$
называют автоморфизмом. Будем  обозначать тем же символом $T$
и  называть автоморфизмом унитарный оператор в $L_2(\mu)$,
отвечающий преобразованию $T$:
$(T f)(x)=f(Tx)$ для $f\in L_2(\mu)$.
Поскольку такие операторы сохраняют неотрицательность функций,
образованная ими группа ${\cal A}$ вложена в полугруппу
${\cal P}$    ограниченных операторов в $L_2(\mu)$, которые
переводят неотрицательные функции в неотрицательные.
 Оператору  $P\in {\cal P}$  соответствует мера $\nu$,
называемая квазиполиморфизмом (\cite{Ver}).
Связь задается формулой:
$$ \int_{X\times X }(f\otimes g) d\nu = \langle Pf,g\rangle,
$$
где $\langle\cdot,\cdot\rangle$ -- скалярное произведение в
$L_2(\mu)$ (иногда $\langle\cdot,\cdot\rangle$ также обозначает
скалярное произведение в  $L_2(\mu\otimes\dots\otimes\mu)$).

{\bf Сплетающие операторы, джойнинги и  кратное перемешивание.}
Особый интерес представляет случай, когда автоморфизм $S$ в формуле
сплетения $PT=SP$ изоморфен тензорной
степени автоморфизма  $T$.
Теория таких сплетений ( и теория джойнингов,
отвечающая этим операторам) имеет приложение к проблеме В.А.Рохлина
о кратном перемешивании.
Рассмотрим частный случай этой проблемы.
Пусть автоморфизм $T$ пространства Лебега перемешивает с кратностью
1, т.е.
для любых множеств $A,B\in{\cal B}$ при $n\to\infty$ выполнено
$$
\mu(T^nA\cap B) \to \mu(A)\mu(B).
$$
Будет ли автоморфизм перемешивать с кратностью 2?
Последнее означает,
что для любых измеримых множеств $A,B,C$ при любых последовательностях
$m,n\to\infty$  имеет место сходимость
$$
\mu(T^mA\cap T^{m+n}B\cap C) \to \mu(A)\mu(B)\mu(C) .
$$
Пусть $T$ перемешивает, но мы не знаем, является ли он перемешивающим
с кратностью 2. Предположим, что
 для некоторых последовательностей $m(i),n(i)\to\infty$
для любых $A,B,C$ выполнено
$$
\mu(T^{m(i)}A\cap T^{m(i)+n(i)}B\cap C) \to
\nu(A\times B\times C),
$$
где
$\nu$ -- мера (полиморфизм) на $X\times X\times X$.
Пусть
$A,B,C$ также обозначают индикаторы соответствующих множеств.  С мерой
$\nu$ связан
бистохастический оператор $P$, действующий из
$L_2(\mu \otimes \mu)$ в $L_2(\mu)$; связь задается равенством:
$$
\nu(A\times B\times C) =
\langle P(\chi_A\otimes \chi_B), \chi_C\rangle .
$$
Проверяется, что для меры $\nu$ выполнены следующие свойства:
$$
\nu(TA\times TB\times TC) =\nu(A\times B\times C)
$$
(мера инвариантна относительно $T\otimes T\otimes T$) и
$$
\nu(X\times A\times B) = \nu(A\times X\times B) = \nu(A\times B\times X)=
\mu(A)\mu(B)
$$
(проекции меры $\nu$ на грани декартова куба стандартны).
Для оператора $P$ точными аналогами этих свойств являются:
$$
TP = P(T\otimes T), \eqno                                              (0.1)
$$
$$
P(f\otimes \1) =
P(\1\otimes f) = Const = \1\otimes \Theta f,  \eqno          (0.2)
$$
где $\Theta$  обозначает
оператор ортопроекции на пространство констант в $L_2(\mu)$.
В дальнейшем тривиальным называется оператор $P$ такой, что
$Im(P)=\{Const\}$ (образ есть одномерное пространство постоянных
функций).

Итак возникают два объекта: мера $\nu$ и оператор $P$.
Если автоморфизм $T$ не обладает свойством перемешивания порядка 2,
то для некоторых последовательностей $m(i),n(i)\to\infty$ получим:
$$
\mu(T^{m(i)}A\cap T^{m(i)+n(i)}B\cap C) \to
\nu(A\times B\times C) \neq \mu(A)\mu(B)\mu(C)
$$
для некоторых $A,B,C$. Следовательно, мера $\nu$ и
оператор $P$ нетривиальны. Это означает, что
$$ \nu\ne\mu\otimes\mu\otimes\mu,\quad Im(P)\neq\{Const\} .$$

 Если существует единственный (тривиальный) оператор $P$,
удовлетворяющий  условиям (0.1),(0.2),
будем говорить, что $T$ обладает свойством $S(2,3)$.
Если перемешивающий автоморфизм $T$ удовлетворяет
свойству $S(2,3)$,
то $T$ будет обладать кратным перемешиванием порядка 2.
Действительно, оператор, отвечающий мере  $\mu\otimes (\mu\otimes \mu),$
удовлетворяет (0.1),(0.2). Из единственности
такого оператора вытекает свойство кратного перемешивания:
$$
\mu(T^{m(i)}A\cap T^{m(i)+n(i)}B\cap C) \to
\langle P(\chi_A\otimes \chi_B)\,, \chi_C\rangle  =
  \mu(A)\mu(B)\mu(C).
$$
{\bf Тензорная простота динамической системы.}
По аналогии со свойством $S(2,3)$ определим свойства $S(n,n+1)$:

существует единственный оператор $Q$, удовлетворяющий
условиям
$$TQ = Q(T_{(1)}\otimes T_{(2)}\otimes\dots\otimes T_{(n)}), \ \ (n>2)
     $$
 и
$$Q(f_{1}\otimes\dots\otimes f_{n-1}) =
Const =\int f_1\int f_2\dots\int f_{n-1},$$
если одна из функций $f_1, f_2,\dots, f_{n-1}$ является постоянной.
Для групповых действий  определения
свойств $S(n,n+1)$ аналогичны.
Если $(n-1)$-кратно перемешивающий автоморфизм $T$ удовлетворяет
свойству $S(n,n+1)$,
то $T$ будет обладать кратным перемешиванием порядка $n$.

\it Тензорно простым \rm называется действие класса $\cap_{n>2}S(n-1,n)$.

Как будет показано, свойство  $S(3,4)$
 влечет за собой каждое из свойств $S(n,n+1)$ и, следовательно,
свойство перемешивания всех порядков
для произвольного перемешивающего коммутативного действия.
Поэтому в диссертации тензорно простым также называется действие
класса $S(3,4)$.
Отметим, что оно эквивалентно следующему свойству:
множество самосопряженных бистохастических операторов $J$,
коммутирующих с $T\otimes T$ и удовлетворяющих условию
$$
  J(I\otimes\Theta) = J(\Theta\otimes I) = \Theta\otimes\Theta,
$$
одноэлементно, т.е. является множеством $\{\Theta\otimes\Theta\}$.

В терминах джойнингов свойство тензорной простоты ввели в рассмотрение
Д.Рудольф и А. дель Юнко   \cite{JR}.
Говорим, что мера $\nu$, заданная на $X^n$, принадлежит  классу
$M(m,n),\; n > m > 1$,
если проекции
$\nu$ на $m$-мерные грани декартова $n$-куба совпадают с мерой
$\mu^{\otimes m}$.
  Если для действия $\{T_g\}$ мера $\mu^{\otimes n}$ является
единственной мерой класса   $M(m,n),$ инвариантной относительно
$T_g\otimes\dots\otimes T_g$, говорим, что
действие принадлежит классу $ID(m,n)$.
Авторы \cite{JR}
ввели в рассмотрение класс $PID=\cap_{n>2}ID(2,n)$
 (их результат: класс $PID$ замкнут относительно
декартовых произведений).

В \cite{Kin} доказано, что
свойство $ID(2,4)$ влечет $ID(2,n)$ для всех $n>4$.
Этот результат стимулировал некоторые обобщения.
Как будет видно, $ID(2p-1,2p)=ID(2,4)$ для всех $p>1$, а
 класс тензорно простых действий совпадает с классом
$PID$.
Открыты вопросы о
совпадении классов $ID(2,3)$ и $ID(3,4)$ и          о
существовании слабо перемешивающего $\Z$-действия с нулевой энтропией,
не принадлежащего классу $ID(2,3)$.

Обозначение   $PID$
использовано Рудольфом и дель Юнко  \cite{JR} как аббревиатура от pairwise
independently determined action.
В диссертации отдается предпочтение операторной терминологии, в которой
словосочетание ``тензорная простота''
подчеркивает, что система и ее тензорная степень
не имеют нетривиальных марковских сплетений.
Придадим этому высказыванию формальный смысл.

{\bf Внутренние операторы.}
Пусть марковский  оператор $P$
(здесь его уместнее назвать бистохастическим)
действует из пространства $L^{\otimes m}$ в   $L^{\otimes n}$,
где $L=L_2(X,\mu)$, $m+n>2$.

Пусть   $$Im P\subset H^{\otimes n}\oplus \{Const\}, \ \
Im P^{\ast}\subset H^{\otimes m}\oplus \{Const\},$$
где $H$ -- пространство функций из $L$ с нулевым средним.
Такой бистохастический оператор назовем $внутренним$.

Система $(T,X,\mu)$ называется   тензорно простой, если для всех
$m,n$, $m+n>2$,
степени    $\T^{\otimes m}$ и   $\T^{\otimes n}$ имеют единственное
(тривиальное) внутреннее сплетение: $PL^{\otimes m}=\{Const\}$
(что равносильно  $PH^{\otimes m}=0$).
Оказывается, чтобы установить свойство тензорной простоты,
 достаточно проверить случай $m=n=2$ (т.е. свойство
S(3,4)). Доказательство этого факта см. в главе 1.

{\bf Джойнинги и сплетения.}
В тексте диссертации  сочетается язык и методология
сплетений и  джойнингов.
В ряде случаев  применение операторов весьма удобно.
В качестве подтверждающего
 примера  мы обсудим новое доказательство теоремы
Леманчика и дель Юнко \cite{JL}:
{\it
2-простая система  $T$ дизъюнктна
с гауссовской системой } (см. также статью Тувено \cite{T00}, где этот
результат обобщается).
Мы получим даже более
общий факт, но доказательство того, что он является более
общим, мы опустим.

Система  $(T, X,\mu)$ является  2-простой,
если любой эргодический джойнинг $\nu\neq \mu\otimes\mu$
двух копий  $T$  лежит на графике автоморфизма $S$, коммутирующего
с $T$.
Говорим, что $T$ дизъюнктен с $G$, если их единственным
марковским сплетением является $\Theta$ --
ортопроекция на пространство констант.

 Пусть система $(G,X)$ представлена $\prod_{i=1}^{n}(G_i,X_i)$,
будем обозначать
 $\breve P_k$ ортопроекцию на пространство $L_2(\prod_{i\neq k}X_i)$.
Говорим, что
система $(G, X)$ является  равномерно делимой, если
для любого  $\eps>0$  можно найти представление
$$(G,X)=\prod_{i=1}^{n}(G_i,X_i)$$
такое, что для всех $k=1,2,...,n$
операторы $\breve P_k$  $\eps$-близки тождественному оператору $I$.
(Считаем, что фиксирована метрика на марковском централизаторе
автоморфизма $G$.)

Пусть $A$ сплетает такой автоморфизм $G$ с  простым $T$.
Мы докажем, что  $A = \Theta$, следовательно,  $G$ и $T$ дизъюнктны.
(Этот результат верен и для квазипростых $T$.)

Известно (см. \cite{1}), что
для    неразложимых сплетений  $A$ и $B$
2-простой системы $T$ и эргодической системы $G$ выполнено
$B^{*}A = \Theta$
или   $BS= A$ для некоторого автоморфизма $S$, коммутирующего
с $T$.

Пусть $A\neq\Theta$, тогда  $A^*A\neq\Theta$.
Из  равномерной делимости $G$ получим
для некоторого $n$ такую факторизацию $G$, что  для всех
$ k=1,2,...n$  выполнено
$(\breve P_k A)^*A = A^{*}\breve P_k\,A\neq \Theta,$
 так как $\breve P_k$ близки к $I$.
Положим $B=\breve P_k\,A$. Поскольку $B$ неразложим, получим
$\breve P_k\,A\,S_k\ = A,$ $(k=1,2,...n).$
Тогда
$$A\ = \ \breve P_n\dots \breve P_1\, A\, S_1\dots S_n =
\Theta A\, S_1\,\dots S_n\
= \ \Theta.$$
Таким образом, $\Theta$ -- единственное сплетение $T$ и $G$.
Хорошо известно, что гауссовская система
является  делимой (см. \cite{JL}), можно показать чуть больше:
она является равномерно делимой.

\newpage

\begin{center}
\bf    0.4.  Структура и основные результаты диссертации
\rm
\end{center}

\bf Краткий обзор полученных результатов. \rm
Значительная часть
 диссертации
посвящена следующей общей задаче: {\it пусть
джойнинг набора $n\geq 3$ копий динамической системы
обладает свойством попарной независимости, верно ли,
что этот джойнинг является произведением мер?   }
В нашей терминологии этот вопрос переформулируется так:
{\it является ли система тензорно простой?}

Исторически этот вопрос возник для систем с минимальным (в более
общем случае -- простым)  централизатором.
Как отмечалось, положительный ответ был получен для перемешивающих
систем  ранга 1 ( Каликов, Кинг, Рыжиков) и систем с сингулярным спектром
(Ост).

В главе 1  излагаются методы теории марковских сплетений и показано
как они применяются для установления тензорной простоты
 и других свойств систем.

В главах 2,3 тензорная простота устанавливается  для систем с минимальным
или простым марковским централизатором  при некоторых дополнительных
условиях. Показано, что для $\Z$-действий вопрос
$MSJ(2)=MSJ(3)$? равносилен проблеме Рохлина в классе  $MSJ(2)$ (классе
систем с минимальными самоприсоединениями порядка 2).
Для потоков доказано совпадение свойства простоты порядка 2 и свойства
простоты всех  порядков
(в частности, этим доказано для потоков совпадение классов
$MSJ(2)$ и $MSJ(3)$),
что дает положительный ответ для потоков на
вопрос Рудольфа и дель Юнко.
Получены обощения этих результатов.
Для некоторых некоммутативных действий обнаружено различие
четной и нечетной тензорной простоты и показано, что
$MSJ(2)\neq MSJ(3)\neq MSJ(4)$.
Решена проблема Рохлина для потоков положительного локального
ранга.

В главе 4 тензорная простота установлена
для  перемешивающих действий конечного
ранга и $\Z^n$-действий положительного локального ранга $\beta > 2^{-n}$,
тем самым установлено и свойство кратного перемешивания.
Как следствие получено равенство классов MSJ(2) и MSJ(3) для действий
конечного ранга.
Установлена бесконечность ранга эргодического автоморфизма
$T\times T$ и
точная оценка локального ранга $T\times T$. В случае, когда он равен
максимальному значению $\frac{1}{4}$, показано,
что $T$ обладает свойством $ \kappa$-перемешивания.
Так как известны перекладывания $T$ такие, что локальный ранг
$T\times T$ равен $\frac{1}{4}$ (примеры А. Катка),
подтверждена гипотеза
Оселедца \cite{Os} о существовании перекладываний со свойством
 $ \kappa$-перемешивания ( $ \kappa =\frac{1}{2})$.

Глава 5 содержит следующие результаты. Дано положительное решение
проблемы Рохлина о непростом однородном спектре  для неперемешивающих и
перемешивающих автоморфизмов.  Для  $ \kappa$-перемешивающих автоморфизмов
$T$ доказано, что изоморфизм $T\times T$ и $S\times S$ влечет за собой
изоморфизм $T$ и $S$.
Рассмотрен асимптотический инвариант (частичное  кратное возвращение),
который может различать некоторые автоморфизмы с их
обратными. Изучен новый класс расширений типа $(T, T^{-1})$-расширений,
 сохраняющих
 свойства тензорной простоты  и  кратного перемешивания.

\newpage
\begin{center}
{\bf ОБЗОР РЕЗУЛЬТАТОВ ДИССЕРТАЦИИ.}
\end{center}
\begin{center}
\bf Глава 1.
   Марковские сплетающие операторы и тензорная простота.
\end{center}

\bf  1.1.     Несколько методологических принципов теории сплетений.
\rm
Этот параграф начинается с  описания естественных
полугрупп операторов, в терминах которых выражаются многие известные
свойства  динамической системы.
Описаны  алгебраические операции над сплетающими операторами:
 композиции, сопряжение, изменение типа сплетения.
Далее излагается  ряд приемов, которые наиболее часто
используются в доказательствах. К ним относятся
сплетения слабых замыканий,
разложимость и  неразложимость сплетений,  принцип симметризации,
дополнительная симметрия
и принцип индуцированных джойнингов (последние два принципа вынесены
в отдельные параграфы).
Чтобы показать преемственность некоторых методов, мы приводим доказательства
теоремы Блюма и Хансона \cite{BH},
 эргодической теоремы  фон Неймана.

\bf  1.2.  Дополнительная симметрия.  \rm
В  1.2. излагаются некоторые приложения дополнительной симметрии
и, в частности, дано новое доказательство теоремы Фюрстенберга
о кратном перемешивании в среднем на прогрессиях.

\bf  1.3.  Индуцированные джойнинги.  \rm
В простейшей модельной ситуации показано, как работают индуцированные
джойнинги (основные  приложения даны в главах 2 и 4).
Индуцированные джойнинги определяются следующим образом.
Пусть для перемешивающих автоморфизмов $R,S,T$  выполнено
тождество (условие    эквивариантности):
$$
  \cal{J}\it (Rx) \equiv T^{\ast }\cal{J}\it (x)S,
$$
где $\cal{J}\it : X \to{} \cal{M}$ -- семейство марковских операторов,
отвечающее некоторому джойнингу $\nu =(R\times S\times T)\nu$.

Для  семейства операторнозначных функций
$ {\cal H}_m :X \to{} {\cal M} $:
$$
  {\cal H}_m(x)  \equiv \J^\ast\it (x)\J (R^mx)\equiv
  \J^\ast\it (x)T^{\ast m} \J\it (x)S^{m}.
$$
будет  выполнено тождество
$$
 {\cal H}_m(Rx) \equiv S^\ast {\cal H}\it _m(x)S .
$$
Функции $\cal H_m$ сопоставим меру $\eta_m$ :
$$
\eta_m(A\times B\times C)= \int_X \chi_A(x)\langle
 H_m(x)\chi_B,\chi_C\rangle  d\mu(x).
$$
Таким образом, мы определили новые джойнинги $\eta_m$ (отметим, что
$\eta_m = (R\times S\times S)\eta_m$), которые будем называть
индуцированными.

Тривиализация индуцированных джойнингов $\eta_m$ влечет за собой
тривиальность исходного джойнинга $\nu$, что дает решение задачи
о тензорной
простоте и проблемы Рохлина для рассматриваемого
класса систем.   В доказательствах индуцированные джойнинги играют
 вспомогательную  роль:
 в конечном итоге  они оказываются тривиальными. Иначе
говоря, с помощью нетривиальных индуцированных джойнингов мы доказываем,
 что в рассматриваемых ситуациях они не существуют.

Содержание  1.1, 1.2, 1.3 опубликовано в работах
автора  \cite{preprint}, \cite{Izv},
 \cite{MZ92}, \cite{1}, \cite{MS}.
\medskip

\bf  1.4.     Примеры тензорно простых систем.   \rm
В этом параграфе (следуя (\cite{MZ92},\cite{F3})  доказывается,
что тензорной простотой обладают
$ \kappa$-перемешивающий автоморфизм  и автоморфизм $T$,
содержащий в слабом замыкании  степеней оператор $\frac{1}{2}(I+T)$;
неперемешивающие автоморфизмы с минимальным
централизатором;   простые действия,
коммутирующие со слабо перемешивающим неперемешивающим автоморфизмом.

\bf  1.5.     Связь типов тензорной простоты.   \rm
Как обобщение теоремы Кинга, доказавшего в терминах джойнингов,
что свойство $S(2,4)$ влечет за собой свойства $\s_{k}$, $k=3,4,...$,
($\s_{k}$ -- краткое обозначение свойства
$S(k-1,k)$), доказана

{\bf Теорема 1.5.1.}(\cite{MS}) \it  Для любых $q\geq 2$ и $k\geq 3$ свойство
  $\s_{2q}$ влечет свойство $\s_{k}$
\rm

Рассмотрен также случай, когда
 свойство $\s_{2n+1}$ влечет за собой свойство
$\s_3$.

\rm

\begin{center} \bf
Глава 2.  Минимальные самоприсоединения, простота и квазипростота групповых
действий.
\end{center}
\bf  2.1.     Простые  системы с несчетным централизатором.
\rm
 Действие $T$ называется $n$-простым ( или простым порядка $n$)
если любой эргодический джойнинг
$\nu\neq \mu^{\otimes n}$ набора из  $n$ копий действия $T$
обладает свойством:
одна из его проекций  на двумерную грань  в
$X\times\dots\times X$ являеся мерой  $\Delta_S = (Id\times S)\Delta$
 для
некоторого автоморфизма $S$, коммутирующего с действием $T$.
Основной результат этого параграфа следующий.
\medskip

{\bf ТЕОРЕМА 2.1.2.}.   \it
Слабо перемешивающий 2-простой поток  является  простым всех порядков.
Перемешивающий 2-простой поток перемешивает с любой кратностью.
\rm
\medskip

{\bf  2.2.     Наследственная независимость и квазипростота действий. }
Джойнинг $\nu$ пары $(T,T)$  называется квазидиагональной мерой, если
для почти всех $x,y$ условные меры  $\nu_x$ и $\nu_y$, возникающие
в представлении
$$
\nu(A\times B)=\int_A \nu_x(B)d\mu(x) = \int_B \nu_y(A)d\mu(y),
$$
имеют вид:
$$
    \nu_x = \frac{1}{p} (\delta_{y_1(x)}+
\delta_{y_2(x)}+\dots + \delta_{y_p(x)}),
$$
$$
    \nu_y = \frac{1}{q} (\delta_{x_1(y)}+
\delta_{x_2(x)}+\dots + \delta_{x_q(y)}).
$$
Говорят, что
действие $T$ квазипростое порядка $n$, если
для любого эргодического джойнинга
 $\nu\neq \mu^{\otimes n}$
набора из  $n$ копий действия $T$ выполнено:
одна из  проекций  на двумерную грань  в
$X\times\dots\times X$ являеся  квазидиагональной мерой.
\rm
Приведенные выше определения автоматически распространяются на произвольные
групповые  действия.

{\bf Теорема 2.2.1.}(\cite{2}) \it
Если слабо перемешивающий автоморфизм обладает свойством
квазипростоты порядка 3, то он является
квазипростым всех порядков.  \rm

Эта теорема обобщает результат Глазнера, Оста и Рудольфа \cite{GHR}
о том, что для $\Z$-действий 3-простота влечет простоту всех порядков.
         \medskip

\bf  2.3.     Минимальные самоприсоединения и  кратная возвращаемость.
\rm
Автоморфизм $T$ обладает минимальными самоприсоединениями порядка $n$,
если он коммутирует только со своими степенями и является простым
порядка $n$.

{\bf Теорема 2.3.1.}(\cite{2}) \it
Если $T\in MSJ(2)$ и $T$ перемешивает с кратностью 2, то
автоморфизм $T$ обладает минимальными самоприсоединениями всех порядков и
кратным перемешиванием всех порядков.
\rm

Эта теорема является следствием более общего утверждения,
в котором фигурирует свойство более слабое, чем  перемешивание
 кратности 2.

{\bf Теорема  2.3.2.} (\cite{2}) \it   Пусть  перемешивающий
автоморфизм $T\in MSJ(2)$ обладает свойством
кратного возвращения:
  для любого множества $A$ положительной меры  и любых
 последовательностей  $k(m)\to \infty$
$|k(m)-m|\to \infty$ таких, что  для всех больших $m$ выполнено
условие
$$
 \mu(T^{-k(m)}A\cap T^{-m}A\cap A) > c > 0.
$$
Тогда  обладает свойством кратного
перемешивания и свойством минимальных самоприсоединений всех порядков.
\rm
\medskip

\bf  2.4.     Четная и нечетная тензорная простота.
\rm

{\bf Теорема 2.4.1.}(\cite{MZ96}) этого параграфа дает
пример  некоммутативного   действия   $\Psi$,
 которое
 обладает свойствами $S(2q, 2q+1)$, но не обладает свойствами
$S(2n-1, 2n)$,  что приводит к понятию нечетной тензорной простоты.
  \rm

Следующий результат показывает, что для некоммутативных систем
         $MSJ(2)\neq MSJ(3)\neq MSJ(4)$.
         \medskip

{\bf Теорема 2.4.2.}(\cite{MZ96}) \it 1.
Действие $\Psi$,  порожденное всеми  автоморфизмами
группы $Y= \Z_3\times\Z_3\times\Z_3\dots$
 и всеми групповыми сдвигами на $Y$, принадлежит   классу  $MSJ(2)$,
но не принадлежит классу $MSJ(3)$.
Действие $\Phi$,  порожденное всеми автоморфизмами
группы $X= \Z_2\times\Z_2\times\Z_2\dots$
 и всеми сдвигами на группе $X$, принадлежит   классу
$MSJ(3)\setminus
MSJ(4)$.
Эргодическими джойнингами как  $(\Phi,\Phi)$, так и $(\Psi,\Psi)$
являются  только меры $\Delta$ и $\mu\otimes \mu$.

2. Классу $MSJ(3)\setminus MSJ(4)$ принадлежит действие  $\Phi'$,
порожденное (бернуллиевским) автоморфизмом
$T$
и инволюциями $Q,R,S$, определенными
на    $\dots\times\Z_2\times\Z_2\times\Z_2\dots$
следующим образом:
$$T(x)_i=x_{i+1},$$
$$Q(\dots x_{-2}, x_{-1}, x_{0}, x_{1},  x_{2},\dots)=
(\dots x_{-2}, x_{-1}, x_{1}, x_{0},  x_{2},\dots),$$
$$R(\dots x_{-2}, x_{-1}, x_{0}, x_{1},  x_{2},\dots)=
(\dots x_{-2}, x_{-1}, x_{0}, x_{1}+ x_{0},  x_{2},\dots),$$
$$S(\dots x_{-2}, x_{-1}, x_{0}, x_{1},  x_{2},\dots)=
(\dots x_{-2}, x_{-1},  x_{0}+1,x_{1}, x_{2},\dots).$$
\medskip
                  \rm

Пункт 2 этой теоремы показывает, что соответствующие примеры
имеются среди действий конечно порожденных групп.
         \medskip
\newpage
\begin{center}    \bf
Глава 3.   Джойнинги и тензорная простота некоторых потоков.
\end{center}
{\bf  3.1. Гладкие джойнинги и  внутренние сплетения потоков. }

М. Ратнер доказала, что любой эргодический
джойнинг  унипотентного потока являются гладким: он сосредоточен
на гладком подмногообразии $Y\subset X^n$ (декартова степень $X$) и абсолютно
непрерывен относительно меры Лебега на  $Y$.
Назовем потоки с таким свойством S-потоками.
\medskip

{\bf Теорема 3.1.1.}(\cite{UMN95},\cite{F3}) \it Слабо перемешивающий S-поток
является тензорно простым. \rm
         \medskip

{\bf  3.2. Кратное перемешивание  $\omega$-простых   потоков.   }
Поток называется $\omega$-простым, если нетривиальный
эргодический джойнинг
второго порядка  сосредоточен
на графике конечнозначного    отображения.

{\bf Теорема 3.2.1.}(\cite{MS}) \it
 Слабо перемешивающий $\omega$-простой поток
является тензорно простым.

{\bf Следствие.} \it  Перемешивающий $\omega$-простой поток является
перемешивающим всех порядков. \rm

\vspace{3mm}
{\bf  3.3. Тензорная простота  потоков положительного локального  ранга.}
Понятие ранга автоморфизма тесно связано со свойством
циклической аппроксимации,
изучавшейся в А.Б. Катком, В.И. Оселедцем, А.М. Степиным \cite{St67},
\cite{KS},   и связано с конструкциями Р. Чакона \cite{Chacon0} и
Д. Орнстейна \cite{Ornstein} (см. также \cite{Abdalaoui1}).
Cистемы локального положительного ранга  являются естественным
подклассом класса стандартных систем, введенных А.Б. Катком и
Е.А. Сатаевым \cite{KSa}.
Интерес к системам  конечного и положительного локального ранга
стимулировался  изучением спектральной кратности эргодических
автоморфизмов (см. \cite{St73},
\cite{Chacon1}).
Новый интерес был связан с изучением  джойнингов  этих систем
(см., работу Кинга \cite{Kin}).

Автоморфизм $S$ пространства Лебега $(X,\mu)$, $\mu(X)=1$,
обладает \it локальным рангом \rm
$\beta(S)$, если  $\beta(S)$ есть максимум чисел $\beta\geq 0$ таких,
что  для некоторой последовательности конечных разбиений
пространства $X$   вида
$$ \xi_j=\{ B_j,  SB_j,    S^2B_j,    \dots,  S^{h_j-1}B_j,\dots\}$$
выполнено:
любое фиксированное измеримое множество аппроксимируется
$\xi_j$-измеримыми множествами при $j\to\infty$ (пишем $\xi_j\to\eps$),
причем $\mu(U_j)\to \beta$, где
$$U_j=\bigsqcup_{0\leq k<h_j}S^kB_j.$$

Аналогично локальный ранг определяется для потоков.
\medskip

{\bf Теорема 3.3.1.}(\cite{Izv},\cite{F4})
\it Перемешивающий поток $\{T_r\},$ $r\in\R^n$, $n \geq 1$,
при $\beta(\{T_r\}) > 0$ обладает свойством перемешивания
всех порядков. \rm
\medskip

{\bf Теорема 3.3.3.} (\cite{Izv},\cite{F4})\it

Эргодический джойнинг $\nu\neq\mu\otimes \mu$ двух
копий перемешивающего потока $\{T_r\}$ $r\in\R^n$, $n \geq 1,$
при $\beta(\{T_r\}) > 0$ является мерой, сосредоточенной на  графике
конечнозначного отображения.

Для набора таких  джойнингов  $\nu_1,\nu_2,\dots,\nu_p$ при
\\ $p\beta(\{T_r\}) > 1$ для некоторых
$i\neq k$   выполнено   $\nu_i=(I\times T_v)\nu_k$.
\rm
\medskip

\begin{center} \bf
Глава 4.   Джойнинги и кратное перемешивание действий конечного
           и положительного локального ранга.
\end{center}      \rm
\medskip

Говорят, что автоморфизм $S$ имеет {\it ранг}  $r=Rank(S)$,
если $r$ есть минимальное число такое, что найдется
последовательность разбиений $\xi_j\to\eps$  вида
$$ \xi_j=\{  B^1_j,\,    SB^1_j,\,
 \dots
    S^{h_j^1}B^1_j,\ \ \dots\ \
  B^r_j,\,  SB^r_j,\,      \dots ,
  S^{h_j^r}B^r_j,\   Y_j  \}$$
(условие $\xi_j\to\eps$
влечет за собой   ${h_j^1}$, $\dots$, ${h_j^r}\to \infty$
и  $\mu(Y_j)\to 0$).

\medskip

{ \bf  4.1.
    D-свойство перемешивающих автоморфизмов  конечного  ранга.  }

Последовательность разбиений множеств $U_j\subset X$
вида
$$ \xi_j=\{ E_j, TE_j,\dots T^{h_j}E_j\}$$
назовем   аппроксимирующей,
если, дополняя разбиение $ \xi_j$ некоторым разбиением  дополнения
$X\setminus U_j$,
получим последовательность разбиений всего фазового пространства $X$,
которая стремится к разбиению на точки.

Будем говорить, что автоморфизм $T$ обладает D-свойством, если
найдутся последовательности аппроксимирующих башен
$$ (U_j,\xi_j),\ \  (U'_j,\xi'_j), \ \ (U''_j,\xi''_j),   $$
где
$$ \xi_j=\{ E_j, TE_j,\dots T^{h_j}E_j\},$$
$$ \xi'_j=\{ E'_j, TE'_j,\dots T^{h_j}E'_j\},$$
$$ \xi''_j=\{ E''_j, TE''_j,\dots T^{h_j}E''_j\},$$
причем  для некоторой последовательности $m_j$, $m_j>h_j$,
выполняются следующие условия:
$$\lim_j \mu(U_j) =a>0, \ \  \mu(E_j)=\mu(E'_j)=\mu(E''_j), $$
$$ E'_j= T^{m_j}E_j,\ \  \mu(T^{m_j}U'_j\Delta U''_j)\to 0, $$
$$ \max_{m>h_j}\mu(T^{m}E'_j\ |\ E''_j)\ \to\ 0.$$
\medskip

{\bf Теорема 4.1.1. } (\cite{F1})\it
 Перемешивающий автоморфизм конечного ранга
обладает D-свойством.
         \rm
\medskip

{ \bf  4.2.      D-свойство  перемешивающих $\Z^n$-действий и локальный ранг.
}

{\bf Теорема 4.2.1.}(\cite{F4})
\it   Если $\Z^n$-действие $\{T_z\}$ обладает свойством перемешивания и
$\beta\{T_z\}>\frac{1}{2^{n}}$, то действие обладает  \it D-свойством. \rm
\medskip

{ \bf  4.3.      Тензорная простота перемешивающих систем с D-свойством.
}

{\bf Теорема 4.3.1.}(\cite{F1}) \it Перемешивающее $Z^n$-действие, обладающее
 D-свойством, является тензорно простым и, следовательно, обладает
 перемешиванием всех порядков.  \rm

{\bf Следствие.} (\cite{UMN91},\cite{F1}) \it Перемешивающие автоморфизмы
конечного ранга обладают перемешиванием всех порядков. \rm
\medskip

{ \bf  4.4.      Локальный ранг и кратное перемешивание.
}
Пусть $T$ -- перемешивающий автоморфизм, предположим, что
он не обладает перемешиванием кратности 2. $A,B,C$ -- некоторые измеримые
множества.
Определим $\e$-отклонение от кратного перемешивания:
$$
  Der(\varepsilon ,A,B,C) = $$ $$
\lbrace (z,w)\in Q(\varepsilon ,h): | \mu (A\bigcap T^z B\bigcap
T^w C) - \mu (A)\mu (B)\mu (C)| > \varepsilon \rbrace,
$$
где
$$
  Q(\varepsilon ,h) = \lbrace (z,w) \epsilon [0,h] :
|z|,|w|,|z-w| > \varepsilon h \rbrace.
$$
Положим
$d(h)=\sharp Der(\varepsilon ,A,B,C)/ h$.
(Аналогичным образом  $d(h)$ определяется для действий групп $\Z^n$.)

Имеет место следующий факт (\cite{R0}, пар. 2):  если  $T$ --
перемешивающий автоморфизм, то $d(h)$ -- ограниченная последовательность.

Если же действие $\lbrace T^z \rbrace $ перемешивает двукратно, то,
очевидно,  $d(h)=0$ начиная с некоторого $h$.

Таким образом, можно рассмотреть свойство, промежуточное между однократным
и двукратным перемешиванием.
Действие $\lbrace T^z : z \epsilon Z^n\rbrace $ назовем
 $(1+\varepsilon)$- перемешивающим, если $$d(h) \rightarrow \, 0$$
для любых измеримых множеств $A,B,C$ и $\varepsilon >0.$
     \medskip

{\bf Теорема 4.4.1.}(\cite{F2}) \it $(1+\varepsilon)$-перемешивающее
$Z^n$-действие $\Psi $
положительного локального ранга является тензорно простым.  \rm

{\bf Следствие.} (\cite{F1}) \it  Перемешивающий с кратностью 2
автоморфизм положительного локального ранга является тензорно простым.
                                \rm \medskip

{\bf  4.5.      Ранги и джойнинги $T\times T$.}
Примеры автоморфизмов $T$
с  положительным локальным рангом
были предъявлены Катком в связи
с изучением спектральной кратности автоморфизмов пространства Лебега.
     \medskip

{\bf Теорема 4.5.1. }(\cite{F5}) \it   a). Ранг декартова квадрата
автоморфизма равен бесконечности: $Rank(T\times T)=\infty$.\

b). Локальный ранг $\beta(T\times T)$ не превосходит $\frac{1}{4}$.
\rm
     \medskip

{\bf Теорема 4.5.4. }(\cite{F5}) \it
Если  локальный ранг эргодического автоморфизма $T\times T$ равен
$\frac{1}{4}$, то $T$ обладает \\
$ \kappa$-перемешиванием при $ \kappa=\frac{1}{2}$:
некоторая последовательность степеней $T$ слабо сходится к
оператору $ \frac{1}{2}I+\frac{1}{2}\Theta$
($\Theta$ -- ортопрекция   на пространства констант). \rm
     \medskip

Последний результат вместе с примером А. Катка
 перекладывания $T$ трех отрезков со свойством
$\beta(T\times T)=\frac{1}{4}$
подтверждает давнюю гипотезу Оселедца  \cite{Os}
 о существовании
перекладываний со свойством  $ \kappa$-перемешивания ($ \kappa=\frac{1}{2}$).
\begin{center}   \bf
Глава 5.   Некоторые спектральные, алгебраические и асимптотические свойства
           динамических систем.
\end{center}         \rm

{\bf  5.1.  Проблема Рохлина об однородном непростом спектре.
}
В спектральной теории динамических систем давно стоял вопрос о существовании
эргодичекого автоморфизма с однородным непростым спектром (см. брошюру
Д.В.Аносова \cite{An3}, посвященную спектральной кратности автоморфизмов).
А. Каток высказал гипотезу о
том, что  нужные системы  можно найти  среди декартовых квадратов
типичных автоморфизмов.
В этом параграфе    предлагается два класса
 автоморфизмов  $T$, декартов квадрат  которых
имеет однородный спектр кратности 2.
Для рассматриваемых  автоморфизмов $T$ устанавлено  свойство
$\sigma\ast\sigma\perp \sigma$, где
$\sigma$ -- спектральная мера автоморфизма  $T$.
     \medskip

{\bf Теорема 5.1.1.}(\cite{3},\cite{SRM}) \it
Пусть $T$ -- эргодический автоморфизм,  и для
некоторой последовательности ${k_i}\to \infty$ и числа $a\in (0,1)$
выполнено
$\ \ {\widehat{T}}^{k_i}\to  (aI+(1-a){\widehat{T}})$.
Тогда

1) для спектральной меры $\sigma$ автоморфизма $T$
  выполнено $\sigma\ast\sigma\perp\sigma$;

2) если $T$ имеет простой спектр, то
$(T\times T)$ имеет однородный спектр кратности 2.
\vspace{5mm}
     \medskip    \rm

Теорему 5.1.1. независимо доказал и одновременно с автором
 опубликовал О.Н. Агеев \cite{Age}.
Автоморфизмы со свойством
$\ \ {\widehat{T}}^{k_i}\to 0.5 (I+{\widehat{T}})$
рассматривались в статье А.Б.Катка и А.М.Степина \cite{KS1}.

Следующая теорема содержит аналогичный результат, но теперь для класса
автоморфизмов, обладающих лучшими перемешивающими свойствами, когда
параметр $a$ близок к 1.
     \medskip

\rm
{\bf Теорема 5.1.2.} (\cite{3},\cite{SRM})\it
Пусть для эргодического автоморфизма $T$
 для некоторой
последовательности ${k_i}\to \infty$ и числа  $a\in (0,1)$  выполнено
$ {\widehat{T}}^{k_i}\to \
{(1-a)}(I+a\widehat{T} +a^2\widehat{T}^2+\dots).$
Тогда

1) для спектральной меры $\sigma$ автоморфизма $T$ выполнено
$\sigma\ast\sigma\perp\sigma$;

2) если $T$ имеет простой спектр, то
автоморфизм $R$, $R(x,y)=(Ty,x)$, имеет простой спектр, а автоморфизм
$(T\times T)$ имеет однородный  спектр кратности 2.
\vspace{5mm}  \rm
     \medskip

{\bf \S 5.2.  Перемешивающие автоморфизмы с однородным непростым спектром.
}
Оказывается, что, улучшая перемешивающие свойства автоморфизма $T$
со свойством
" $(T\times T)$ имеет однородный  спектр кратности 2 ", мы можем
в пределе сохранить это свойство и получить перемешивающий автоморфизм.
Таким образом, проблема Рохлина получает решение в классе перемешивающих
систем.

\medskip
{\bf Теорема 5.2.1.} (\cite{SRM})
\it Существует перемешивающий автоморфизм $T$ такой,
что симметрический квадрат ${T\odot T}$ имеет простой спектр.
\medskip

{\bf Следствие.} Соответствующий перемешивающий
автоморфизм  $T\times T$
имеет однородный спектр   кратности 2.   \rm
\medskip

{\bf  5.3.  Изоморфизм  декартовых степеней преобразований
и $ \kappa$-перемешивание. }
Для типичных автоморфизмов $T$  мы отвечаем положительно на вопрос
 Тувено: " влечет ли изоморфизм  $T\times T$ и $S\times S$
за собой изоморфизм $T$ и $S$?"\    При этом мы пользуемся
результатом Степина \cite{St86} о типичности $ \kappa$-перемешивания
и следующим результатом.
\medskip

 \bf Теорема 5.3.1 (\cite{MZ96a})\it Если автоморфизм $T$ обладает
свойством     $ \kappa$-перемешивания, $0< \kappa<1$,
то изоморфизм $T\times T$ и $S\times S$
 влечет за собой изоморфизм $T$ и $S$.
\rm
\medskip

{\bf  5.4.  Асимметрия прошлого и будущего динамической системы.
}
\medskip

Примеры преобразований, не изоморфных своему обратному,
известны давно. Цель этого параграфа -- предложить новый инвариант,
который может различить автоморфизм $T$ и $T^{-1}$.   Таковым
является свойство кратной возвращаемости на последовательностях.
\medskip

{\bf Теорема 5.4.1.}(\cite{TM})\it
Существует автоморфизм $T$, обладающий
свойством:  для некоторой последовательности
\\ $n(i)\to\infty$ для любого множества  $A\in\B$  выполнено
$$
 \lim_{i\to \infty} \mu(A\cap T^{n(i)}A\cap T^{3n(i)}A) \geq
\frac{1}{5} \mu(A),
$$
при этом для некоторого множества $A', \ \mu(A') > 0,$ имеет место
$$
\lim_{i\to \infty} \mu(A'\cap T^{-n(i)}A'\cap T^{-3n(i)}A') =  0 .
$$
В качестве такого множества $A'$ годится
любое множество, удовлетворяющее условию
$$\mu(A'\cap TA')=\mu(A'\cap T^2A')=0.$$  \rm

{\bf Следствие.} \it   Автоморфизм $T$ асимметричен. \rm
\medskip

{\bf  5.5.  Расширения, сохраняющие
тензорную простоту и кратное перемешивание.
}
\medskip

Следующая теорема является обобщением   того факта, что тензорно
простая перемешивающая  система обладает перемешиванием любой кратности.
\medskip

{\bf Теорема 5.5.1.}(\cite{UMN94},\cite{F3})
\it  Пусть автоморфизм $S$ перемешивает с кратностью $k$,
а косое произведение $R$,
$R(x,y) = (S(x),T_x(y))$, является перемешивающим. Если
все преобразования  $T_x$ коммутируют с некоторым тензорно простым
действием $\Psi$ (быть может, некоммутативным), то косое произведение
$R$    перемешивает с кратностью $k$.
\rm
\medskip

Тип расширений, который приведен ниже, сохраняет свойство тензорной простоты
даже в случае, когда автоморфизм
$T$ этим свойством не обладает. В случае
бернуллиевского $T$ косое произведение $R$ обладает континуальной системой
факторов.
\medskip

{\bf Теорема 5.5.2.}(\cite{F3}) \it Пусть $R,T$ -- перемешивающие
преобразования,
где $R$ есть косое произведение над $S$ следующего вида:
$$
  R(x,y)=(S(x),T^{n(x)}(y)), \quad \int n(x)d\mu = 0.
$$
Если автоморфизм $S$ перемешивает с кратностью $k$, то  косое произведение
$R$ также обладает перемешиванием  кратности $k$.
Если автоморфизм $S$ является тензорно простым,
то  $R$ также является тензорно простым.
\rm
\medskip

Доказательство  теоремы 5.5.2 использует хорошо известную
теорему Аткинсона \cite{Atk} о возвратности случайных блужданий,
которую можно назвать теоремой  Крыгина-Аткинсона, так как
 она тривиально вытекает из результата А.Б.Крыгина \cite{Kry}
о консервативности   цилиндрического каскада,
ассоциированного с косым произведением
 $R$ при $\int n(x)d\mu = 0$.
\medskip
\medskip

\newpage

\begin{center}   { \Large  \bf  ГЛАВА 1 }
\end{center}
\begin{center}
{ \bf   МАРКОВСКИЕ СПЛЕТАЮЩИЕ ОПЕРАТОРЫ И ТЕНЗОРНАЯ ПРОСТОТА ДИНАМИЧЕСКИХ
СИСТЕМ   }
\end{center}
\medskip
\medskip

Эта глава в основном посвящена методологии теории марковских
сплетений. Основная цель -- показать приемы, использующиеся
для доказательства
тензорной простоты и кратного перемешивания динамических систем.
Приведены результаты, касающиеся связей свойств тензорной простоты
разных порядков.

Приведем для удобства читателя некоторые из обозначений,
используемых в работе.
\\$A\,,B\,,C$ (и те же символы с индексами) -- измеримые подмножества
фазового пространства динамической системы,
а иногда ( в главе 1) -- индикаторы этих множеств.
\\$Id$ -- тождественное преобразование (обычно множества $X$).
\\$\langle\cdot,\cdot\rangle$ -- скалярное произведение в $L_2$.
\\$I$ -- тождественный оператор (обычно на $L_2(X,\mu)$).
\\$\Theta$ --  оператор ортопроекции
на пространство констант:
  $$\Theta f = Const\equiv \int f(x)d\mu(x).$$
\\ Иногда $\Theta$ также обозначает бистохастический
  оператор с образом, являюшимся одномерным пространством констант.
В рамках нашей работы бистохастический  оператор отличается
от марковского лишь тем, что первый действует из одного пространства
в другое, а второй -- из пространства в себя.
\\$R,S,T$ (и те же символы с индексами) обозначают автоморфизмы
  пространства Лебега или, когда это видно из контекста, унитарные операторы,
  отвечающие этим автоморфизмам.
\\$J\,,P\,,Q$ -- марковские (или бистохастические) операторы
(обычно они  сплетают автоморфизмы).
\\$\eta,\,\nu, \lambda$ -- меры, которые являются джойнингами
относительно  тензорных произведений динамических систем. Эти меры
заданы на декартовых произведениях фазового пространства изучаемой
(исходной) динамической системы.     

\begin{center}
{\bf  1.1. Несколько методологических принципов теории сплетений}
\end{center}
\medskip

{\bf Полугруппы положительных  операторов}.
Определим серию полугрупп
 положительных ограниченных  операторов в  пространстве $L_2({X},\mu)$.
Следует сказать,  что операторы
 называются положительными в том смысле, что они переводят
неотрицательные функции   в неотрицательные.

$\cal{P}$ -- полугруппа ограниченных положительных  операторов
на $L_2({X},\mu)$.
Оператор ${P}$ называется  положительным если, для любой
$f\in L_2({X},\mu),  f(x)\geq 0 $ выполнено
$$
  f(x)\geq 0 \Rightarrow (Pf)(x)\geq 0 .
$$

$\cal{R}$ --  полугруппа регулярных операторов.
Оператор ${R}$ называется  регулярным, если он является интегральным
оператором с ограниченным
$\mu\otimes\mu$-измеримым ядром $K(x,y)$:
$$
 Rf(y)= \int_{X}   K(x,y) f(x) d\mu (x) .
$$

$\cal{M}$  -- полугруппа марковских   операторов, т.е. все
операторы из $\cal{P}$, которые удовлетворяют условию
$$
                           P1 =P^{\ast}1=1.
$$

$\cal{D}$  -- полугруппа всех  операторов из  $\cal{M}$, удовлетворяющих
условию: для некоторого $a > 0$  выполнено
$$
                J^\ast J \geq aI \leq JJ^\ast.
$$

$ \cal{A}$ -- максимальная группа в  $\cal{M}$.
Эта группу составляют  операторы $T\in\cal{M}$ такие , что
$$
                T^\ast T = I = TT^\ast .
$$
Операторы $T\in\cal{A}$ и только они отвечают автоморфизмам
$\check T$ пространства Лебега $({X}, \cal{B},\mu)$:
$
      Tf(x) =f(\check Tx).
$

{\bf Естественные структуры и топология на полугруппах. }
На полугруппе $\cal{P}$ можно задать порядок $\leq$:
$$
 P\leq Q,$$ если  для всех $f\geq 0, f\in L_{\infty}$
 выполнено (mod0) неравенство
 $(Pf)(x)\leq (Qf)(x).$

На полугруппе $\cal{M}$  определим бинарное отношение
$\perp$: пишем $J_1 \perp J_2$, если
для любого оператора $P\in {\cal P}, \ P\neq 0,$ условие   $P< J_1$
несовместимо с условием $P< J_2$ (соответствующие полиморфизмы
взаимно сингулярны как меры).

Слабая операторная топология на $\cal{M}$ превращает эту полугруппу
в компакт. Любая последовательность марковских операторов имеет
предельную точку.  Слабая сходимость $ P_j\,\to \ {P}$ означает
$$\forall f,g\in L_2 \ \ \ \
  \left< P_j f\ |\ g\right>\ \to \ \left< Pf|g\right>.
$$

{\bf Эргодичность, слабое перемешивание, перемешивание.}
Оператор $T\in \cal{A}$ эргодический,  если   выполнено
$$T f=f \,\Rightarrow \, f=\Theta  f, $$
где $\Theta$ -- оператор ортопроекции на пространство констант.
Известно, что эргодичность эквивалентна следующему статистическому
свойству (перемешиванию в среднем):
$$
P_N := \frac{1}{N} \sum_{i=1}^N T^i \to \Theta,
$$
что легко доказать, пользуясь компактность полугруппы $\M$.
Если последовательнось $P_N$ имеет предельную точку
$P$, то $T P=P$, следовательно, для любой функции $f$ выполнено
$$T(Pf)=Pf=const=\Theta f,$$
так как
$T$ -- эргодический оператор.
\medskip

Оператор $T\in \M$ называется перемешивающим, если
$
                   T^{n} \to  \Theta  $
при $ n \to\infty$.
Для автоморфизмов это определение эквивалентно  следующему:
для любых измеримых
множеств $A,B$ выполнено
$$
\mu(A\cap T^{n}B) \to  \mu (A)\mu (B)\ \ \ (n \,\to\,\infty). $$
\medskip

Автоморфизм $T$ обладает свойством слабого перемешивания, если выполнено
одно из (эквивалентных) условий:
$$
\exists n_i \to \infty \quad T^{n_i} \to \, \Theta \, .
$$
$$
\quad (T \otimes T)F=F\quad \Rightarrow \quad F=Const ;
$$
$$
T R = RT,\quad R\in{\cal R} \quad \Rightarrow R=с\Theta ;
$$

Из первого условия легко следует второе, а третье есть переформулировка
второго (рассмотрим $F$ как ядро регулярного оператора $R$).

Из последнего условия можно вывести, что отклонения от перемешивания
происходят на множестве нулевой плотности.
Пусть $I_n$ -- последовательность возрастающих по длине интервалов, и
$D_n\subset I_n$ -- последовательность подмножеств таких, что
 $\frac{| D_n|}{|I_n|}\to c>0$.
Из эргодичности $T$ имеем
$$\frac{1}{ |I_n|}\sum_{k\in I_n} T^k \to \Theta. $$
Пусть
$$\frac{1}{|I_n|}\sum_{k\in D_n} T^k \to R,$$
тогда $R \leq \Theta$, следовательно $R =c \Theta$.
Отсюда стандартным образом получаем,
что большинство $T^k$ при $k\in I_n$ близки   к   $\Theta$
(иначе нашли бы $R \neq c \Theta$).

{\bf Алгебраические операции над сплетающими операторами.}
 В доказательстве следующей леммы мы воспользуемся такими
 операциями как
 композиция, сопряжение, изменение типа сплетения.

{\bf ЛЕММА 1.1.1.} \it Пусть  внутренний  оператор
$P: L_2 \to L_2\otimes L_2$ сплетает автоморфизм
$T$ с автоморфизмом $S\otimes V$.  Известно, что
$S$ обладает свойством $S(3,4)$ (является тензорно простым).
Тогда оператор
$P$ тривиален:  $Im P=\{Const\}$.    \rm

Доказательство. Рассмотрим внутренний оператор
$PP^\ast : L_2\otimes L_2 \to L_2\otimes L_2$, он
сплетает $S\otimes S$ с ним же
(иначе говоря, коммутирует с ним):
$$PP^\ast(S\otimes V)=(S\otimes V)PP^\ast. $$

Теперь определим новый оператор
$Q: L_2\otimes L_2 \to L_2\otimes L_2$ следующим образом:
$$
  \left< Q(f\otimes g), f_1\otimes g_1\right>
= \left< PP^\ast (f\otimes f_1) , g\otimes g_1) \right> .
$$
Для этого оператора выполнено
$$Q(S\otimes S)=(V\otimes V)PP^\ast.$$
Также верно, что $Q$ -- внутрений оператор. Наконец, рассмотрим
внутренний оператор $Q^\ast Q$, коммутирующий с автоморфизмом
$(S\otimes S)$.  Если $S$ обладает свойством $S(3,4)$, то
$Q^\ast Q$ является тривиальным: $Im (Q^\ast Q)=\{Const\}$.
Отсюда вытекает тривиальность оператора $Q$ (см. ниже принцип
симметризации). Тривиальность оператора $Q$ влечет за собой
тривиальность $PP^\ast$.  Тем самым получили
тривиальность оператора $P$. 
(Эта лемма применяется в \S 5.5.)

               \it
{\bf ТЕОРЕМА 1.1.2.}  Свойство S(3,4) влечет за собой свойство  S(n,n+1).
               \rm

Доказательство. Пусть $T$ обладает свойством $S(3,4)$.  Пусть
внутренний оператор $P$ сплетает $T$ и $T\otimes\dots\otimes T$
($(n-1)$ сомножителей).
В качестве $V$ из леммы 1.1.1 рассмотрим автоморфизм, изоморфный
$T\otimes\dots\otimes T$ ($(n-2)$ сомножителей), и применим лемму.
Получим, что $P$ -- тривиальный оператор.

{\bf Сплетения слабых замыканий.}
Пусть некоторая последовательность степеней $T^{n_i}$ слабо сходится к
оператору $R$, а последовательность $\s^{n_i}$ к оператору $Q$.
Если выполнено $TP_1 = P_2S$, то будет выполняться
$RP_1 =P_2Q$. Этим  утверждением
 мы будем многократно пользоваться в дальнейшем.

{ \bf Разложимость и  неразложимость сплетений.}
    Говорим, что  оператор $A$ неразложим в полугруппе $\M$, если
условия $A = aA_1+(1-a)A_2$ и $A_1,A_2\in \M$ влекут за собой  $A_1= A_2$ при
$a(1-a)\neq 0.$
     Марковский оператор $P$, сплетающий  эргодические
автоморфизмы $T\in \cal{A}$ и
$S\in \cal{A}$, однозначно (mod 0)
представляется в виде интеграла
$$
 \int P_{\alpha} d\lambda(\alpha ),
$$
где $\lambda$ -- некоторая вероятностная мера на полугруппе ${\cal M}$,
сосредоточенная на
$\{P_\alpha \}$ -- множестве неразложимых операторов, сплетающих $T$ и
 $S$.

Сказанное  -- прямой аналог разложения джойнинга на эргодические компоненты.
\medskip

{\bf Симметризация сплетений.}
Принцип симметризации: \it  пусть $P\in \M$ и
 $P^\ast P=\Theta $, тогда $P=\Theta$.\rm

Это совсем простое утверждение имеет
ряд замечательных применений. Приведем  некоторые из них.
Следующий факт -- хорошо известное утверждение, доказательство
которого обычно использует спектральные аргументы. Наше доказательство
работает для любых групповых действий.

\bf ТЕОРЕМА. \it
Если $T\otimes T$ и $S$ эргодические автоморфизмы,
то $T\otimes S$
также эргодический  автоморфизм
(подразумевается эргодичность $T\otimes T$ относительно
меры $\mu\otimes\mu$).   \rm

Доказательство.  Пусть  $\mu\otimes\mu(F)>0$ и $(T\times S)F=F$.
Рассмотрим  интегральный оператор
$R:L_2 \to L_2 \otimes L_2$ с ядром $K(x,y)=\chi_F(x,y)/\mu\otimes\mu(F)$:
$$
     Rf(y)=\int_X K(x,y)f(x)d\mu(y)
$$
(оператор $R$ принадлежит полугруппе $\cal{R}$).
Так как $K(Tx,Sy)=K(x,y)$, оператор $R$ сплетает $T$ и $S$:
$RT =S R$. Получим
$$
                R^\ast RT =R^\ast S R= T R^\ast R,
$$
следовательно ядро $H(x,y)$ оператора $R^\ast R$ является
$T\otimes T$-инвариантным.
Из эргодичности  $T\otimes T$ получим  $H(x,y)=1$,
   $R^\ast R=\Theta $. Но это влечет за собой $R=\Theta$, т.е.
ядро интегрального оператора $R$ есть константа, равная 1. Таким образом,
мы доказали, что любое множество  $F$, удовлетворяющее
условиям $\mu\otimes\mu(F)>0$ и $(T\times S)F=F$ совпадает
(mod0) с $(X\times X)$. Следовательно, $T\times S$ --  эргодический
автоморфизм.   \\ 

Принцип симметризации вытекает из следующего утверждения.
\medskip

{\bf ЛЕММА 1.1.3.} \it Пусть $P_j^\ast P_j \,\to \Theta$,
тогда $P_j f\,$ сходится к
$\Theta \! f$ по норме $L_{2}(\mu)$.
\medskip

\bf
Следствие. $P^\ast P=\Theta \Rightarrow P=\Theta.$  \rm
\medskip

Доказательство. $$\|P_j (f-\Theta {f})\|^2
= \, \left< P_j^\ast P_j (f-\Theta {f}),(f-\Theta {f})\right>  \,\to 0.
$$    

Приведем утверждение, которое обобщает
теорему Блюма и Хансона из \cite{BH}.
 Последовательность
 $\lbrace a_{z}^j \rbrace , z\,\in\,Z ,j\,\in\,\bf N$,
будем назывть диссипативной, если
$$
\sum_z a_{z}^j =1, a_{z}^j \geq 0,$$

 $$ max_z \lbrace a_{z}^j \rbrace
\,\to \,0,\,j\,\to \,\infty .
$$
\medskip

\bf ТЕОРЕМА 1.1.4. \it Пусть $\T \in \,\cal{A}$ -- перемешивающий оператор:
$\T^n\to\Theta$, $n\to\infty$.
Для любой диссипативной последовательности
 $\lbrace a_{z}^j \rbrace$ выполнено
$$
  \| \sum_z a_{z}^j \T^zf - \Theta \! f\, \| \,\to \,0\quad .
$$
\\ \rm
\medskip

Доказательство. Положим $P_j = \sum_z a_{z}^j \T^z.$ Имеем
$$P_j^\ast P_j = \sum_w b_{w}^j \T^w ,$$
где последовательность $\lbrace b_{w}^j \rbrace $ также диссипативна.
Действительно,
$$   b_{w}^j \leq \sum_z a_{w-z}^j a_{z}^j \leq max_z {a_z^j } \to \,0
.$$
Получили $P_j^\ast P_j \to \Theta $, следовательно,
$\|P_jf-\Theta f\|\to 0$. \\ 
\medskip
\medskip

Изложим принцип симметризации в терминах джойнингов.
Мера  $\nu$ на $X\times Y$ со стандартными проекциями на сомножители
$X$ и $Y$
имеет следующее представление:
$$
    \nu(A\times B) =\int_B\nu_y(A)d\mu(y),
$$
где $\nu_y$ -- условные меры на $X$.
Пусть $\nu'$ -- другая мера на $X\times Y$.
Тогда можно определить их относительно независимое произведение
$$\nu\times_Y \nu':= \int_Y (\nu_y\otimes\nu_y')d\mu(y). $$
Если $\nu$ и $\nu'$ являются джойнингами пары $(T,S)$, то
$\nu\times_Y \nu'$ есть  джойнинг пары $(T,T)$.
Джойнингу $\nu$ мы сопоставляем  оператор   $P:L_2(X,\mu)\to L_2(Y,\mu)$:
$$
  P f(y)=\int_X f(x)d\nu_y(x).
$$
Ввиду равенств
$$
    \nu\times_Y \nu'(A\times B)= \int_Y \nu_y(A)\nu_y'(B) d\mu(y)=
$$ $$   = \left< P\chi_A | P'\chi_B \right> =
    \left< {P}'^{\ast} P\chi_A | \chi_B \right> .
$$
получим
$$\nu\times_Y \nu'= \mu\otimes\mu \
\Longleftrightarrow \ {P}'^{\ast} P=\Theta.$$

 Принцип симметризации теперь выглядит так:
$$\nu\times_Y \nu =\mu\otimes\mu\ \Longleftrightarrow \nu=\mu\otimes\mu,
$$
что вытекает из эквивалентностей
$$\nu\times_Y \nu=\mu\otimes\mu  \Longleftrightarrow  P^\ast P=\Theta
  \Longleftrightarrow  P=\Theta
 \Longleftrightarrow  \nu=\mu\otimes\mu .
$$

\begin{center}
{\bf \S 1.2. Дополнительная симметрия}
\end{center}
\medskip
Принцип дополнительной симметрии связан со следующим простым утверждением:
{\it если $QT=Q$ для  эргодического оператора $T$,   то $Q=\Theta $.}
Сказанное -- одна из формулировок эргодичности $T$.
  Ниже предлагается небольшое обобщение этого утверждения.
\medskip

{\bf ЛЕММА.}
    \it Если $T$ --  эргодический  оператор, и для некоторой
последовательности $P_j \in \cal{M}$ выполнено $(T P_j -  P_j) \to  0$,
то $P_j \to \Theta$.  \rm
\medskip

 Пусть $P_{j'} \to P$ (пользуемся компактностью
полугруппы $\M$). Тогда $T P=P$. Как отмечалось выше,
в силу   эргодичности $T$ получим
$P=\Theta$.

Следующее утверждение -- вариант хорошо известной операторной эргодической
теоремы  фон Неймана.
\medskip

{\bf ТЕОРЕМА. }\it Пусть $V_N = \frac{1}{N} \sum_{i=0}^{N-1} T^i$,
где  $T$ -- эргодический оператор,  тогда $V_N\,f$ сходится
$\Theta\! f$ по норме в $L_2(\mu).$  \\   \rm
\medskip

Доказательство. Рассмотрим последовательность операторов
$$
  V_{N}^\ast V_N = \frac{1}{N^2} \sum_{i=-N}^N (N-|i|)T^i .
$$
Легко проверить, что
$$
(T V_{N}^\ast V_N - V_{N}^\ast  V_N )\,\to \,0.
$$

Применив для последовательности $P_N = V_{N}^\ast\,V_N$
лемму, получим $V_{N}^\ast\,V_N \to\Theta$.
Принцип симметризации теперь дает утверждение теоремы. \\
\medskip

 {\bf ТЕОРЕМА 1.2.1.} (Принцип дополнительной симметрии)
\it Пусть оператор $J : L_2\otimes L_2\to L_2$ удовлетворяет
условию
$$
J({\bf 1}\otimes L_2) = J(L_2\otimes {\bf 1}) =
                             \lbrace Const \rbrace\,
$$
и  выполняется одно из равенств:
$$
                         RJ = J(I\otimes S),
$$
$$
                         J = J(R\otimes S).
$$
Если автоморфизм $R$ или $S$ является слабо перемешивающим,
то  $Im J=\{Const\}$ (оператор $J$ тривиален).     \rm
\medskip

 Доказательство. Обозначим ${\tilde Q} = J^\ast J$.
Пусть $S\otimes S$  эргодический автоморфизм
( $S$ является слабо перемешивающим).
Рассмотрим оператор $Q : L_2\otimes L_2\to L_2\otimes L_2$,
 который определен следующим образом:
$$
   \left< Q(f\otimes f_1),(g\otimes g_1)\right>
 = \left< {\tilde Q}(f\otimes g), f_1\otimes g_1\right>  .
$$
Получим  $Q(S\otimes S) = Q$, значит  $Q=\Theta\otimes\Theta$,
${\tilde Q} = \Theta\otimes\Theta$, $$J(L_2\otimes L_2) = \{Const\}.$$

В качестве приложения принципа дополнительной симметрии
приведем наше доказательство следующей теоремы Фюрстенберга.

ОБОЗНАЧЕНИЕ. До конца этого параграфа $A_i$ одновременно обозначает
измеримое множество и его индикатор.
\medskip

\bf ТЕОРЕМА 1.2.2.  \it Пусть $T$ --
слабо перемешивающий автоморфизм
($T\otimes T$ эргодичен), тогда  \rm
$$
      \lim_{N\to\infty} \frac{1}{N}\left< A_0,\sum^{N}_{i=1} (T^{i}A_1)\dots (T^{ki} A_k)\right>  =
\mu (A_0)\mu (A_1)\dots \mu (A_k).
$$
\rm
\medskip

Доказательство. Пусть    $N_i$ -- такая  последовательность, что
предел
$$
 \lim_{N_i\to \infty}\frac{1}{N_i}
\left< A_0,\sum^{N_i}_{i=1} (T^{i}A_1)\dots (T^{ki} A_k )\right> ,
$$
существует для всех
$A_0,\dots , A_k$.

Используя индуктивное предположение, имеем
$$
      \lim_{N\to\infty} \frac{1}{N}\left< \1,
\sum^{N}_{i=1} (T^{i}A_1)\dots (T^{ki} A_k)\right>  =
$$
$$
     \lim_{N\to\infty} \frac{1}{N}\left< A_1,
\sum^{N}_{i=1} (T^{i}A_2)\dots (T^{(k-1)i} A_k)\right>  =
\mu (A_0)\mu (A_1)\dots \mu (A_k).
$$
Положим $F=(T\otimes T^2\dots \otimes T^k)$.
Сказанное выше позволяет определить
 оператор $Q$:
$$
 \left< QA_0,  A_1 \otimes \dots \otimes A_k \right>  =
 \lim_{N_i\to \infty}\frac{1}{N_i}
\left< A_0,\sum^{N_i}_{i=1} (T^{i}A_1)\dots (T^{ki} A_k )\right> .
$$
Получим
$Q = FQ$, причем автоморфизм $F$  эргодичен  относительно меры $\mu^{\otimes k}$,
так как $T$ слабо перемешивающий.
Таким образом, из принципа дополнительной симметрии получим, что
 $Im(Q)$ -- одномерное пространство постоянных
функций.  Получили
$$
\left< QA_0\, ,\,  A_1 \otimes \dots \otimes A_k \right>  =  $$ $$
= \mu(A_0)\left< \1 \otimes \dots \otimes \1 , A_1
\otimes \dots \otimes A_k \right>  =
    \mu(A_0)\mu(A_1 )\dots \mu(A_k).
$$
Так как из любой последовательности мы можем выбрать  подпоследовательность
$N_i$ такую, что  для всех     $A_0,\dots , A_k$ существует предел
$$
 \lim_{N_i\to \infty}\frac{1}{N_i}
\left< A_0,\sum^{N_i}_{i=1} (T^{i}A_1)\dots (T^{ki} A_k )\right>,
$$
мы тем самым доказали, что
$$
      \lim_{N\to\infty} \frac{1}{N}\left< A_0,\sum^{N}_{i=1} (T^{i}A_1)\dots (T^{ki} A_k)\right>  =
\mu (A_0)\mu (A_1)\dots \mu (A_k).
$$

Аналогичный подход работает в доказательстве следующих теорем.
\medskip

 \bf ТЕОРЕМА 1.2.3. \it Для слабо перемешивающего потока
$\{T_t\}$ и для любых последовательностей
$a_1(i),a_2(i),...,a_k(i)$ таких , что $|a_p(i) - a_q(i)| \to \infty$ при
$1\leq p<q\leq k$     выполняется
$$
\int^{1}_0\mu (A_0\cap T_{sa_1(i)}A_1\cap\dots\cap T_{sa_k(i)}A_k)ds \to
      \mu (A_0)\mu (A_1)\dots\mu (A_k).
$$
\medskip

\rm Доказательство. Определим оператор $Q$:
$$
       \left< QA_0,  A_1 \otimes \dots \otimes A_k \right> =
   \lim_{i}\int^{1}_0\mu (A_0\cap T_{sa_1(i)}A_1\cap\dots\cap T_{sa_k(i)}A_k)ds
$$
Имеем
$$
   (I \otimes T_{c_1} \otimes \dots \otimes T_{c_k})Q = Q,
$$
где $c_m = \lim_{i'} \left(a_m(i')/a_k(i')\right)$
(предел определен для некоторой последовательности $i'$).
Из принципа дополнительной симметрии получаем $Im(Q)=\{Const\}$.
 Следовательно,
$$
  \left< QA_0, A_1 \otimes \dots \otimes A_k\right>
=\mu (A_0)\mu (A_1)\dots\mu (A_k) .
$$
\\  

Следующая теорема дает некоторую полезную информацию о кратном
перемешивании потоков, входящих в действие группы Гейзенберга.
\medskip

 \bf ТЕОРЕМА 1.2.4. \it Пусть сохраняющие меру эргодические потоки
$\Psi_a ,\Phi_b , T_c$  удовлетворяют соотношению
$$\Psi_a\Phi_b= T_{ab}\Phi_b\Psi_a,$$
где  $T_c$ -- слабо перемешивающий поток, коммутирующий с
 $\Psi_a$ и $\Phi_b$. Тогда $\Psi_a$ и $\Phi_b$ обладают перемешиванием
всех порядков. \rm
\medskip

 Доказательство. Определим  оператор $P$:
$$
  \left< A,P(B\otimes C)\right>  =
\lim_{i'\to\infty} \left< A,\,\Psi_{m(i')}B \Psi_{n(i')}C)\right>.
$$
Предположим, что $ 0 < m(i') < n(i')$ и существует предел
 $$\lim_{i''} m(i'')/n(i'') = a$$
для $i''$ --  подпоследовательности
 $i'$.

Далее переобозначим    $i''=i$.
Заметим, что при $s\to 0$ выполнено
$$
\left< \Phi_sA,\,(\Psi_{m(i)}\Phi_sB)(\Psi_{n(i)}\Phi_sC)\right>  \to
\left< A,\,(\Psi_{m(i)}B)(\Psi_{n(i)}C)\right>.
$$
Следовательно,
$$
\lim_{\e\to 0}\lim_{i\to\infty} \frac{1}{\e}\int_{0}^\varepsilon
\left< \Phi_sA,\,(\Psi_{m(i)}\Phi_sB)(\Psi_{n(i)}\Phi_sC)\right> ds =
$$
$$
\lim_{i\to\infty} \left< A,\,(\Psi_{m(i)}B)(\Psi_{n(i)}C)\right>.
$$
Ввиду коммутационных соотношений  получим
$$
\left< \Phi_sA,\,(\Psi_{m(i)}\Phi_sB)(\Psi_{n(i)}\Phi_sC)\right>  =
\left< A,\,(\Psi_{m(i)}T_{sm(i)}B)(\Psi_{n(i)}T_{sn(i)}C)\right>.
$$
Таким образом,
$
\left< A,P(B\otimes C)\right> =
$
$$
\lim_{\e\to 0}\lim_{i\to\infty} \frac{1}{\e}\int_{0}^\varepsilon
\left< A,\Psi_{m(i)}T_{sm(i)}B
\Psi_{n(i)}T_{sn(i)}C)\right>  ds.
$$
 Выражение
$$ \frac{1}{\e}\int_{0}^\varepsilon
\left< A,\Psi_{m(i)}T_{sm(i)}B
\Psi_{n(i)}T_{sn(i)}C)\right>  ds
$$
при фиксированном $r$ и больших $i$ мало отличается от
$$ \frac{1}{\e}\int_{0}^\varepsilon
\left< A,\Psi_{m(i)}T_{(s-\frac{r}{n(i)})m(i)}B
\Psi_{n(i)}T_{(s-\frac{r}{n(i)})n(i)}C)\right>  ds.
$$
Поэтому при фиксированном произвольном $r\in\bf R$ получим
$
\left< A,P(B\otimes C)\right> =
$
$$
\lim_{\e\to 0}\lim_{i\to\infty} \frac{1}{\e}\int_{0}^\varepsilon
\left< A,\Psi_{m(i)}T_{sm(i)}T_{ar}B
\Psi_{n(i)}T_{sn(i)}T_{r}C)\right>  ds.
$$
Это  приводит к  равенству
$$
P(T_{ar}\otimes T_{r} )= P,
$$
что влечет за собой
$$\left< A,P(B\otimes C)\right> = \mu(A)\mu(B)\mu(C) .
$$
\\  
\begin{center}
{ \bf \S 1.3. Индуцированные джойнинги}
\end{center}
\medskip

Пусть
бистохастический оператор $J:L_2(X,\mu)\otimes L_2(X,\mu)\to L_2(X,\mu)$
удовлетворяет условиям:
$$
J({\bf 1}\otimes L_2) = J(L_2\otimes {\bf 1}) =
                             \lbrace Const \rbrace\, ,
$$
$$
                         RJ = J(S\otimes T),
$$
где  $R,S,T$ -- некоторые  автоморфизмы  $L_2(X,\mu)$.
Оператору $J$  сопоставим измеримую функцию
${\cal J}:X\to {\cal M},$ (где ${\cal M}$ -- полугруппа бистохастических
операторов, действующих из $L_2(\mu)$ в $L_2(\mu)$)
$$
\int_A\left< {\cal J}(x)B,C\right>  d\mu(x) =
\left< J(A\otimes B), C\right> .
$$
Получим
$$
{\cal J}(Rx) \equiv T^{-1}{\cal J}(x)S.      \eqno (\ast)
$$
Запись $\dots\equiv \dots$ означает в дальнейшем равенство для почти всех
$x\in X$.

Если  для перемешивающих автоморфизмов $R,S,T$  выполнено
тождество  $(\ast)$, то
измеримому семейству операторов $\J : X \to{} \M$
отвечает  джойнинг $\nu =(R\times S\times T)\nu$.
Связь  задается формулой
$$
\nu(A\times B\times C)= \int_A\left<\J(x)B,C\right>  d\mu(x).
$$
Теперь рассмотрим семейство операторов
$ {\cal H}_m :X \to{} \M $:
$$
  {\cal H}_m  \equiv \J^\ast (x)\J (R^mx)\equiv
  \J^\ast (x)T^{-m} \J (x)S^{m}.
$$
Выполняется тождество
$$
 {\cal H}_m(Rx) \equiv S^{-1}{\cal H}_m(x)S .
$$
Функции $ {\cal H}_m$ сопоставим меру $\eta_m$:
$$
\eta_m(A\times B\times C)= \int_A\left< {\cal H}_m(x)B,C\right>  d\mu(x).
$$
Таким образом мы определили новые джойнинги
$\eta_m$,   которые будем называть \it индуцированными. \rm Отметим,
что для них выполнено
$\eta_m = (R\times S\times S)\eta_m$
($\eta_m$ является   джойнингом набора   $(R, S, S)$).
Одно из применений  индуцированных джойнингов состоит в следующем:
 \it
если для почти всех $x\in X$ выполнено
 ${\cal H}_m(x) \to \Theta$, то
${\cal J}(x)\equiv\Theta.$                   \rm

Основные приложения индуцированных джойнингов изложены в главах 2, 4.

Напомним, что действие $\{T_g: g\in G\}$ на $(X,\mu)$ называется простым,
если  коммутирующий с ним оператор $P\in {\cal M}$
имеет представление
$$
    P=  a\int_{C(\{T_g\})} Sd\theta(S) + (1-a)\Theta ,      \eqno      (1.1)
$$
где    $C(\{T_g\})$ обозначает групповой централизатор действия;
$\Theta$ --  оператор ортопроекции на пространство констант.
Напомним, что действие называется слабо перемешивающим,
если мера $\mu\otimes\mu$
эргодична   относительно $T\otimes T$ (ее нельзя представить в
виде суммы  различных  $T\otimes T$-инвариантных мер).
Доказательство следующего утверждения покажет причины,
по  которым "тривиализуются" некоторые сплетения для
2-простых систем.
\medskip

{\bf УТВЕРЖДЕНИЕ. {}} \it Пусть $\{V_j\}$ -- некоторая
последовательность
слабо перемешивающих автоморфизмов, коммутирующих с 2-простым
слабо перемешивающим действием
$\{T_g\}$ группы $G$, причем $V_j\to Id$.
Пусть ${\cal J}(x)\in {\cal A}$ для всех $x$,
и для всех $g\in G$ выполнено
$$T_{g^{-1}}\J(x)T_g \equiv {\cal J}(T_g(x)).$$
Тогда найдется автоморфизм $S$, коммутирующий с
действием $\{T_g\}$ такой, что
$${\cal J}(x)\equiv S.$$
\rm

{ Доказательство.  }
Рассмотрим операторнозначную функцию
$$
 {\cal H}_j(x) \equiv V_j \J(x) V_{j}^{-1} \J^{-1}(x) .
$$
С функцией $ {\cal H}_j$ можно связать динамическую систему
$\{T_g\otimes T_g\otimes T_g, X\times X'\times X'', \eta_j\},$
где мера $\eta_j$ определена  формулой
$$
\eta_j(A\times B\times C)= \int A(x)\left<   {\cal H}_j(x)B,C\right>  d\mu(x).
$$
Полученная динамическая система изоморфна  системе
$\{T_g\otimes T_g, X\times X',\mu\otimes\mu\}$.
Изоморфизм осуществляет отображение
$F: X\times X'\to X\times X'\times X''$,  определенное формулой
$$
F(x,x') =(x,x', {\cal H}_j (x)(x')),
$$
где $ {\cal H}_j (x)$ сейчас рассматривается как преобразование, действующее на
точку $x'$.
Следовательно, эта система эргодична вместе с (ее факторсистемой)
$\{T_g\otimes T_g, X'\times X'',\pi\eta_j\},$ где мера $\pi\eta_j$
есть проекция
меры $\eta_j$ на $X'\times X''$. Эргодической мере $\pi\eta_j$
отвечает
бистохастический оператор, коммутирующий с действием $\{T_g\}$.
Так как этот оператор имеет представление (1.1), а отвечающая ему мера
эргодична, имеем два случая:
$$
\int  {\cal H}_j (x) d\mu(x) = \Theta
$$
или
$$
\int  {\cal H}_j (x) d\mu(x) = S ,
$$
где $S$ -- автоморфизм, коммутирующий с нашим действием.
Так как при $j\to \infty$
$$
\int  {\cal H}_j (x) d\mu(x)  \to \int \J(x) \J^{-1}(x)  d\mu(x) = I ,
$$
($I$ -- тождественный оператор),  при достаточно больших $j$ первый
случай исключается. Из второго случая вытекает, что
$$
  {\cal H}_j (x)\equiv S,
$$
так как  оператор $S$ есть
крайняя точка в $\M$.

 Из $ {\cal H}_j (x)\equiv S$ получим
$$  \J (x)  \equiv V_{j}^{-1}S\J (x)V_j ,
$$
что эквивалентно условию
$$
J = J(V_j\otimes S^{-1}V_j )
$$
для оператора $J$, связанного с функцией  ${\cal J}$ формулой $(\ast)$.
Так как   $V_j$ слабо перемешивает, из принципа дополнительной симметрии
(теорема 1.2.1) получим, что  оператор $J$ тривиален.
Тогда из тривиальности $J$ вытекает, что $\J(x)\equiv \Theta$, что
противоречит $\J(x)\in \A$.  
\medskip
\medskip
\begin{center}
{\bf  1.4.\ Примеры тензорно простых систем}
\end{center}
\medskip

{\bf Слабое кратное перемешивание и $ \kappa$-перемешивание.}
Оператор   $S$ называется $ \kappa$-перемешивающим, $ \kappa\in (0,1)$,
если для некоторой последовательности
$n(i)\to\infty$ выполняется
$$
  S^{n(i)}\to (1- \kappa)I +  \kappa\Theta,
$$
где $I$ обозначает тождественный оператор,
$\Theta$ -- оператор ортопроекции на одномерное пространство $\{Const\}$
постоянных функций.   Докажем, что                       \it

  $ \kappa$-перемешивающее  действие   при $0<  \kappa < 1$
    принадлежит    классам $S(n-1,n)$.                   \rm

Рассмотрим случай $n=3$ (при $n>3$ рассуждения аналогичны).
Пусть
для некоторой последовательности  $\{T_{g(i)}\}$ выполнено
$$T_{g(i)}\to P_ \kappa = (1- \kappa)I +  \kappa\Theta  .$$

Обозначим через $J:L\to L^{\otimes 2}$ внутренний оператор, сплетающий
действие с ее тензорным квадратом:
$$
   (T_{g(i)}\otimes T_{g(i)}) J = JT_{g(i)}.
$$
Имеем
$$
     JH \subset H^{\otimes 2}, \quad
   (P_{ \kappa}\otimes P_{ \kappa})J = JP_{ \kappa}.
$$
Отсюда для $f\in H$ получим
$$
     (1- \kappa)^2Jf=(1- \kappa)Jf,
$$
следовательно, $JH=\{0\}$. Таким образом, $J=J\Theta$, значит наше
действие обладает свойством $S(2,3)$. \\ 

Как следствие получаем следующее утверждение:\\
 \it
  $ \kappa$-перемешивающий  автоморфизм $T$  обладает
свойством слабого премешивания кратности 2:
если $T^{m_i}$, $T^{n_i}$,  $T^{m_i-n_i}\to\Theta$,
то выполняется
$$
  \mu(A\cap T^{m_i}B \cap T^{n_i}C)\  \to \  \mu(A)\mu(B)\mu(C)
$$
для  любых измеримых множеств $A,B,C$.
 \rm

\bf  Другой пример $\Z$-действия класса $S(3,4)$. \rm
Используя метод аппроксимаций \cite{KS}, несложно построить
эргодический автоморфизм $T$, удовлетворяющий условию:
для некоторой последовательности $k(i)\to \infty$
$$
\forall A,B \ \ \   \mu(T^{k(i)}A\cap B)\to \frac{1}{2}\mu(A\cap B) +
\frac{1}{2}\mu(TA\cap B),
$$
что эквивалентно слабой сходимости
$T^{k(i)}\to Q= \frac{1}{2}I+\frac{1}{2}T$.
Отметим также, что
 этим свойством обладает популярный в теории
джойнингов автоморфизм  Чакона.
(Из эргодичности такого автоморфизма $T$ следует, что он
обладает свойством слабого перемешивания.)

Докажем, что автоморфизм $T$ является тензорно простым.
Пусть $\nu$ -- его некоторое эргодическое  самоприсоединение класса $M(3,4)$.
Pассмотрим оператор
$P:L_2 (\mu\otimes\mu\otimes\mu)\to L_2(\mu)$,
соответствующий мере $\nu$:
$ \forall f_1,f_2,f_3, f\in L_2(\mu)$
$$
\left< P(f_1\otimes f_2 \otimes f_3)\,,\,f\right> =
\int_{X\times X\times X\times X} f_1\otimes f_2 \otimes f_3\otimes f d\nu .
$$

Из проекционных свойств меры $\nu$ вытекает,
что $ P(f_1\otimes f_2 \otimes f_3)$ является постоянной функцией,
если постоянна одна из функций $ f_1, f_2, f_3$.
Из инвариантности $\nu$ относительно $T\otimes T\otimes T\otimes T$
вытекает условие сплетения
$
 TP =  P(T\otimes T\otimes T),
$
следовательно,
$$
T^{k}P =  P(T\otimes T\otimes T)^{k},  \ \ \
QP =  P(Q\otimes Q\otimes Q).
$$
Теперь    получим
$$
\frac{1}{2}(I + T)P =
\frac{1}{2} P(I\otimes I\otimes I) +
\frac{1}{2}P(T\otimes T\otimes T) =
$$
$$ = P(Q\otimes Q\otimes Q)=  \frac{1}{8}\left( P +
P(I\otimes T\otimes T)+ \dots +  P(T\otimes T\otimes T)\right).
$$
Приняв во внимание эргодичность  меры $\nu$ (эргодичность эквивалентна
неразложимости  оператора $P$ в сумму двух
различных сплетающих марковских операторов), имеем равенства вида
$$ P =  P(T\otimes I\otimes I) $$
или
$$ P =  P(T\otimes T\otimes I). $$
Так как $T^{m(i)}\to \Theta$, получим в первом случае
$$ P(f\otimes g\otimes h) =
P(\Theta f\otimes g\otimes h)= \Theta f\Theta g\Theta h,$$
что влечет  $\nu= \mu\otimes\mu\otimes\mu\otimes\mu$. Второй случай
рассматривается аналогично и также приводит к
равенству $\nu= \mu\otimes\mu\otimes\mu\otimes\mu$.
Таким образом, автоморфизм $T$ является тензорно простым.
\medskip

{\bf ТЕОРЕМА 1.4.1.} {\it Неперемешивающий автоморфизм с минимальным
централизатором является тензорно простым.}
\medskip

 Доказательство. Предположим, что
$$
                         TP = P(T\otimes T),
$$
где $P$ -- нетривиальное сплетение.
Если $T$ не обладает свойством перемешивания, то
найдется последовательность $T_i=T^{n_i}\to Q\neq \Theta$. Из
$$
                         T_iP = P(T_i\otimes T_i),
$$
получим
$$
                         QP = P(Q\otimes Q).
$$
Так как $Q$ коммутирует с автоморфизмом $T$, обладающим свойством
минимальности централизатора, имеем
$$Q=\sum_z a_z T^z + a\Theta ,$$
где $a<1$ и $ \sum a_z = 1-a$.  Рассмотрим случай $a=0$.
 Из  $QP = P(Q\otimes Q)$ получим
$$
  (\sum a_zT^z )P = P(\sum a_zT^z \otimes T^z ) .
$$
Так как  $ T^zP = P(T\otimes T)^z$, получим
$$
   P(\sum_v a_vT^v\otimes T^v) = P((\sum_z a_zT^z)\otimes(\sum_w a_wT^w)) .
$$
Так как  $P$ неразложим,
оператор $P(T^z\otimes T^w)$ также будет неразложимым.
Если
$a_z,a_w >0$ для некоторого $z\ne w$, то найдется $v$, для которого
$$
   P(T^v\otimes T^v) = P(T^z\otimes T^w) .
$$
Это равенство, ввиду  принципа дополнительной симметрии,
показывает, что оператор $P$ тривиален. Действительно, условие
$$
   P(T^{v-z}\otimes T^{v-w}) = P
$$
( $T^{v-z}$ или $T^{v-w}$ является слабо перемешивающим автоморфизмом)
вместе с условием
$$
   P({\bf 1}\otimes L_2) = P(L_2\otimes {\bf 1}) = \{ Const \}
$$
приводят к  равенству $P(L_2\otimes L_2)= \{ Const \}$.

Таким образом, для некоторого
$z$ выполнено $a_z = 1$, и, следовательно,
 $T^{n_i-z}\to I$ (автоморфизм $T$ является жестким).
Но жесткие эргодические автоморфизмы имеют
несчетный централизатор (теорема А.М.Степина, см. \cite{KSS}).

Это противоречит условию того, что централизатор нашего
 $T$ есть множество
 $\{T^n : n\in {\bf Z}\}$.
Таким образом, $P$ тривиален.

Осталось привести доказательство в случае $0<a<1$.
Имеем
$$(\sum_z a_z T^z + a\Theta)P\ =\ P (\sum_z a_z T^z + a\Theta)\otimes
(\sum_z a_z T^z + a\Theta).$$
Так как оператор $P$ внутренний и мы предполагаем, что
$P$ -- нетривиальное сплетение, получим:
 в левой части  тривиальная компонента имеет вес $a$, а
в правой вес тривиальной компоненты равен $2a-a^2$.
Отсюда для некоторых $z,k,l$ имеем
  $$T^z\Theta P = P(T^k\otimes T^l),$$
что влечет  тривиальность $P$:  $P = \Theta P$,
 $P(L_2\otimes L_2)= \{ Const \}$.

Теперь мы обобщим предыдущую теорему, пользуясь принципами
слабого замыкания, разложимости сплетений, симметризации и
дополнительной симметрии.

\vspace{3mm}
{\bf ТЕОРЕМА 1.4.2.} \it  Если простое действие
коммутирует со слабо перемешивающим неперемешивающим автоморфизмом,
то это действие является тензорно простым.
\rm   \\

\vspace{3mm}
\bf  ЛЕММА 1.4.3. \it Если $S_{g(i)}\to P$, то найдется такая целочисленная
последовательность $k(i)$, что  $S_{g^{-1}(i)g(m(i))}\to P^\ast P.$
\rm

\vspace{3mm}
Доказательство.    Фиксируем функции  $f_1,f_2\in L_2(\mu)$,
натуральное
$i$ и $\eps >0$. Из определений слабой сходимости имеем, что
для некоторого натурального числа $k(i)$  для всех $m>k(i)$ выполнено
$$
 |\left< S_{g(i)}f_1, S_{g(m)}f_2 \right>  -
\left< Pf_1,Pf_2\right> |< \eps .
$$
Таким образом, для некоторой последовательности $m'(i)$ имеем
$$
\left< S_{g(i)}f_1, S_{g(m'(i))}f_2 \right>  \to \left< Pf_1,Pf_2\right> .
$$
Используя сепарабельность пространства $L_2(\mu)$, применяем диагональную
процедуру и получаем, что  для всех   $f_1,f_2 \in L_2(\mu)$
для некоторой диагональной последовательности $m(i)$

$
\left< S_{g(i)}f_1, S_{g(m (i))}f_2 \right> $
$$ =\left< S_{g^{-1}(m (i))}S_{g(i)}f_1,f_2 \right>
 \to \left< P^\ast Pf_1,f_2\right> .
$$

Теперь  докажем теорему 1.4.2.
Пусть действие коммутирует со слабо перемешивающим автоморфизмом $S$,
который не является перемешивающим.
Будем считать, что наша система является $\Z$-действием, порожденным
автоморфизмом $T$.
Предположим, что $J$ сплетает $T$ и
$T\otimes T\otimes T,$ причем   $J$ -- внутренний неразложимый нетривиальный
оператор.

Пусть последовательность $n(i)$ такова, что $S^{n(i)}\to Q\ne\Theta $,
тогда по лемме  имеем   $S^{m(j)}\to Q^\ast Q \neq\Theta .$
Так как
$$
J^\ast(I \otimes S^{m(j)}\otimes I)J \to J^\ast(I\otimes Q^\ast Q\otimes I)J,
$$
нам достаточно доказать, что
$$U= J^\ast(I\otimes Q^\ast Q\otimes I)J\neq\Theta  .$$
Учитывая, что  $V=(I\otimes Q \otimes I)J\neq\Theta $ влечет
$$
U=V^\ast V\neq\Theta  ,
$$
достаточно убедиться в том, что оператор $V$ нетривиален. Это легко вытекает
из свойства простоты централизатора $T$. Действительно, стохастический
оператор $Q\neq\Theta $ коммутирует с простым действием, следовательно,
он имеет вид
$$
      Q = a\left(\int  R d\sigma(R)\right) + (1-a)\Theta ,
$$
где $a>0$.
Если оператор  $J$ нетривиален и неразложим, то операторы
$J_R = (I\otimes R\otimes I)J$
также неразложимы и нетривиальны.
Представление тривиального оператора в виде
интеграла $\int_{C(T)} J_R d\sigma(R)$ невозможно,
так как тривиальный оператор
является неразложимым. (Напомним, что тривиальному оператору отвечает
мера $\mu^{\otimes 4}$, эргодическая относительно $T^{\otimes 4}$.)

Таким образом, доказано, что
$J^\ast(I\otimes Q^\ast Q\otimes I)J\neq\Theta  .$
Так как $S^{m(j)}\to Q^\ast Q $, для некоторого $m=m(j)\neq 0$ получим
$$J^\ast(I\otimes S^m\otimes I)J\neq\Theta  .$$

 Определим меру $\eta$ на $X_{(1)}\times X_{(2)}\times Y$, $Y=X^3$ формулой
$$
\eta(A\otimes B\otimes C) = \int_Y (J\chi_A)\,(J'\chi_B)\,\chi_C\, d\lambda ,
$$
где $J'=(I\otimes S^m\otimes I)J$, $\lambda =\mu^{\otimes 3}$, $C\subset Y$.
Ввиду равенства
$$
 J^\ast J' = a\left(\int_{SC} Rd\sigma_{n}(R)\right) +(1-a)\Theta   ,
$$
выполненного вместе с условием $a>0$,
найдется  эргодическая компонента $\beta$ меры $\eta$
такая, что  $$\pi_{12}\beta = (R\times Id)\Delta ,$$
где $\pi_{12}$ обозначает
 проекцию на $X_{(1)}\times X_{(2)}$, а $\Delta$ -- диагональная мера
 на $X_{(1)}\times X_{(2)}$.

Пусть $\nu$ -- проекция меры $\beta$ на $X\times Y$, а
$\nu'$ -- проекция меры $\beta$ на $X'\times Y$.  Тогда выполнено
$$ (R\times Id)\nu=\nu'.$$
Мерам $\nu$  и $\nu'$ соответствуют операторы  $J$  и $J'$.
Следовательно, эти операторы связаны равенством $J'=JR$,
которое переписывается в виде
$$
(I\otimes S^m\otimes I)J= JR .
$$
Так как $S^m$ -- слабо перемешивающий оператор, а $J$ -- внутренний
оператор,
по принципу дополнительной симметрии  получим,
что $J$ есть тривиальный оператор.

Замечание. В приведенном только что рассуждении содержится доказательство
следующей важной леммы (\cite{GHR}), переформулированной нами на
 языке сплетений.
\medskip

{\bf ЛЕММА 1.4.4.} \it Пусть $S$ -- простая система, $T$ -- эргодическая, а
$P$,$Q$ -- их неразложимые марковские сплетения.
Тогда $Q^\ast P=\Theta$ или  $P=QR$ для   некоторого автоморфизма $R$,
коммутирующего с $S$.  \rm
\medskip

Действительно, рассмотрим  джойнинг
$$
\eta(A\otimes B\otimes C) = \int_Y (P\chi_A)\,(Q\chi_B)\,\chi_C\, d\lambda .
$$
Если   $Q^\ast P\neq\Theta$, то найдется  эргодическая компонента
меры $\eta$ с проекцией вида $(R\times Id)\Delta$. Как отмечалось выше,
последнее влечет равенство    $P=QR$.

Этой леммой мы воспользуемся в следующей главе.
\medskip
\begin{center}
{ \bf  1.5. Связь типов тензорной  простоты}
\end{center}
\medskip

Для краткости  обозначим свойство $S(n-1,n)$ через  $\s_{n}$.
Следующая теорема обобщает основной результат работы Кинга
\cite{King}. Это обобщение нам понадобится   в главе 3 при изучении
гладких джойнингов (теорема 3.1.1.)
\medskip

\bf ТЕОРЕМА 1.5.1. \it a) Для любых $q\geq 2$ и $k\geq 3$ свойство
  $\s_{2q}$ влечет свойство $\s_{k}$.
\\ b) Для всех $q,r>1$ выполнено $\s_{2q}=\s_{2r}$.
\rm
\medskip

 Доказательство.
Предварительно установим, что
свойство $\s_{4}$ влечет свойство $\s_{k}$.  Рассмотрим случай
$k=5$. Нам нужно показать, что для внутреннего
оператора $P$, сплетающего тензорный квадрат и тензорный куб действия
из класса $\s_{4}$, выполнено
$$
  \left< P(f_1\otimes f_2)\,,\, f_3\otimes f_4\otimes f_5)
  \right> =  0,
$$
как только одна функция из набора $f_1,f_2,\dots ,f_5$ имеет нулевое
среднее. Это равносильно следующему:
 для любых $f_1,f_2$ при $\int f_1 \int f_2 =0$ выполнено
$P(f_1\otimes f_2)= 0$ (оператор $P$ тривиален).
Чтобы доказать это, установим
$$
  \left< P(f_1\otimes f_2)\,,\, P(f_1\otimes f_2)
  \right> =  0
$$
для всех таких пар  $f_1,f_2$. Заметим, что
$$
  \left< P^\ast P(f_1\otimes f_2)\,,\, (f_1'\otimes f_2')
  \right> =  0,
$$
если одна из функций $f_1,f_2,f_1',f_2'$ является постоянной функцией.
Следовательно, оператор $P^\ast P$ является внутренним. Так как он
коммутирует с тензорным квадратом действия, а действие из класса
 $\s_{4}$, то оператор $P^\ast P$ тривиален. Последнее влечет за
собой, что тривиален и $P$. Так как выбор $P$ был произвольным,
получили, что действие принадлежит классу $\s_{5}$. Те же
рассуждения показывают, что  $\s_{4}$ влечет $\s_{k}$.
Для доказательства теоремы осталось установить, что
 $\s_{2q}$ влечет $\s_{4}$ для всех $q >2$.

Пусть $n\geq 4$ -- четное число. Установим для всех таких
$n$, что каждое из свойств
$  \s_4, \s_6,\dots , \s_{2n-2}$
влечет $ \s_{n}$.
Пусть $l+m = n, \, m\geq 2$ и $T\in \s_{2m}$. Наша цель -- доказать,
что тогда $T\in \s_{n}$.
Пусть  $ P:L^{\otimes m}\to L^{\otimes l}$ является внутренним оператором,
сплетающим  $T^{\otimes r}$ и $T^{\otimes m}$.
Тогда оператор $P^\ast P:L^{\otimes r}\to L^{\otimes r}$
также является внутренним. Поясним это. Если одна из функций
$f'_1,\dots , f'_r$ является константой, то
$P(f'_1\otimes \dots \otimes f'_r )$ есть постоянная функция.
При условии, что одна из функций $f_1,\dots ,f_r$ имеет нулевое среднее,
получим, что  $P(f_1\otimes \dots \otimes f_r)$ также имеет нулевое среднее
значение. Тогда выполнены равенства
$$
  \left< P^\ast P(f_1\otimes \dots \otimes f_r),\,
  f'_1\otimes \dots \otimes f'_r \right>  =   $$ $$
  \left< P(f_1\otimes \dots \otimes f_r),\,
P(f'_1\otimes \dots \otimes f'_r )\right> =   $$ $$
    \left< P(f_1\otimes \dots \otimes f_r),\, Const\right>  = 0 .
$$
Так как внутренний оператор $P^\ast P$ сплетает
$T^{\otimes m}$ c $T^{\otimes m}$, а наше действие принадлежит
классу $\s_{2m}$, получим, что  $P^\ast P$ -- тривиальный оператор.
Последнее эквивалентно  тривиальности оператора
$P$.

Теперь, пользуясь тем, что
$\s_4, \s_6,\dots , \s_{2n-2}$  влекут за собой $\s_{n}$,
докажем, что  $\s_{20}$ влечет $\s_{4}$.
Действительно, мы имеем
$$\s_{20}\Rightarrow \s_{12}\Rightarrow \s_{8}\Rightarrow \s_{6}
\Rightarrow \s_{4}.
$$
Таким образом, установлено, что свойство
$\s_{2m}$ влечет свойство $\s_{4}$, следовательно, они эквивалентны.

Вопрос: может ли динамическая система обладать нечетной тензорной
простотой, но при этом не являться системой класса $\s_4$ ?
Внутри некоторых классов динамических
систем установлено совпадение этих инвариантов.
 \bf Случай, когда  $\bf \s_{2p+1} \Rightarrow \s_{2q+1}.$  D-системы. {}\rm

Действие $\Psi$ группы $G$ назовем  D-системой,
если для любого неразложимого стохастического оператора $P$,
сплетающего $\Psi$ с другим
произвольным действием $\Phi$, выполняется одно  из условий:

1.\quad $ImP^\ast = \{Const\}$, что равносильно  $P=\Theta $,

2.\quad Образ $ImP^\ast$ плотен в  $L_2(\mu)$.
\\

Примерами D-систем служат автоморфизмы с
минимальным стохастическим централизатором.
Сформулируем  более общий факт.
\medskip

 \bf УТВЕРЖДЕНИЕ 1.5.2. \it Простая динамическая система,
не имеющая нетривиальных
факторов, является D-системой.
\medskip

 \rm
Доказательство.  Пусть $P$ -- неразложимое сплетение систем $T_g$ и
$R_g,$ а $\nu$ обозначает полиморфизм, отвечающий оператору $P$. Мера $\nu$
является джойнингом, эргодическим относительно $T_g\otimes R_g$.

Предположим, что $T_g$ -- простая система.
Рассмотрим меру $\eta$ на $X\times X'\times X''$, где $X=X'=X''$,
определенную следующим образом:
$$
\eta(C\times A\times B) = \int \chi_C\,(P\chi_A)\,(P\chi_B) d\mu.
$$
Из определения видно, что мера $\eta$ инвариантна относительно
$R_g\otimes T_g\otimes T_g.$

Отметим другое свойство меры $\eta$:
проекции $\pi'\eta$ и $\pi''\eta$ меры
$\eta$ соответственно на $X\times X'$ и $X\times X''$  совпадают с мерой
$\nu$.
Так как мера $\nu$ эргодична
относительно $R_g\otimes T_g$, получим, что почти все эргодические
(относительно $R_g\otimes T_g\otimes T_g$) компоненты меры $\eta$
обладают этим свойством.
Выберем одну из таких компонент, обозначив ее через $\sigma$.

Проекция эргодической
меры $\sigma$ на $X'\times X''$
есть эргодическая (относительно $T_g\otimes T_g$)
мера $\pi\sigma$. В силу простоты действия имеем альтернативу:
или эта мера равна  $\mu\otimes\mu$, или она сосредоточена на
графике автоморфизма, коммутирующего с действием.
В первом случае имеем:
$$\pi\sigma= \mu\otimes\mu, \quad P^*P=\Theta , \quad P=\Theta .$$
Так как мы предполагаем, что $P\neq\Theta ,$  первый случай
исключен.

Рассмотрим второй случай. Мера $\sigma$ сосредоточена на множестве
$X\times \Gamma_S$, где   $\Gamma_S$ -- график автоморфизма $S$,
коммутирующего с действием $T_g.$  Так как проекции меры $\sigma$ на
$X\times X'$ и $X\times X''$   совпадают с $\nu$, получим, что
$\nu = (I\times S)\nu,$ что эквивалентно  $S^{-1}P=P.$
Если автоморфизм $S$ эргодичен, то
$P=\Theta ,$ что противоречит предположению.
Если же $S$ не эргодичен, но отличен
от $I$, получим, что
нетривиальная алгебра $S$-инвариантных множеств является
$T_g$-инвариантной сигма-алгеброй
(следствие того, что $S$ коммутирует с нашей системой).
Таким образом, наше действие имеет нетривиальный фактор, что
противоречит условиям теоремы.

Подводя итог сказанному выше, получим, что мера $\pi\eta$
есть сумма мер $\mu\otimes\mu$ и диагональной меры $\Delta$, отвечающей
тождественному преобразованию:
$$\Delta(A\times B)=\left<\chi_A\,,\,\chi_B\right>.$$
Имеем
$$
\pi\eta(A\times B) = \left< PA,\,PB\right>  =
\left< P^*PA,B\right>  = $$ $$\left< (a(\Theta )+(1-a)I)\chi_A,\,\chi_B\right> .
$$
Таким образом, каждый неразложимый сплетающий оператор $P$
удовлетворяет условию       $P^*P=a\Theta+(1-a)I$,
что влечет D-свойство нашей системы.  
\newpage

 \bf ТЕОРЕМА 1.5.3. \it  Для D-систем свойство $\s_{2n+1}$ влечет свойство
$\s_3$.
\rm
\medskip

 Доказательство.  Покажем, что  $\s_5$ влечет $\s_3$. Предположим,
что динамическая система $\Psi$ не обладает свойством $\s_3$. Тогда она
не обладает    свойством $\s_4$, что пояснялось выше.
Следовательно, найдутся внутренние
операторы  $J:L^{\otimes 2}_2\to L_2$ и $P:L^{\otimes 3}_2\to L_2$,
сплетающие соответственно $(\Psi\otimes\Psi)$ и $(\Psi\otimes\Psi\otimes\Psi)$
c действием $\Psi$, причем образы операторов $J$ и
$P$ будут всюду плотны в $L_2$.
Тогда, очевидно,  оператор $P^{\ast}J$ является внутренним
нетривиальным оператором,  сплетающим $\Psi\otimes\Psi$ и
$\Psi\otimes\Psi\otimes\Psi$.

Докажем, что $P^{\ast}J$ в наших предположениях должен быть нетривиальным.
Так как оператор $P^{\ast}$ нетривиален,
некоторую функцию $f$ он переводит в непостоянную функцию;
у оператора $J$ образ всюду плотен в $L_2(\mu)$, следовательно,
некоторую функцию $f_1$ он переводит в функцию, близкую по норме к $f$.
Таким образом, получим, что $P^{\ast}Jf_1$ не является постоянной функцией.
Но это противоречит тому, что система
$\Psi$ обладает свойством $\s_5$.
Итак, получили, что оператор $P^{\ast}J$ нетривиален.
Общий случай $n>5$ доказывается аналогично.         .
\\ \bf Независимый фактор. \rm Если   выполнено
$$PT=SP ,\quad P^\ast P = I, $$
говорим, что $T$ является фактором $S$. А в случае, когда
при этом $P$ является внутренним оператором и  $ S=T\otimes T$, говорим,
что  $S$ имеет независимый $T$-фактор (независимый по отношению к
координатным факторам). Из определений вытекает, что
 $T\notin \s_3$.
\medskip

 \bf ТЕОРЕМА 1.5.4. \it  Если  $S=T\otimes T$ имеет независимый
$T$-фактор, то $T$ не обладает ни одним из свойств $\s_p$.
\rm
\medskip

Доказательство.
 Проведем рассуждение при  $p=5$. Случай нечетных $p>5$
рассматривается аналогично. Для четных $p$ уже все доказано:
$\s_{2m}$ влечет свойство $\s_{3}$, которым $T$ не обладает.

Нам нужно показать, что из свойства $\s_5$
 системы $\Psi$ вытекает отсутствие независимого $\Psi$-фактора
у системы $\Psi\otimes\Psi$.

 Предположим, что такой фактор найдется и ему отвечает внутренний
оператор $P:L^{\otimes 2}_2\to L_2$.
Рассмотрим новый оператор $J:L^{\otimes 4}_2\to L_2$, определенный через
$P$ при помощи подстановок формулой
$$
  \left< J(f_1\otimes f_2\otimes f_3\otimes f_4), f\right>  =
\left< P(f_1\otimes P(f_2\otimes P(f_3\otimes f_4))), f\right>.
$$

 Проверим, что  $J$ есть внутренний оператор (для этого
нам достаточно знать только то, что $P$ -- внутренний оператор.)
 Если, скажем, $f_4$ есть константа, то получим
$$
 P(f_3\otimes Const) = Const\Theta  f_3 ,
$$
так как оператор $P$ внутренний. Теперь последовательно получим
$$
 P(f_2\otimes P(f_3\otimes Const)) = Const\Theta  f_2\Theta  f_3,
$$
$$
P(f_1\otimes P(f_2\otimes P(f_3\otimes Const))) =
Const\Theta  f_1\Theta  f_2\Theta  f_3 .
$$
Таким образом, если одна из функций  $f_1, f_2, f_3, f_4$ есть константа,
имеем
$$
  J(f_1\otimes f_2\otimes f_3\otimes f_4) =
\Theta  f_1\Theta  f_2\Theta  f_3\Theta  f_4.
$$
Это означает, что оператор $J$ внутренний.

Следующая лемма показывает, что построенный оператор $J$ нетривиален.
\medskip

{\bf ЛЕММА.} \it  Пусть  образы внутренних операторов
$P:L_2(\mu)\otimes L_2(\mu)\to L_2(\mu)$
и $Q:L_2(\mu)\otimes L_2(\mu')\to L_2(\mu)$
плотны в $L_2(\mu)$. Тогда оператор   $J$, определенный подстановкой
$$
  \left< J(f_1\otimes f_2\otimes f')\,, f\right>  =
\left< P(f_1\otimes Q(f_2\otimes f'))\,, f\right>
$$
также имеет образ, плотный в  $L_2(\mu)$.
 \rm
\medskip

Доказательство леммы.
Линейные комбинации функций
вида  $Q(f\otimes f')$ плотны в  $L_2(\mu)$ ( то же самое
выполнено для  оператора $P$).
Так как образ $Q$ всюду плотен в $L_2(\mu)$, получим, что любая функция
из   $L_2(\mu)$ с любой точностью приближается линейными
комбинациями вида
$$
  \sum_{m,n}^{N} P(h_m\otimes Q(\sum_{i,j}^{N}f_{i,n}\otimes f_{j,n}))
  = \sum_{m,n,i,j}^{N} J(h_m\otimes f_{i,n}\otimes f_{j,n}) .
$$
Таким образом,  $Im(J)$ всюду плотен в $L_2(\mu)$.

 Теперь приходим к противоречию с тем, что $\Psi$ принадлежит классу $\s_5$.
Действительно,  если внутренний оператор $J$ сплетает
$\Psi^{\otimes 4}$ с  $\Psi$, то оператор $J$ тривиален.
Это завершает доказательство теоремы 1.5.4.
\newpage
\begin{center}
{\Large  \bf ГЛАВА 2 }
\end{center}
\begin{center}
\bf   СИСТЕМЫ С МИНИМАЛЬНЫМ, ПРОСТЫМ И КВАЗИПРОСТЫМ МАРКОВСКИМ
           ЦЕНТРАЛИЗАТОРОМ.
\end{center}
\rm
\medskip
\medskip

В этой главе устанавливается
тензорная простота систем с минимальным, простым и квазипростым
централизатором при некоторых дополнительных условиях.
В частности,  доказано, что   2-простой поток
 обладает
свойством простоты всех порядков;\ показано, что
квазипростые действия порядка 3 обладают свойством
квазипростоты всех порядков;\  для автоморфизмов установлено, что
свойства перемешивания
кратности 2 и свойство MSJ(2) (минимальных самоприсоединений порядка 2)
влекут за собой свойство  минимальных самоприсоединений
всех порядков. Таким образом, если  контрпример к проблеме
 ``MSJ(2) = MSJ?'' для $\Z$-действий существует, то
 он  дает отрицательное
 решение проблемы Рохлина о кратном перемешивании.
 Для потоков эта проблема решена положительно.
 Однако, для некоторых некоммутативных
действий доказано различие свойств  MSJ(2), MSJ(3) и MSJ(4).

\medskip
\begin{center}
{\bf  2.1. Простые  системы с несчетным централизатором}
\end{center}
\medskip

{\bf Простые самоприсоединения.}
Мера   $\Delta_S  = (Id\times S)\Delta$, где $S$ -- автоморфизм,
называется сдвигом диагональной меры $\Delta$.
 Действие $T$ называется $n$-простым (или простым порядка $n$),
если любой эргодический джойнинг
$\nu\neq \mu^{\otimes n}$ набора из  $n$ копий действия $T$
обладает свойством:
одна из его проекций  на двумерную грань  в
$X\times\dots\times X$ являеся мерой  $\Delta_S$ для
некоторого автоморфизма $S$, коммутирующего с действием $T$.

Неизвестно, влечет ли 2-простота $\Z$-действия за собой простоту
порядка 3.
Мы даем частичный ответ:
\medskip

{\bf УТВЕРЖДЕНИЕ 2.1.1.} \it
Если автоморфизм $T$ является простым порядка 2 и
коммутирует с несчетной группой слабо перемешивающих
автоморфизмов (исключая тождественный),
то  $T$ является тензорно простым.
\rm
\medskip

Непосредственным следствием этого утверждения является

{\bf ТЕОРЕМА 2.1.2.}   \it
Слабо перемешивающий 2-простой поток  является  простым всех порядков.
Перемешивающий 2-простой поток перемешивает с любой кратностью.
\rm
\medskip

На множестве марковских операторов, сплетающих автоморфизмы
$T$ и $S$, введем эквивалентность  $P\sim P'$: если
 для    некоторого автоморфизма $R$, коммутирующего с $T$,
выполнено $P'=PR$.
\medskip

{\bf ЛЕММА 2.1.3. }   \it
  Число классов эквивалентности на множестве
неразложимых сплетений пары $(T,S)$,
где  $T:X\to X$ 2-простой, а  $S:X\to X$ -- эргодический автоморфизм,
не более, чем счетно.
\rm
\medskip

 Доказательство леммы.
Воспользуемся следующим вспомогательным утверждением
( лемма 1.4.4.): если $T$ -- простая система, $S$ -- эргодическая, а
$P$,$P'$ -- их неразложимые марковские сплетения,
то $P^\ast P'=\Theta$ или  $P'=PR$ для   некоторого автоморфизма $R$,
коммутирующего с $T$.

Если  выполнено $P'=PR$, то
$P L_2^0 = P RL_2^0 = P 'L_2^0$ (мы обозначаем через $L_2^0$  пространство
функций с нулевым средним).

Если  $P^\ast P'=\Theta$ и $f,g$ -- функции
с нулевым средним, то выполнено
$$
0=\langle P^\ast P' f\ |\ g \rangle  =  \langle  P' f\ |\ Pg \rangle .
$$
 Утверждение теоремы теперь вытекает из
того, что  число попарно ортогональных
подпространств сепарабельного пространства $L_2$ не более, чем счетно.
\\  
  \medskip

 Доказательство утверждения  2.1.1.
Пусть $\nu$ -- попарно независимый нетривиальный эргодический джойнинг
набора $(T,T\times T)$. Предположим, что мера $\nu$ сингулярна
относительно меры $\mu\otimes\mu\otimes\mu$.
Рассмотрим континуальное семейство
$(Id\times (Id\times S_g))\nu$, где $\{S_g\}$ -- несчетная группа
слабо перемешивающих автоморфизмов, исключая тождественный.
Так как  число классов эквивалентных
джойнингов не более, чем счетно,
для несчетного  множества  автоморфизмов $S_g$
найдется  джойнинг $\lambda$ такой, что
$$
    (Id\times (Id\times S_g))\nu =  (R_g\times (Id\times Id))\lambda.
$$
Получим
$$
 (R_g^{-1}\times (Id\times S_g))\nu =  (R_h^{-1}\times (Id\times S_h))\nu,
$$
откуда вытекает равенство
$$
 (R_hR_g^{-1}\times (Id\times S_h^{-1}S_g))\nu = \nu.
$$
Автоморфизм  $S=S_h^{-1}S_g$ слабо перемешивающий: $S^{i_k}\to \Theta$
для некоторой последовательности ${i_k}$.
Последнее равенство влечет за собой
$$
\nu= (R_h^iR^i_{g^{-1}}\times (Id\times S^i_{h^{-1}}S^i_g))\nu.
$$
Теперь можно доказать, что для некоторого марковского оператора
$\bar{R}$  выполнено
$$
\nu = (\bar{R}\times (Id\times \Theta))\nu .
$$
Действительно,
 обозначая $R=R_hR_g^{-1}$
 имеем для всех $i\in \N$
$$
  \nu = (R^i\times (Id\times S^i))\nu .
$$
Используя компактность марковской полугруппы,
для некоторой последовательности $i_k\to\infty$ получим
 слабые сходимости $S^{i_k}\to \Theta$ и  $R^{i_k}\to\bar{R}$,
где $\bar{R}$ -- некоторый марковский оператор.
Тогда выполнено $R^{i_k}\otimes S^{i_k}\to \bar{R}\otimes \Theta$,
и мы получаем
$$
    \int\chi_A\otimes\chi_B\otimes\chi_C d\nu =
    \int \bar{R}(\chi_A)\otimes\chi_B\otimes \Theta(\chi_C) d\nu =
    \mu(A)\mu(B)\mu(C),
$$
$$\nu =  \mu\otimes\mu\otimes\mu.
$$

Таким образом, возможен лишь случай $\nu = \mu\otimes\mu\otimes\mu$.

Аналогично, если  $\nu$ -- попарно независимый эргодический джойнинг
набора $(T,T\times T\times T)$, мы показываем, что мера $\nu$ есть
 $\mu\otimes\mu\otimes\mu\otimes\mu$.   Последнее
 означает, что $T$ является тензорно простым.
\\  
\medskip
\medskip

\begin{center}
{ \bf  2.2. Наследственная независимость и квазипростота действий}
\end{center}
\medskip

{\bf Квазипростые самоприсоединения.}
Джойнинг $\nu$ пары $(T,T)$  называется квазидиагональной мерой, если
для почти всех $x,y$ условные меры  $\nu_x$ и $\nu_y$, возникающие
в представлении
$$
\nu(A\times B)=\int_A \nu_x(B)d\mu(x) = \int_B \nu_y(A)d\mu(y),
$$
имеют вид:
$$
    \nu_x = \frac{1}{p} (\delta_{y_1(x)}+
\delta_{y_2(x)}+\dots + \delta_{y_p(x)}),
$$
$$
    \nu_y = \frac{1}{q} (\delta_{x_1(y)}+
\delta_{x_2(x)}+\dots + \delta_{x_q(y)}).
$$
Говорят, что
\it
действие $T$ квазипростое порядка $n$, если
для любого эргодического джойнинга
 $\nu\neq \mu^{\otimes n}$
набора из  $n$ копий действия $T$ выполнено:
одна из  проекций  на двумерную грань  в
$X\times\dots\times X$ являеся  квазидиагональной мерой.
\rm
Приведенные выше определения автоматически распространяются на произвольные
групповые  действия.

В операторных терминах  свойство 2-квазипростоты формулируется так:
кроме $\Theta$ все неразложимые марковские операторы $P$,
коммутирующие с системой, принадлежат
полугруппе ${\cal D}$: $P^\ast P\geq aI$.

Мы рассмотрим системы с так называемой \it наследственной
независимостью, \rm называя их  HI-системами.

ОПРЕДЕЛЕНИЕ.
Действие $T$ по определению обладает HI-свойством, если для любой
эргодической  динамической  системы, порожденной тремя факторами
$\F, \F',\F''$, каждый из которых изоморфен системе $T$,
выполнено
$$
    \F\bot\F' \ \ \& \ \  \F\bot\F'' \ \
\Rightarrow \ \ \ \F\ \bot\ (\F'\bigvee \F''),
$$
\rm
(где $\bot$ обозначает независимость факторов).

Это определение эквивалентно следующему.
Действие $T$ обладает HI-свойством, если \it для любого эргодического
джойнинга  $\eta$ тройки  $(T,T,T)$  выполнено
$$
    \pi_{12}\eta =\mu\otimes\mu \ \ \& \ \ \pi_{13}\eta =\mu\otimes\mu
\ \ \Rightarrow \ \ \eta = \mu\otimes \pi_{23}\eta,
$$
где  $\pi_{ij}\eta$ --  проекция
 меры  $\eta$
на грань  $X_i\times X_j$ куба $X_1\times X_2\times X_3.$
\rm

Нетрудно видеть, что 3-простая  система $T$ обладает  HI-свойством.
Действительно, если $\pi_{23}\eta=\mu\otimes\mu$, то из 3-простоты мы
получаем $\eta=\mu\otimes\mu\otimes\mu$.
Иначе, если $\pi_{23}\eta$ -- сдвиг диагональной меры,
из явного описания всех эргодических джойнингов порядка 3 получим
$$
\eta(A\times B\times C) = \mu(A)\mu(B\cap SC),
$$
где   $S$ -- некоторый автоморфизм, коммутирующий с $T$.
\medskip

{\bf ЛЕММА 2.2.1.} \it
Квазипростота порядка 3 влечет за собой  HI-свойство.
\rm
\medskip

{  Доказательство}. Пусть $T$ -- 3-квазипростое действие и $\nu$ --
эргодический джойнинг   тройки $(T,T,T)$.
Обозначим $\eta=\pi_{23}\nu$. Если $\eta=\mu\otimes\mu$, из определения
квазипростоты порядка 3 получим
$\nu = \mu\otimes\mu\otimes\mu$.

Теперь рассмотрим случай, когда $\eta$ является квазидиагональной мерой.
Так как условная мера  $\eta_y$ дискретна, получим, что для
почти всех относительно меры $\mu\otimes\mu$  пар $(x,y)$
условнная мера
$\nu_{xy}$ является дискретной, следовательно, она является компонентой
меры $\eta_{y}$. Отсюда вытекает, что мера  $\nu$ абсолютно непрерывна
относительно меры     $\mu\otimes\eta$.

Но мера  $\mu\otimes\eta$  эргодична относительно
$T\times(T\times T)$, так как
произведение  слабо перемешивающей системы $(T,\mu)$ и эргодической
системы $(T\times T,\eta)$ эргодично    относительно меры
$\mu\otimes\eta$.
Таким образом, $\nu=\mu\otimes\eta$.
\\  
\medskip

{\bf ТЕОРЕМА 2.2.2.} \it
HI-система является тензорно простой.
\medskip

{\bf СЛЕДСТВИЕ.} \it
Квазипростота порядка 3 эквивалентна квазипростоте всех порядков.
\rm
\medskip

 Пусть $\nu$ -- попарно независимый джойнинг набора $(T,T,S)$,
где $S$ -- некоторый эргодический
автоморфизм, а  $T$ является  автоморфизмом с  HI-свойством.

Для доказательства  теоремы 2.2.2 мы воспользуемся индуцированными
джойнингами: с джойнингом $\nu$ мы ассоциируем некоторую последовательность
новых джойнингов  $\eta_m$. Дадим определение последних.

Пусть семейство марковских операторов $\{P_x\}$ отвечает мере
$\nu$ :
для всех $A,B,C\in\B$
$$
  \int \! \chi_A(x)\langle P_x\chi_B ,\chi_C \rangle d\mu(x) \
= \ \nu (A\times B\times C),
$$
где $\chi_A $ -- индикатор  множества $A$,
$\langle\ , \ \rangle$ обозначает скалярное произведение
в пространстве $L_2(X,\mu)$.
Из инвариантности  меры  $\nu$ относительно $T\times T\times T$
вытекает тождество
$$S^{-1}P_{T^{-1}(x)}T \equiv P_x .$$

 Рассмотрим новые  меры (джойнинги, индуцированные мерой $\nu$),
обозначаемые $\eta_m,$ $ m\in Z,$ и заданные  равенствами
$$
    \eta_m (A\times B\times C) = \int \! \chi_A(x)
\langle P_x\chi_B ,P_{T^{m}x}\chi_C \rangle d\mu(x).
$$
 Мера $\eta_m$  инвариантна относительно
$T\times T\times T$, что оправдывает название ``индуцированный джойнинг''.
\medskip

{\bf ЛЕММА 2.2.3. }\it
Если $T$ является HI-системой, то равенство
$$
P_{T^{m}x}^\ast P_x = \int_X P_{T^{m}x}^\ast  P_x d\mu(x)
$$
выполнено для почти всех  $x$.
\rm
\medskip

{ Доказательство}.
Фиксируем $m$, пусть $\eta =\eta_m$. Из определения
$\eta$ имеем
$$
\eta (B\times C\times X) = \int \chi_B d\mu\int \chi_C d\mu
=\mu(B)\mu(C)= \eta (B\times X\times C),
$$
следовательно,
$$
    \pi_{12}\eta =\mu\otimes\mu, \  \ \pi_{13}\eta =\mu\otimes\mu.
$$
Так как $T$  есть HI-система, мы получаем
$$
\eta = \mu\otimes\ \pi_{23}\eta,
$$
что влечет за собой выполнеие  равенств
$$\eta_{x}=\pi_{23}\eta$$
для почти всех условных  мер $\eta_{x}$ (соответствующих
 операторам $P_{T^{m}x}^\ast P_x$). Поэтому
для почти всех $x\in X$ выполнено
$$P_{T^{m}x}^\ast P_x = \int_X P_{T^{m}x}^\ast  P_x d\mu(x).$$
\\  

 Доказательство теоремы 2.2.2.     Последовательность $m_i\to\infty$ будем
называть  перемешивающей (конечно, правильнее называть перемешивающей
последовательность  $T^{m_i}$), если
$$
\forall A,B\in \B \ \ \ \mu(A\cap T^{m_i}B) \to  \mu(A)\mu(B).
$$
Существование  перемешивающей последовательности
вытекает из свойства  слабого  перемешивания ( и может служить
определением слабого перемешивания).
\medskip

\bf ЛЕММА  2.2.4. \it
Для перемешивающей последовательности $\{m_i\}$ имеет  место
следующая слабая операторная сходимость:
$$
\int_X P_{T^{m_i}(x)}^\ast P_x  d\mu(x) \to \Theta,
$$
\rm (\it $\Theta$ -- ортопроекция на пространство констант\rm).
\rm
\medskip

{ Доказательство.}
Для заданных функций $f,f'\in L_2(\mu)$,
$\ 0\leq f(x),\, f'(x)\leq 1,$ и $\eps >0$
найдется  набор
марковских операторов $\{P_1,\dots , P_N\}$  такой, что
фазовое пространство $X$ можно представить в виде
$X=B\cup\cup_{k=1}^N A_k$, где $\mu(B)<\eps$ и
$
  \forall k \ \forall x,x'\in A_k $
$$\|P_xf-P_kf\|<\eps,
  \ \ \|P_xf'-P_kf'\|<\eps.
$$
Тогда с учетом свойство перемешивания
$$\mu(T^{-m_i}A_k\cap A_l)\to \mu(A_k)\mu(A_l)$$
и равенства $$\int_x P_x d\mu(x)=\Theta$$
  получаем:
$$
\int_X P_{T^{m_i}(x)}f P_x f' d\mu(x) \ \leq  \
\sum_{k,l}\int_{T^{-m_i}A_k\cap A_l} P_{T^{m_i}(x)}f P_x f' d\mu(x)
\   +   \ 2\eps  \ \leq
$$
$$
\leq \sum_{k,l}\int_{T^{-m_i}A_k\cap A_l} P_kf P_l f' d\mu(x)
 \  +    \ 4\eps  \ \leq  $$
$$
\leq \sum_l\int_{A_l}(\sum_k \mu(A_k)P_kf) P_l f' d\mu(x)
\   +   \ 5\eps  \ \leq
$$
$$
\leq \sum_l\int_{A_l}(\int fd\mu) P_l f' d\mu(x)
\   +    \ 6\eps  \ \leq
\int_X fd\mu \int_X f'd\mu
\   +  \   7\eps.
$$
\medskip

\bf ЛЕММА  2.2.5. \it Для почти всех  $y$   найдется  подпоследовательность
$m_i(y)\to\infty$ такая, что имеет место слабая сходимость
$$
 P_{T^{m_i(y)}(y)}^\ast P_y  \to  P_y^\ast P_y.
$$
{ Доказательство.}     \rm
Для любых $\eps >0$ и $f\in L_2(\mu)$
найдутся множества
$A_1, \, A_2,\,\dots, A_p,\dots $ положительной меры
такие, что
$$
\forall p \ \ \  \sup_{x,y\in A_p}  \|P_xf-P_yf\| \leq \eps.
$$
Так как  $\{ m_i\}$ является  перемешивающей последовательностью,
для почти всех $y$ найдется  подпоследовательность
 $\{m_{i'}\}$  (зависящая от $y$)
такая, что $T^{m_{i'}}(y)\in A_p$. Тогда выполнено
$$
         \|P_{T^{m_{i'}}(y)}f-P_yf\| \ \leq \eps.
$$
Используя обычную диагональную процедуру и сепарабельность
пространства $L_2(\mu)$), получим
$$ P_{T^{m_{i''}}(y)}^\ast P_y \to  P_y^\ast P_y.$$
\\ 
Завершим доказательство теоремы 2.2.2.
Так как, ввиду леммы 2.2.4, выполнено
$$
\int_X P_{T^{m_i(y)}(x)}^\ast P_x  d\mu(x) \to \Theta ,
$$
то из лемм 2.2.3 и 2.2.5  вытекает, что для почти всех
$y$ $\ P_y^\ast P_y =\Theta$, $P_y =\Theta$. Следовательно,
$\nu=\mu\otimes\mu\otimes\mu$.
\\  
\medskip
\medskip

\begin{center}
 { \bf  2.3. Минимальные самоприсоединения и  кратная возвращаемость}
\end{center}
\medskip

{\bf Минимальные самоприсоединения.}
Автоморфизм $T$ пространства $(X,\B,\mu)$
 обладает свойством минимальных самоприсоединений
порядка
$n$ ($T\in MSJ(n)$), если
любой эргодический  джойнинг  $n$ копий $T$,
исключая меру $\mu^{\otimes n}=\mu_{(1)}\otimes\dots\otimes\mu_{(n)}$,
обладает следующим свойством:
одна из его проекций  на двумерную грань  в
$X\times\dots\times X$ являеся мерой  $\Delta_{T^i}$ (сдвиг диагональной
меры).
Неформально говоря,  такой автоморфизм $T$ имеет только
очевидные джойнинги.

Для  $\Z$-действий мы  покажем, что  в классе
$MSJ(2)$ проблема Рохлина эквивалентна  открытому
вопросу  терии джойнингов:
\it совпадает ли класс $MSJ(2)$ с классом  $MSJ(3)$? \rm

Хотя имеются некоммутативные контрпримеры, т.е. для групповых действий
возможно, как мы показали в главе 1, несовпадение классов
$MSJ(2)$ и $MSJ(3)$ (и даже $MSJ(3)\neq MSJ(4)$),
случай $\Z$-действий остается нерешенным.
\medskip

{\bf ТЕОРЕМА  2.3.1.} \it  Если $T\in MSJ(2)$ и $T$ перемешивает с кратностью 2, то
автоморфизм $T$ обладает минимальными самоприсоединениями всех порядков и
кратным перемешиванием всех порядков.
\rm
\medskip

Этам теорема является непосредственным следствием более общего утверждения.
\medskip

{\bf ТЕОРЕМА  2.3.2.} \it   Пусть перемешивающий
 автоморфизм $T\in MSJ(2)$ обладает свойством
кратного возвращения:
  для любого множества $A$ положительной меры  и любой
 последовательности  $k(m)\to \infty$,
$|k(m)-m|\to \infty$  для всех больших $m$ выполнено
условие
$$
 \mu(T^{-k(m)}A\cap T^{-m}A\cap A) > 0.
$$
Тогда $T\in MSJ(3)$ и, следовательно, обладает свойством кратного
перемешивания всех порядков.
\rm
\medskip

ЗАМЕЧАНИЕ.  Сформулированная теорема, в частности, утверждает
следующее:
для перемешивающего
автоморфизма $T\in MSJ(2)$ свойство
$$\lim_{m\to\infty} \mu(T^{-k(m)}A\cap T^{-m}A\cap A) > 0$$
для любого $A$, $\mu(A) > 0$
влечет за собой
$$\mu(T^{-k(m)}A\cap T^{-m}A\cap A) \to \mu(A)^3.$$
\medskip

Сейчас мы сформулируем техническое утверждение,
играющее ключевую роль в доказательстве.
\medskip

\it  {\bf УТВЕРЖДЕНИЕ  2.3.3.}
Пусть автоморфизм $T$ принадлежит классу $ MSJ(2)\setminus MSJ(3)$,
тогда найдется число $a>0$,  множество $M\subset \N$
положительной  плотности,   семейство марковских операторов
 $\{\J_x\},\, \J_x:L_2(\mu)\to L_2(\mu)$, отвечающих
некоторому нетривиальному эргодическому
джойнингу со свойством парной независимости,
и  последовательность $k(m)$ такая, что
  $k(m)\to \infty$,   $|k(m)-m|\to \infty$ и
для любых множеств $A',B$ положительной
меры  для некоторого $A\subset A'$, $\mu(A)>0$
неравенство
$$
   \langle \J_{T^{k(m)}(x)} \chi_B |\chi_B \rangle
     \ >\ a\mu(B)
$$
выполнено для всех   $x\in {A\cap T^{-m}A}$ при  $m\in M$.
\rm
\medskip


Утверждение будет доказано позже, а сейчас мы выведем из него
теорему 2.3.2.

{ Доказательство} теоремы 2.3.2.
Пусть $T\in MSJ(3)$ не выполняется. Ввиду утверждения 2.3.3 имеем:
для любого фиксированного $\e >0$ и   множества $B\in \B$  пространство
$X$ представляется как объединение
некоторых дизъюнктных множеств $A_1,A_2, \dots ,$ таких, что   выполнено
$$
\forall j\   \forall x,x'\in A_j \ \ \  \|\J_x\chi_B -\J_x'\chi_B\|<\e,
$$
причем для всех точек $x\in A_j\cap T^{-m}A_j$ имеет место неравенство
$$
\langle \J_{T^{k(m)}(x)}\chi_B |\chi_B \rangle
     \ >\ \frac{a}{2} \mu(B).
$$
Кратное возвращение обеспечивает следующее:
найдется $x\in A_j\cap T^{-m}A_j$ такая, что ${T^{k(m)}(x)}\in A_j$.
Поэтому выполнено
$$
\forall x'\in A_j \ \ \
\langle  \J_{x'}\chi_B \,|\,\chi_B \rangle\ \ > \ \frac{a}{2} \mu(B) - \eps.
$$
Таким образом, при $\frac{a}{10}\mu(B)>\eps>0$  и
$\mu(B)< \frac{a}{2}$ мы получаем  противоречие:
$$
 \mu(B)^2 = \int_X \langle  \J_{x'}\chi_B \,|\,\chi_B \rangle d\mu(x') =
 \langle  \Theta\chi_B\, |\,\chi_B \rangle
     \ \geq\ \frac{a}{2}\mu(B).
$$
Следовательно, предположение $T\in MSJ(2)\setminus MSJ(3)$ неверно.

Теперь приступим к доказательству утверждения 2.3.3.

Основная  идея состоит в следующем.
Индуцированный  джойнинг $\eta_m$
(индуцированный некоторым эргодическим джойнингом $\nu$ с попарной
независимостью) имеет следующее представление:
     $$
       \eta_m = \frac{1}{p}( \nu_{i(m,1)}+\nu_{i(m,2)}+\dots +\nu_{i(m,p)}),
     $$
где $\nu_{i(m,j)}$ --  эргодические джойнинги класса $M(2,3)$
(конечно, при  $m\neq 0$).
Можно доказать, что одна из компонент, скажем, $\nu_{i(m,1)}$,
стремится к $\Theta$. Это влечет за собой тривиальность исходного
джойнинга  $\nu$.
\medskip

{\bf ЛЕММА  2.3.4. }\it    Пусть $\nu$ -- эргодический джойнинг набора
$(T,T,T)$ и выполнены условия: $\nu\in M(2,3)$,
$\nu\neq \mu\otimes\mu\otimes\mu$
 и $T\in MSJ(2)$.
Пусть $\{\P_x\}$ -- марковский  оператор, отвечающий   мере $\nu_x$,
где $\{\nu_x\}$ ($x\in X_{(1)}$) -- семейство  условных мер
на $X_{(2)}\times X_{(3)}$, отвечающих
джойнингу $\nu$.

Тогда найдутся целые числа $p,q\geq 1$ такие, что

(i) для почти всех $x$ выполнено
$$\P_x^\ast\P_x \ \geq \ \frac{1}{q}I; $$

(ii)
для всех $m\in\N, \ m\neq 0$ имеет  место тождество
     $$
       \P_{T^m(x)}^\ast \P_x  = \frac{1}{p}
          ( \J_x^{i(m,1)}+\J_x^{i(m,2)}+\dots +\J_x^{i(m,p)});
     $$

(iii)
Пусть $B$ --  множество  положительной  меры.
 Для любого множества $A'$, $\mu(A')>0$  найдется множество $A\subset A'$
положительной меры  такое, что для всех $x\in A\cap T^{-m}A$
выполнено

$$
  \langle \P_{T^m(x)}^\ast \P_x \chi_B |\chi_B \rangle
     \ >\ \frac{1}{2q} \mu(B);
$$

(iv)
   для любого $m$ найдется  число $r(m), \ 1\leq r(m)\leq p$, такое, что
для всех $x\in A\cap T^{-m}A$
$$
  \langle \J_x^{i(m,r(m))}\chi_B |\chi_B \rangle
     \ >\ \frac{1}{2q} \mu(B);
$$

(v)     для некоторого множества $M\subset \N$  положительной плотности
для  эквивариантного семейства
        $\{\J_x\}$, соответствующего некоторому джойнингу класса $M(2,3)$,
        выполнено
$$
 \forall m\in M \ \ \
   \J_x^{i(m,r(m))} =\J_{T^{k(m)}(x)}.
$$
\rm
\medskip

{ Доказательство} пункта (i).

 Равенство
  $$\int \P_x^\ast\P_x d\mu(x)= \Theta$$
влечет за собой
$$\P_x^\ast\P_x \equiv \Theta,  \ \ \P_x \equiv \Theta, \ \
\nu=\mu\otimes\mu\otimes\mu. $$
Так как $T\in MSJ(2)$, то оператор   $\int \P_x^\ast\P_x d\mu(x)$ является
выпуклой суммой оператора $\Theta$ и степеней $T^i$.

Так как $\nu\neq\mu\otimes\mu\otimes\mu$,
 для некоторого целого $m$ и числа $a>0$ имеем
$$\int \P_x^\ast\P_x d\mu(x) = aT^i +\dots .$$
Заметим, что случай $i\neq 0$ невозможен. Действительно,
из $\int \P_x^\ast\P_x d\mu(x) \geq aT^i$ вытекает, что
для почти всех $x$ операторы $\P_x$ и $\P_xT^i$
имеют ``общую часть'', т.е.  мера
$\nu$ и мера $(Id\times Id\times T^i)\nu$ имеют
общую компоненту. Но эти меры эргодичны относительно
$T\times T\times T$, следовательно, они совпадают.
Таким образом, мы получили равенство
$$\nu=(Id\times Id\times T^i)\nu. $$
Пусть $i\neq 0$.  Так как $\nu$ принадлежит $M(2,3)$, а
 $T^i$ -- эргодический автоморфизм (автоморфизм класса $MSJ(2)$ обязан быть
слабо перемешивающим), мы получаем $\nu =\mu\otimes\mu\otimes\mu$
(см. принцип дополнительно симметрии).
Таким образом, возможен только  случай $i=0$.
\medskip

{ Доказательство} пункта (ii). Имеем
$$\P_x^\ast\P_x \geq \frac{1}{q}I;
\ \ \P_x\P_x^\ast  \geq \frac{1}{r}I. $$
Тогда для $ {\cal H}_{mx}= \P_{T^m(x)}^\ast \P_x$ мы получаем  неравенства
$$ {\cal H}_{mx}^\ast {\cal H}_{mx} \geq \frac{1}{qr}I,
\ \  {\cal H}_{mx} {\cal H}_{mx}^\ast \geq \frac{1}{qr}I . $$

Теперь  для фиксированного $m$ рассмотрим  систему
$(T\times T\times T, \eta_m)$, где
$\eta_m$ --  джойнинг,  отвечающий семейству $\{ {\cal H}_{x}\}$.
Эта система для некоторого $s\leq qr$ является  $\Z_{s}$-расширением
системы    $(T\times T, \mu\otimes \mu)$.

Так как  последняя  эргодична,  число эргодических компонент
системы $(T\times T\times T, \eta_m)$ не превосходит числа $qr$.
Все компоненты являются джойнингами класса $M(2,3)$.
 Иначе мы получим
$$\nu=(T^m\times Id\times T^i)\nu,  \ m\neq 0,$$
но это влечет за собой равенство $\nu=\mu \otimes\mu \otimes\mu$,
так как  джойнинг $\nu$ принадлежит классу $M(2,3)$, а автоморфизм
$T^m\times T^i$
эргодичен при $i\neq 0$.
\medskip

{ Доказательство} пункта (iii).

    Зафиксируем некоторые множества $A'$ и $B$  положительной меры, для
которых выполнено условие
$P_x\chi_B\neq const$ при
$x\in A'$. Обозначим $ \hat{B}=\chi_B-\Theta\chi_B$.
 Рассмотрим  множество $A'$   положительной меры  такое, что
для некоторого $c,\ 0< c <\frac{1}{q}$  неравенство
$\|P_x \hat{B}\|>c$ выполнено  для всех $x\in A'$.  Для  $\eps>0$
выберем  $f\in L_2(\mu)$
и   множество $A\subset A'$ положительной меры
такое, что  $\|P_x \hat{B}- f\| < 0.1c\mu(B)$ для всех $x\in A$.
Существование  такой   функции $f$ следует из
 сепарабельности  пространства $L_2(\mu)$.
Получаем
$$
 \forall x,x'\in A\ \ \ \|\P_x\hat{B}- \P_{x'} \hat{B}\| < 0.2c\mu(B),
$$
следовательно,
$$
 \forall x,x'\in A\ \ \ \
\|\P_{x'}^\ast\P_x \hat{B}- \P_{x}^\ast\P_x \hat{B}\| < 0.2c\mu(B).
$$
Поскольку выполнено
$$\P_x^\ast\P_x \chi_B\geq \frac{1}{q}\chi_B, $$
при $c <\frac{1}{q}$ мы получаем
$         \forall\  x\in A\cap T^{-m}A$
$$
  \langle \P_{T^m(x)}^\ast \P_x \chi_B |\chi_B \rangle
   \ >\ \frac{1}{q} \mu(B) - 0.2c\mu(B)  \ >\ \frac{1}{2q} \mu(B).
$$
\medskip

{ Доказательство} пункта (iv).

Из пункта (ii) вытекает представление
     $$
       \P_{T^m(x)}^\ast \P_x  = \frac{1}{p}
          ( \J_x^{i(m,1)}+\J_x^{i(m,2)}+\dots +\J_x^{i(m,p)}).
     $$
С точностью до перестановки членов в приведенной выше сумме
для всех $x\in A\cap T^{-m}A$ выполнено  неравенство
$$
  \langle  \J_x^{i(m,1)}\chi_B |\chi_B \rangle
     \ >\ \frac{1}{2q} \mu(B).
$$
\medskip

{ Доказательство} пункта (v).

Теперь   рассмотрим  оператор
$$\jj_m:L_2(X,\mu)\to L_2(X,\mu)\otimes L_2(X,\mu),$$ определенный
формулой
$$
  \langle \jj_m\chi_A |\chi_B\otimes\chi_C \rangle  =
  \int_A\langle \J_x^{i(m,1)}\chi_B |\chi_C \rangle d\mu(x).
$$
Для различных $m,k$ выполнено: или  пространства $\jj_m L^2_0$, $\jj_k L^2_0$
совпадают,
или  эти  пространства ортогональны (см. доказательство
теоремы 2.1.2).

Наша задача -- доказать, что множество попарно ортогональных пространств,
взятых из набора $\{\jj_m L^2_0\}$, должно быть конечным.

Положим

$\bar{f}_m=\chi_{A\cap T^{-m}A} - \mu(A\cap T^{-m}A)\1$,

$\bar{\chi}_B=\chi_B -\mu(B)\1$.
\\
Так как выполнено
$$
     \langle \jj_m\bar{f}_m |\1\otimes\bar{\chi}_B \rangle \ = \
     \langle \jj_m\bar{f}_m |\bar{\chi}_B\otimes\1 \rangle \ = \
    \langle \1 |\bar{\chi}_B\otimes\bar{\chi}_B \rangle  \ =\  0,
$$
получаем
$$
    \langle \jj_m\bar{f}_m |\bar{\chi}_B\otimes\bar{\chi}_B \rangle  =
  \langle \jj_m\chi_{A\cap T^{-m}A} |\chi_B\otimes\chi_B \rangle -
  \mu(A\cap T^{-m}A)\mu(B)\mu(B).
$$
Таким образом, для всех больших $m$ выполнено неравенство
$$
    \langle \jj_m\bar{f}_m |\bar{\chi}_B\otimes\bar{\chi}_B \rangle  \ > \
\frac{1}{2q}\mu(A)^2\mu(B) - \mu(A)^2\mu(B)^2 >0.
$$
Теперь предположим, что для бесконечного множества $N$ для  всех $m\in N$
 пространства $\{\jj_m L^2_0\}$ попарно ортогональны.
Покажем, что предположение приводит к противоречию.

Действительно,
функции $\jj_m\bar{f}_m$ попарно ортогональны, причем
для некоторого положительного числа $c_1$
выполнено   $$\|\jj_m\bar{f}_m\|> c_1.$$
Но из этого вытекает следующее: для некоторой положительной константы $c_2$
для всех $m\in N$ выполнено неравенство
$$
    |\langle \jj_m\bar{f}_m |\bar{\chi}_B\otimes\bar{\chi}_B \rangle|  \ > c_2,
$$
что невозможно, так как из-за попарной ортогональности $\jj_m\bar{f}_m$
получается, что
$$\|\bar{\chi}_B\otimes\bar{\chi}_B\| =\infty.$$

Следовательно, найдется  множество $M$ положительной плотности
(не меньшей, чем число $\frac{1}{|N|}$)
такое, что  все пространства $\{\jj_m L^2_0\}$ совпадают при $m\in M$.
Отсюда  следует, что  для некоторого $\jj=\jj_{m_0}$
$$
\forall m\in M \ \ \    \jj_m =\jj T^{k(m)}.
$$
Мы доказали, что  выполнено
$$
 \forall m\in M \ \ \ \ \ \
            \P_{T^m(x)}^\ast \P_x  = \frac{1}{p}(\J_{T^{k(m)}(x)}+\dots),
$$
причем  для всех $x\in A\cap T^{-m}A$
$$
  \langle  \J_{T^{k(m)}(x)}\chi_B |\chi_B \rangle
     \ >\ \frac{1}{2q} \mu(B).
$$
\\  

Чтобы завершить  доказательство утверждения 2.3.3, покажем необходимость
следующих условий:
 $$k(m)\to \infty, \ \ |k(m)-m|\to \infty.$$
Если для  бесконечного множества, элементы которого обозначим через $m'$,
выполнено  $k(m')=s$, где $s$ фиксировано, то получим
$$
       \Theta = \lim_{N\to\infty} \frac{1}{N}\sum_{i=1}^{N}
    \int \P_{T^{m'(i)}(x)}^\ast \P_x d\mu(x) =
       \frac{1}{p}(\J_{T^{s}(x)}+\dots).
$$
Следовательно, выполнено
$$\J_{T^{s}(x)}\equiv \Theta , \ \J_x \equiv \Theta.$$

Если   для  бесконечного множества различных $m'$
выполнено $k(m')-m'=s$, то
$$
       \Theta = \lim_{N\to\infty} \frac{1}{N}\sum_{i=1}^{N}
    \int \P_x^\ast \P_{T^{-m'(i)}(x)}d\mu(x)  = \frac{1}{p}(\J_{T^{s}(x)}+\dots).
$$
Получаем, что $\J_x \equiv \Theta$. Получили,  что  меры
$\nu_{T^{m}(x)}^\ast \nu_x$ (меры, отвечающие операторам
$\P_{T^{m'(i)}(x)}^\ast \P_x$) имеют в качестве компоненты
меру  $\mu\otimes\mu$.
Это противоречит тому,
что  меры  $\nu_{x'}^\ast \nu_x $ являются квазидиагональными мерами
(напомним, что композиция  квазидиагоналей является квазидиагональю).
Таким образом, $k(m)\to \infty$  и $|k(m)-m|\to \infty$.
  
\medskip
\medskip

\begin{center}
{\bf  2.4. Четная и нечетная тензорная простота}
\end{center}
\medskip
В этом параграфе приводится пример некоммутативного действия $\Psi$,
обладающего свойством  $\s_{2q+1}$, но не обладающего
свойством $\s_{2p}$.
В этом случае говорим, что действие $\Psi$ обладает \it нечетной тензорной
простотой. \rm
Нечетная тензорная простота означает, что
${\Psi}$ не допускает нетривиальных марковских сплетений с четными
тензорными степенями системы ${\Psi}$.

Определим действие $\Psi$.
Предварительно отметим, что автоморфизмы
группы $X= \Z_2\times\Z_2\times\Z_2\dots$  и
сдвиги на $X$ сохраняют меру Хаара $\mu$
на $X$.
Для  перестановки $\sigma$ натурального ряда $\N$ определим автоморфизм
$T_\sigma$
группы $X$ равенством $T_{\sigma}(\{x_i\})=\{x_{\sigma(i)}\}$.
Последовательности $\alpha=\{a_i\}\in X$ сопоставим
 сдвиг: $S_{\alpha}(x)=x+\alpha$, где $+$ обозначает групповую операцию
в $X$.   В качестве $\Psi$ рассмотрим действие, порожденное
всевозможными $T_\sigma$ и $S_{\alpha}(x)$.
\medskip

\medskip
\bf ТЕОРЕМА 2.4.1. \it При $q\geq 1$  действие   $\Psi$
 обладает свойствами $\s_{ 2q+1}$, но не обладает свойствами
$\s_{ 2q+2}$.
\medskip
  \rm

Доказательство для случая $q=1$.
$\Psi\notin \s_4$, так как мера
$\eta$, определенная формулой
$
\eta(Y_1\times Y_2\times Y_3\times Y) =
$
$$
\mu^3(\{(a,b,c): a\in Y_1, b\in Y_2, c\in Y_3,\,\, a+b+c\in Y\}),
$$
отлична от $\mu^4$,  инвариантна относительно
$\Psi\otimes \Psi\otimes \Psi\otimes \Psi$
и принадлежит  классу $M(3,4)$.

Покажем, что  действие $\Psi$  обладает свойством  $\s_3$.
Пусть мера $\nu\in M(2,3)$ инвариантна относительно
${\Psi}\otimes{\Psi}\otimes{\Psi}$:
$$
\forall \psi\in \Psi \ \ \
\int \psi f\otimes \psi g\otimes \psi hd\nu \
= \int f\otimes g\otimes h d\nu .
$$
Характеры группы $X$ имеют вид  $\chi_F =\prod_{n\in F}\chi_n$,
где $F$ - конечное подмножество натурального ряда, а
$\chi_n(x)=1$,  если $x_n=0$,   и $\chi_n(x)=-1$,  если $x_n=1$.
Отметим свойства действия ${\Psi}$:
$S_{\alpha}\chi_{n}=-\chi_{n}$ при $a_n=1$;
$S_{\alpha}\chi_{n}=\chi_{n}$ при $a_n\neq 1$;
$\chi_{\sigma(n)}(x) = \chi_n(T_{\sigma}(x)).$
Пусть $A,B,C$ -- конечные подмножества $\N$.
Если $A\cap B\cap C \neq \phi$ и
   $j\in A\cap B\cap C$, то
для сдвига $S_j$, определенного формулой
$(S_j(x))_i=x_i +\delta_{ij}$,
имеем
$$
\int  (\chi_A\otimes \chi_B\otimes \chi_C)  d\nu  =
\int(S_j\chi_A\otimes S_j\chi_B\otimes S_j\chi_C)  d\nu =
$$ $$ - \int (\chi_A\otimes \chi_B\otimes \chi_C) d\nu  = 0.
$$
То же самое получим при $j\in A$, когда $j\notin  B\cup C$.

Теперь рассмотрим случай, когда  $j\notin A$ и $j\in  B\cap C$.
Найдется перестановка $\sigma$, удовлетворяющая условиям:
 a)\ $\sigma(i)=i,$ если $i\in A$;
b)\ $\sigma^n(j)\to\infty$ при $n\to\infty$.
Тогда, учитывая инвариантность меры $\nu$ и свойства a),b), получим
$$
\int (\chi_A\otimes \chi_B\otimes \chi_C)  d\nu =
\int  (T_\sigma \chi_A\otimes T_\sigma \chi_B\otimes T_\sigma \chi_C)d\nu  =
$$
$$
 \int  (\chi_A\otimes T_\sigma \chi_B\otimes T_\sigma \chi_C) d\nu  =
{1\over N}\int  (\chi_A\otimes \sum_{n=1}^N
 (T_\sigma^n \chi_B\otimes T_\sigma^n \chi_C)) d\nu.
$$
Последнее выражение стремится к 0, что  вытекает из попарной
ортогональности  функций
$$
Y_n = T_\sigma^n \chi_B\otimes T_\sigma^n \chi_C =
\chi_{\sigma^n(B)}\otimes \chi_{\sigma^n (C)}, \ \ n=0,1,2,\dots\ .
$$
Подводя итог сказанному выше, заключаем, что неравенство
$\int (\chi_A\otimes \chi_B\otimes \chi_C)d\nu \neq 0$
может иметь место только в случае $A=B\sqcup C$ (дизъюнктное объединение).
Однако, по тем же причинам получим $B=C\sqcup A$, следовательно,
$A=(C\sqcup A)\sqcup C= A\sqcup C$, что невозможно. Доказано,
что  $\int  (\chi_A\otimes \chi_B\otimes \chi_C) d\nu= 0$
  для всех наборов
$\chi_A,\chi_B,\chi_C$, следовательно,  $\nu=\mu^3$.  Таким образом,
$\Psi\in \s_3$.  Аналогичные рассуждения показывают, что
$\Psi\in \s_{2q+1}$.

 Действие $\Psi$ содержит много поддействий класса
$\s_{2q+1}\setminus \s_{4}$.
 В частности, к ним относится действие счетной подгруппы,
порожденной всеми финитными   перестановками
($supp(\sigma)$ -- конечное множество)  и финитными сдвигами
(сдвигами, для которых последовательности $\alpha=\{a_i\}$ финитны).

\medskip
 \bf Действия с минимальным и простым централизатором. \rm
Следуя терминологии [3], действие $\Phi$ называем
 2-простым, если любой эргодический джойнинг $(\Phi,\Phi)$ есть мера
$\mu\otimes \mu$ или мера вида $(Id\otimes R)\Delta$, где
$R$ -- обратимое преобразование, коммутирующее с $\Phi$,
а $\Delta$ - диагональная мера на $X\times X$.

Действие $\Phi$ называется $n$-простым, если  $\Phi$ является
 2-простым и принадлежит классу $\cap_{m=3}^{n} \s_m$.

Если все эргодические джойнинги пары $(\Phi,\Phi)$ суть $\mu\otimes \mu$  и
меры вида $(Id\otimes T_h)\Delta$, где $T_h$ входит в действие
$\Phi=\{T_g:g\in B\}$, а $h$ лежит в центре группы $B$, говорят, что
действие $\Phi$ имеет минимальный марковский цетрализатор.
В теории джойнингов принято писать в этом случае $\Phi\in MSJ(2)$.
Обозначим  $MSJ(n)=MSJ(2)\cap\cap _{m=3}^{n} \s_{m}$.

Для $\Z$-действий выполнено $MSJ(3)=MSJ(4)$, для них имеет место
также  более общий факт:
3-простота влечет за собой $n$-простоту.
Однако для некоммутативных динамических систем аналогичное утверждение,
вообще говоря, неверно.

\bf ТЕОРЕМА 2.4.2. \it 1.
Действие $\Psi$,  порожденное всеми  автоморфизмами
группы $X'= \Z_3\times\Z_3\times\Z_3\dots$
 и всеми групповыми сдвигами на $X'$, принадлежит   классу  $MSJ(2)$,
но не принадлежит классу $MSJ(3)$.
Действие $\Phi$,  порожденное всеми автоморфизмами
группы $X= \Z_2\times\Z_2\times\Z_2\dots$
 и всеми сдвигами на группе $X$, принадлежит   классу
$MSJ(3)\setminus
MSJ(4)$.
Эргодическими джойнингами как  $(\Phi,\Phi)$ так и $(\Psi,\Psi)$
являются  только меры $\Delta$ и $\mu\otimes \mu$.

2. Классу $MSJ(3)\setminus MSJ(4)$ принадлежит действие  $\Phi'$,
порожденное (бернуллиевским) автоморфизмом
$T$ на   $\dots\times\Z_2\times\Z_2\times\Z_2\dots$
и инволюциями $Q,R,S$, определенными следующим образом:
$$T(x)_i=x_{i+1},$$
$$Q(\dots x_{-2}, x_{-1}, x_{0}, x_{1},  x_{2},\dots)=
(\dots x_{-2}, x_{-1}, x_{1}, x_{0},  x_{2},\dots),$$
$$R(\dots x_{-2}, x_{-1}, x_{0}, x_{1},  x_{2},\dots)=
(\dots x_{-2}, x_{-1}, x_{0}, x_{1}+ x_{0},  x_{2},\dots),$$
$$S(\dots x_{-2}, x_{-1}, x_{0}, x_{1},  x_{2},\dots)=
(\dots x_{-2}, x_{-1},  x_{0}+1,x_{1}, x_{2},\dots).$$

\bf СЛЕДСТВИЕ. \it
Простота порядка 2, вообще говоря, не влечет простоту
порядка 3.
Простота порядка 3, вообще говоря, не влечет простоту
порядка 4.
\rm
\medskip

Теорема 2.4.2. доказывается аналогично теореме 2.4.1.
Пункт 2 показывает,
что контрпримерами могут служить действия конечно порожденных групп.

Действие  $\Psi$ не обладает свойством $\s_3$, так как мера
$\nu$, определенная формулой
$$
\nu(Y_1\times Y_2\times Y) =
\mu^2(\{(a,b): a\in Y_1, b\in Y_2,  a+b\in Y\}),
$$
является нетривиальным самоприсоединением.
Для действия $\Phi$  нетривиальный  джойнинг $\eta$ был предъявлен ранее
при доказательстве теоремы 2.4.1.

\newpage

\begin{center}   { \Large  \bf  ГЛАВА 3 }
\end{center}
\medskip
\begin{center}
{ \bf ДЖОЙНИНГИ И ТЕНЗОРНАЯ ПРОСТОТА НЕКОТОРЫХ ПОТОКОВ }
\end{center}
\medskip

Поток является коммутативным и непрерывным действием. Изучая
джойнинги потоков, можно рассматривать их малые возмущения.
Это существенно отличает случай потоков от случая каскадов
(действий группы $\Z^n$) и позволяет для потоков получить результаты
более сильные, чем для каскадов.
В этой главе мы установим свойство тензорной простоты для
слабо перемешивающих потоков,
допускающих только гладкие джойнинги, и для потоков,
марковский централизатор которых порожден счетным
набором операторов из полугруппы ${\cal D}$.
Для перемешивающих потоков положительного локального ранга
установим свойство тензорной простоты и свойство кратного
перемешивания.
\medskip


\begin{center}
{\bf  3.1. Гладкие джойнинги и  внутренние сплетения потоков            }
\end{center}
\medskip
М.Ратнер \cite{Rat} доказала, что любой эргодический
джойнинг  унипотентного потока является гладким: он сосредоточен
на гладком подмногообразии $Y\subset X^n$ (декартова степень $X$) и абсолютно
непрерывен относительно меры Лебега на  $Y$.
Назовем потоки с таким свойством S-потоками.

В \cite{Rat} доказано больше:
эргодическим   джойнингом  унипотентного потока является $H$-инвариантная мера,
сосредоточенная
на замкнутой орбите $Hx\subset X^n$ связной подгруппы $H$
в группе $G^{\otimes n}$. А.Старков, используя этот факт,
решил проблему о кратном перемешивании
для однородных потоков \cite{Star}. Доказательство в \cite{Star} можно
модифицировать: после того, как проведена
редукция к случаю унипотентных потоков, достаточно
доказать следущее утверждение.
\medskip

{\bf ТЕОРЕМА 3.1.1. } \it Слабо перемешивающий S-поток
является тензорно простым.

{\bf СЛЕДСТВИЕ.}
Унипотентный перемешивающий поток является тензорно простым
и перемешивающим всех степеней.  \rm
\medskip

Доказательство. Как было показано в  1.5, свойство  тензорной простоты
эквивалентно   свойству $S(2p-1,2p)$ для любого $p>1$.
Пусть размерность пространства $X$ равна $d$, и на нем действует S-поток $\{T_t\}$.
Покажем, что поток принадлежит
классу $S(n,n+1)$ при $n>d$ (из этого
следует тензорная простота потока $\{T_t\}$).

Фиксируем $n>d$. Пусть $\nu$ -- эргодическая относительно
$\{T_t\otimes\dots\otimes T_t\}$ ($n+1$ сомножителей) мера класса
$M(n,n+1)$. Предположим, что $\nu$ сосредоточена на гладком
многообразии  $Y\subset X^{n+1},\, dim(Y)< dim(X^{n+1})$.
Для некоторого фиксированного $\delta >0$ рассмотрим множества
$$
V_r= \bigcup_{0<t_1,t_2,\dots,t_{r}<\delta}
 (T_{t_1}\otimes T_{t_2}\otimes \dots\otimes
T_{t_{r}}\otimes I\otimes I\dots) Y .
$$
Так как множество $V_r$  инвариантно относительно
эргодического потока $T_t^{\otimes n+1}$ (поток $T_t$
слабо перемешивающий относительно $\mu$), мера (объем) множества
$V_r$ равна 0 или 1.  Последнее невозможно, если $\delta$
достаточно мало.

Таким образом, считаем, что $\mu^{n+1}(V_r) = 0$,
следовательно, $dim(V_r) < dim(X^{n+1})$.
Начиная с некоторого $r$ размерности  множеств $V_r$ стабилизируются:
$dim(V_r)= dim(V_{r+1})$. Это означает, что
для достаточно малых $s_1,s_2,\dots, s_{r}, s$, $\,s>0$
$$
(T_{s_1}\otimes T_{s_2}\otimes \dots\otimes
T_{s_{r}}\times T_s\otimes I\otimes I\dots) Y
\subset V_r.
$$
Переходя от многообразий к джойнингам, определим
 меры:
$$
\nu_{(t_1,t_2,\dots,t_{n+1})}(A_1\times A_2 \times\dots\times A_{n+1}) =
\nu(T_{-t_1}A_1\times  \dots\times T_{-t_{n+1}}A_{n+1}),
$$
$$
\overline{\nu} = \int_{0<t_i<\delta}\nu_{(t_1,t_2,\dots,t_{r},0,0,\dots)}
                dt_1\,dt_2\dots dt_r  \ .
$$
При достаточно малом $s> 0$   получим, что мера $\overline{\nu}$ и ее сдвиг
под действием
$I_{1}\otimes \dots \otimes I_{r} \otimes T_s \otimes I_{r+2}
\otimes I_{r+3}\dots \otimes I_{n+1}$
имеют общую часть, так как обе эти меры абсолютно непрерывны относительно
меры Лебега на $V_r$,  а $s$ мало.
Так как разложение этих мер на эргодические
компоненты известно (это свиги эргодической меры $\nu$), получим
$$
\nu_{(s_1,s_2,\dots,s_{r},s,0,0,\dots)} =
   \nu_{(t_1,t_2,\dots,t_{r},0,0,0,\dots)}
$$
(слева и справа выписаны эргодические компоненты).
Таким образом,
$$\nu=(T_{s_1-t_1}\otimes\dots\otimes T_{s_{r}-t_r}\otimes T_s\dots\otimes
I\otimes I\otimes\dots)\nu.
$$
Но из того, что мера  $\nu$ принадлежит  классу $M(n,n+1)$ и
 инвариантна относительно
$(\dots\otimes T_s\otimes\dots\otimes I\otimes \dots)$ ($n+1$ сомножителей),
где $T_s$ -- слабо перемешивающий автоморфизм,
из принципа дополнительной симметрии следует $\nu=\mu^{n+1}$.

Таким образом, для достаточно больших $n$ мера $\mu^{n+1}$ --
единственная
$(T_t\otimes\dots\otimes T_t)$-инвариантная мера класса
$M(n,n+1)$. Следовательно, поток является
тензорно простым.   
\medskip
\medskip

{\bf  3.2. Кратное перемешивание  $\omega$-простых   потоков   }
\medskip

Будем писать $P\geq Q$, если для операторов $P,Q$ для всех
неотрицательных функций  $f,g$ выполненено $Pf\geq Qf$.
\\ ОПРЕДЕЛЕНИЕ. Действие $\{T_g\},\, g\in G$ называется
$\omega$-простым, если
выполнено следующее условие:  бистохастический централизатор действия
содержит счетный набор  операторов $\{P_i\}$ и счетный набор
положительных чисел $a(i)>0$, для которых имеет место
$$
    P_i^\ast P_i > a(i)I,\,\,\,
 P_iP_i^\ast > a(i)I ,     \eqno (3.1)
$$
причем каждый  неразложимый   оператор $J\neq\Theta$
из стохастического централизатора действия $\{T_g\}$
имеет  вид  $J = SP_i$ для некоторого автоморфизма $S$, коммутирующего
с действием.

Свойство (3.1) эквивалентно конечнозначности полиморфизма, отвечающего
оператору $P_i$.  Мера $\theta$, отвечающая такому оператору, является
относительно дискретной. При разложении
$$
\int(f\otimes g)d\theta =
 \int f(x)\left(\int g(y)d\theta_x(y)\right) d\mu(x)
$$
меры $\theta$ в систему условных мер
получим, что почти все условные меры $\theta_x$ являются дискретными
мерами.  Более того, для некоторого $m = m(x)$ условные меры имеют вид
$$
\theta_x= \frac{1}{m}(\delta_{y(1,x)} +\delta_{y(2,x)} +
\dots + \delta_{y(m,x)}),
$$
где $\delta_{y(i,x)}$ -- нормированная точечная мера с носителем в
точке $y(i,x)$.  Из эргодичности потока  вытекает, что почти
все меры $\nu_x$ устроены одинаково в следующем смысле:
они имеют одно и то же число точечных компонент ($m$ не зависит от $x$).
Действительно, пусть $X_n$ обозначает множество
$\{x\in X : m(x)=n\}$, тогда получим $\mu(X_n)=0$ или $\mu(X_n)=1$,
так как $X_n$ инвариантно относительно эргодического действия
$\{T_g\}$.
\\ \bf Пример $\omega$-простого некоммутативного действия. \rm
Рассмотрим двумерный тор $X$ как группу с мерой Хаара $\mu$.
Автоморфизмы двумерного тора образуют
$\omega$-простое действие $\Psi$ группы $GL(2,\Z)$. Это
непосредственно вытекает из результатов
\cite{Par}  о джойнингах второго порядка действия $\Psi$.
Типичным примером неразложимого оператора служит оператор $P$,
заданный условием $(Pf)(2x)=f(5x)$ (для всех функций $f$ на $X$).
Однако, действие $\Psi$ не является  тензорно простым. Действительно,
оператор $J:L_2\otimes L_2\to L_2$, которому отвечает мера,
распределенная  равномерно на множестве
$\{(x,y,z) : y-x=z-x\}$,  является нетривиальным внутренним оператором,
сплетающим $\Psi$ и $\Psi\otimes\Psi$.  Такие операторы можно задавать
явно. Например, оператор $J:L_2\otimes L_2\to L_2$, заданный равенством
$$
(Jf)(y) = \int_X f(py-qx,my-nx)dx
$$
($f$ -- функция из $L_2\otimes L_2$), при подходящем выборе $m,n,p,q\in\Z$
является внутренним.

Отметим, что рассмотренное действие $\Psi$ является дискретным
и некоммутативным.
\\  \bf Коммутативный непрерывный случай (потоки).
 \rm
Следующая теорема дополняет теорему 2.1.2 о том, что
простота порядка 2 для потоков влечет за собой простоту всех порядков.

Напомним, что потоком называется действие группы $\R^n$, являющееся
непрерывным: операторы $T_r$ слабо стремятся к тождественному оператору
$I$, если $r\to 0$.
\medskip

{\bf ТЕОРЕМА 3.2.1.} \it  Слабо перемешивающий $\omega$-простой поток
является тензорно простым.
\rm
\medskip

{\bf СЛЕДСТВИЕ.} \it  Перемешивающий $\omega$-простой поток является
перемешивающим всех порядков.
\medskip    \rm

{\bf ЛЕММА 3.2.2.} \it  Пусть $J$ и $J_s,\, s\in (0,\eps)$,  суть
неразложимые операторы, сплетающие $\omega$-простой поток $T_t$
с некоторым потоком $S_t$. Если выполнено $J\neq \Theta$ и $J_s\to J$
при $s\to 0$, то для некоторых различных $s,s'\in (0,\eps)$
найдется такой автоморфизм $R$, коммутирующий с $T_t$,
что выполнено равенство $J_{s'}=J_{s}R.$
 \rm
\medskip

Доказательство леммы.
 Пусть  $X'= X$, поток $\{T_t\}$ действует на пространстве $(X,\mu)$, а поток
$\{S_t\}$  -- на  $(Y,\lambda)$. Рассмотрим
меру $\eta_s$ на $X\otimes X'\otimes Y$, определенную следующим образом:
$$
\eta_s(A\times B\times C) =
\int_Y (J\chi_A)\,(J_{s}\chi_B)\,\chi_C\, d\lambda ,
$$
где операторы $J,J_{s}$ действуют из пространства
$L_2(\mu)$ в пространство $L_2(\lambda)$.
Операторам $J_s$ соответствуют меры  $\nu_s$ на $X'\times Y$, связь
между ними задается формулой
$$
  \nu_s(B\times C)=\langle J_s\chi_A\,,\,\chi_B \rangle .
$$
Наша  цель -- доказать, что для некоторых различных
$s',s$  и некоторого автоморфизма $R$, коммутирующего с потоком,
выполнено
$$
\nu_{s'}=(R\times Id)\nu_s,                     \eqno               (3.2)
$$
откуда непосредственно вытекает утверждение леммы.

Проекции $\pi\eta_s$ и $\pi'\eta_s$ меры
$\eta_s$ соответственно на $X\times Y$ и $X'\times Y$
суть меры  $\nu$ и  $\nu_s$, отвечающие сплетениям $J$ и $J_{s}$.
Так как эти операторы неразложимы,
меры $\nu$ и  $\nu_s$ эргодичны
относительно $T_t\otimes S_t$. Следовательно, почти все
эргодические  компоненты меры $\nu_s$ обладают такими же проекциями
(те меры, которые не попали в "почти все", мы не рассматриваем).

Так как бистохастический оператор $J_s^\ast J$ коммутирует с потоком,
из определения $\omega$-простоты потока получим
$$
 J_s^\ast J = a_s\left(\sum_i J_i(\int_{SC} Rd\sigma_{(s,i)}(R))\right)
 +(1-a_s)\Theta,
$$
где $\sigma_{(s,i)}$ -- некоторые положительные,
не обязательно нормированные
меры на   стохастическом централизаторе $SC$ нашего потока.
Для всех достаточно малых $s$ выполнено $a_s\neq 0$,
так как $J_s^\ast J\to J^\ast J\neq \Theta$.
Напомним, что неразложимому оператору $\Theta$ отвечает мера $\mu\otimes\mu$.

Заметим, что проекция $\pi_{12}\eta_s$ меры $\eta_s$ на $X\times X'$
отвечает оператору $J_s^\ast J$, следовательно,
$\pi_{12}\eta_s$ имеет вид
$$
\pi_{12}\eta_s = a_s\left(\sum_i \int_{SC} (Id\times R)\theta_i
d\sigma_{(s,i)}(R)\right) + (1-a_s)\mu\otimes\mu,
$$
где  $\theta_{i}$ -- полиморфизмы, отвечающие
операторам $P_i$, которые участвуют в определении $\omega$-простоты  потока.
Таким образом, для достаточно малых $s$
инвариантные относительно $\{T_t\times T_t\}$ меры $\pi_{12}\eta_s$
содержат эргодические компоненты вида
$$(Id\times R_s)\theta_{i(s)}.          \eqno     (3.3)    $$
Считаем, что сказанное выше выполнено при всех $s\in (0,\eps)$
(если нет, то желаемое достигается при уменьшении $\eps$).
Сопоставим каждому $s$ эргодическую компоненту $\beta_s$
меры  $\eta_s$, потребовав, чтобы $\beta_s$ имела вид (3.3).
Так как неэквивалентных полиморфизмов
$\theta_i$ только счетное число,  найдется
несчетное множество $W\subset (0,\eps)$
такое, что для всех $w\in W$ для некоторого фиксированного $n$
выполнено $i(s)=n$. Тогда   для всех $w\in W$ имеем
$$\pi_{12}\beta_w =(Id\times R_w)\theta_n .$$
Из относительной дискретности меры $\theta_s$ вытекает
относительная дискретность меры $\beta_s$, т.е. мера $\beta_s$
имеет следующее представление. Пусть $Z={X\times Y}$, тогда
$$
\beta(F\otimes f) =
\int_Z \left(\int_{X'} f(x')d\beta_{s,z}(x')\right) F(z) d\pi_{12}\eta(z) ,
$$
где $z\in Z$, а условные меры $\beta_{s,z}$ являются дискретными мерами.

Вспомним, что  $\pi_{13}\eta_s =\nu_s$, и заметим, что из эргодичности
(неразложимости) меры $\beta_s$ вытекает, что
для $\nu$-почти всех $z$ условные меры $\beta_{s,z}$
имеют вид
$$\beta_{s,z}= \frac{1}{p}(\delta_{x'(1,z)} + \dots + \delta_{x'(p,z)}),$$
где $\delta_{x'(i,z)}$ обозначает точечную меру на $X'$, т.е.
$\delta_{x'(i,z)}(\{x'(k)\})=1$ (подразумевается, что $x'(i,z)$ зависят от
$s$).

Пусть $\Gamma$ -- график полиморфизма $\theta=\theta_n$. Рассмотрим
меру $\tilde{\theta}$ на  $X\times X'\times Y$, которая
является поднятием меры $\nu$ на множество $\Gamma\times Y$.
Формально мера $\tilde{\theta}$ определяется следующим образом:
$$\tilde{\theta}(F\otimes f) =
\int_Z \left(\int_{X'}f(x')d\tilde\theta_{(x,y)}(x')\right)  F(x,y) d\nu ,
$$
где по определению $\tilde\theta_{(x,y)}=\theta_{(x)}$, а
$\theta_{(x)}$  есть условная мера, возникающая при разложении меры
$\theta$.

Заметим, что мера $\beta_s$ абсолютно непрерывна относительно
меры $\tilde{\theta}_s =(Id\times R_s^{-1}\times Id)\tilde{\theta}$.
Это есть следствие следующих фактов:

 1) проекции мер на $X\times Y$    совпадают;

 2) проекции мер на $X\times X'$ совпадают;

 3) проекции мер на $X\times X'$ относительно дискретны.
\\ Если $\theta$ является мерой, отвечающей
(однозначному) автоморфизму,
то $\beta_s$ и $\tilde{\theta}_s $ совпадают. Если же
$\theta$ отвечает конечнозначному полиморфизму,
график последнего можно представить в виде конечного объединения
графиков $\Gamma_j$ обычных отображений.
Мера $\beta_s$ есть сумма дизъюнктных мер $\beta_{s,j}$ таких, что
$\Gamma_j$ является носителем меры  $\pi_{12}\beta_{s,j}$.
Очевидно, что $\pi$-проекция меры  $\beta_{s,j}$ на $X\times Y$
абсолютно непрерывна относительно меры $\nu$. "Подняв" меру $\nu$ на график
$\cup_j \Gamma_j\times Y$, и получив тем самым  меру $\tilde{\theta}_s,$
замечаем, что меры $\beta_{s,j}$ абсолютно непрерывны относительно
$\tilde{\theta}_s.$

Вспоминая определение множества $W$, замечаем, что все меры
$\beta_w$ при  $w\in W$ абсолютно непрерывны относительно
$(Id\times \tilde R_w\times Id )\tilde{\theta}$, где $\tilde R_w =R_w^{-1}$.
Так как таких мер несчетное число, то среди них найдутся две
недизъюнктные меры.  (На самом деле, в наборе из  $p+1$
таких мер найдется недизъюнктная пара.)
В силу эргодичности эти меры совпадут.

Подведем итог. Мы доказали, что для некоторых различных $s,s'\in W$
мера $(\tilde Id\times R_s\times Id)\beta_s$
равна мере $(\tilde Id\times R_{s'}\times Id)\beta_{s'}$.
Тогда равны их проекции на $X\times Y$:
 $$(\tilde Id\times R_s)\nu_s = (Id\times \tilde R_{s'})\nu_{s'}.$$
Отсюда непосредственно вытекает (3.2), что завершает доказательство леммы.
\\ 
\medskip

Доказательство теоремы 3.2.1. Пусть внутренний оператор
$J:L_2(\mu)\to L_2(\mu^3)$ есть неразложимое сплетение потока
с его тензорным кубом.
Положим
$$
 J_s = (T_s\otimes I\otimes I)J ,
$$
тогда выполнены условвия леммы 5.2. Следовательно, для некоторого
$s\neq 0$ выполнено
$$
(T_{s}\otimes \tilde{I})JR = J,
$$
где $\tilde{I}=I\otimes I$ -- тождественный оператор на пространстве
$(\tilde{X},\tilde{\mu})=(X\times X,\mu\otimes \mu)$.
Учитывая, что  $T_{s}$ -- слабо перемешивающий оператор
и  что оператор $J$ внутренний, по лемме 2.2
получим, что $J$ -- тривиальный оператор.  
\\ ЗАМЕЧАНИЕ.
Выше мы рассматривали оператор вида
$J^\ast (I\otimes T_s)J = P$, где $J$ --- гипотетический нетривиальный
оператор, сплетающий поток $T_t$  с потоком $T_t\otimes T_t$.
Так как оператор $P$ коммутирует c потоком, мы располагали информацией
о виде оператора $P$, затем делали выводы о том, каким может быть $J$.
Если же вместо потока рассмотреть
некоммутативную динамическую систему,
то оператор $P$, вообще говоря,
не коммутирует со всей системой.
Это препятствует
перенесению результатов
на некоммутативный случай.

{\bf  3.3. Тензорная простота  потоков положительного локального  ранга }
\medskip

{\bf Локальный ранг.}
Говорят, что автоморфизм $T$ обладает \it локальным рангом \rm
$\beta > 0$, если
для некоторой последовательности $U_j=\bigsqcup_{k\in Q_j}T^zB_j$ башен
Рохлина-Халмоша, где $Q_j=\{0,1,\dots,h_j\}$, выполнены условия:

$\mu(U_j)\to \beta$, и для каждого $A\in\B$ пересечение  $U_j\cap A$
асимптотически близко к
объединению некоторго набора этажей в башне. Это означает, что
для некоторой последовательности множеств $S_j\subset   \{0,1,\dots, h(j)-1\}$
(последовательность зависит от $A$) выполнено
$$\mu((U_j\cap A)\Delta\bigsqcup_{z\in S_j}T^zB_j)\ \to \ 0,  \ j\to\infty .
$$
О системах локального положительного ранга  см. обзор \cite{F}.
Для  $\Z^n$-действий и $\R^n$-действий (потоков) определение аналогично,
конфигурации $Q_j$ имеют    вид
$
Q_j = \{0,1,\dots, h_1(j)\}\times
\dots \times \{0,1,\dots, h_n(j)\}$, причем   $h_m(j)\to\infty$
$(1\le m\le n)$.

В определении локального ранга потока $\{T_r\}$, $r\in\R^n,$
соответствующие башни имеют вид
$U_j=\bigsqcup_{z\in Q_j}T_{t_jz}B_j,$   где $t_j\in \R$,
причем  для всех $m$, $1\le m\le n$, выполнено
$t_jh_m(j) \to \infty$ при  $j \to \infty$ и $t_j\to 0$.
\medskip

ЗАМЕЧАНИЕ. Отметим, что
при линейной замене времени  ранг потока сохраняется, а у
$\Z^n$-действия ранг меняется. Однако
 потоки, полученные линейной заменой времени из
перемешивающего потока $\{T_t\}$,
попарно  неизоморфны в случае $\beta(\{T_t\})>0$.
\medskip

\bf ТЕОРЕМА 3.3.1. \it Перемешивающий поток $\{T_r\},$ $r\in\R^n$, $n \geq 1$,
при $\beta(\{T_r\}) > 0$ обладает свойством перемешивания
всех порядков. \rm

\medskip
Для  доказательства этой теоремы мы  установим,
что рассматриваемый  поток допускает лишь конечное число неэквивалентных
эргодических джойнингов второго порядка. Эквивалентными мы называем
джойнинги $\nu'$ и $\nu$, если для некоторого $v\in \R^n$ выполнено
$\nu'=(I\times T_v)\nu$.
Причем эти джойнинги   являются конечнозначными
полиморфизмами.   В терминологии предыдущего параграфа это означает,
что перемешивающий поток положительного локального ранга
является $\omega$-простым, следовательно, по теореме 3.2.1 является
тензорно простым  и обладает свойством кратного перемешивания.
 В конце параграфа мы приведем также набросок
прямого доказательства теоремы 3.3.1.

Предлагаемые ниже утверждения ( лемма 3.3.2 и теорема 3.3.3)
сформулированы для потоков,
однако они справедливы для $\Z^n$-действий и доказываются анологично.
Мы воспользуемся  этим в следующей главе.
\medskip

\bf ЛЕММА 3.3.2. \it Пусть поток $\{T_r\},$ $r\in\R^n,$
положительного локального ранга обладает свойством перемешивания.
Пусть для джойнинга
$\rho$ двух копий  $\{T_r\}$
 выполнено
$$\limsup_{j\to\infty}\rho(U_j\times U_j) = c>0, \eqno           (3.4)
$$
где $U_j$ --  башни, фигурирующие в определении локального ранга.

Тогда  для некоторого $a\geq 0$
и некоторого джойнинга $\rho'$ имеет место равенство
$$
  \rho=c\left(a\int_{\R^n}(Id\times T_{v})\Delta d\sigma(v)  +
(1-a)\mu\otimes \mu\right) +   (1-c)\rho',                  \eqno (3.5)
$$
где $\sigma$ -- некоторая нормированная мера на  $\R^n$.
\rm
\medskip

ЗАМЕЧАНИЕ. Для  $\Z^n$-действий  имеет место аналогичная лемма, в
формулировке которой $\sigma$ является мерой на  $\Z^n$.
\medskip

Доказательство.  Пусть $K_N$ обозначает куб в $\R^n$ с центром в нуле,
причем линейный размер $K_N$ равен $N$ (для нас важно, чтобы $K_N$
были компактами, а их объединение совпадало с $\R^n$).  Определим множество
$D^N_j\subset (U_j\times U_j)$ как объединение множеств  вида
$T_{t_jq}B_j\times T_{t_jq'}B_j$, для которых выполнено
$(t_jq-t_jq')\in K_N$, $q,q'\in Q_j$.
Пусть для некоторого $N$ величина
$
 d(N)=\limsup_{j\to\infty}\rho(D^N_j)
$
положительна, покажем, что мера $\rho$  имеет представление
$$
\rho = d(N)\int_{K_N} (Id\times T_{v})\Delta d\sigma(v)
 + (1-d(N))\rho'', \ \ \
$$
где $\sigma$ -- некоторая нормированная мера на ${K_N}$,
а $\rho''$ -- некоторый джойнинг.
Обозначим
$C_{w,j} = \bigsqcup_{q}T_{(q+w)t_j}B_j\times T_{t_jq}B_j$ и
определим меру $\Delta^w_j$ равенством
$$
\Delta^w_j(A\times B)= \mu(Y_{w,j}\cap T_{wt_j}A\cap B)/\mu(Y_{w,j}),
$$
где $Y_{w,j}$ есть проекция  множества  $C_{w,j}$
на второй сомножитель произведения $X\times X$.

Для фиксированных $A, B$ значение
$\mu(Y_{w,j}\cap T_{wt_j}A\cap B)/\mu(Y_{w,j})$ близко
к $\mu(T_{t_jw}A\cap B)$ (это вытекает из
эргодичности сдвига диагональной меры,
ограниченности последовательности $t_jw$ и
 свойства
$\mu(Y_{w,j}\Delta T_sY_{w,j}) \to 0$ при
фиксированном $s\in\R^n$).
Суммируя  по  $w$ таким, что $w=q-q'$, $q,q'\in Q_j$ и
$(t_jq-t_jq')\in K_N$,     получим
$$
\rho(\ \ |\,D^N_j) = \sum_{w} \rho(C_{w,j}\ |\ D^N_j) \rho(\ \ |\, C_{w,j})
\approx
\sum_{w} \rho(C_{w,j}\ |\ D^N_j) \Delta^w_j .
$$
Из сказанного получим, что мера $\rho(\ \ |\, D^N_j)$
аппроксимируется суммами   сдвигов  диагональной меры
(подобные аппроксимации  для $\Z$-действий рассматривались в \cite{R0}).
Поэтому предельной точкой для мер  $\rho(\ \ |\, D^N_j)$  будет  интеграл
$
\int_{K_N}(Id\times T_{v})\Delta \, d\sigma_N(v),
$
где мера $\sigma_N$  получается как предельная мера
последовательности
дискретных мер $\sigma_{N,j}$, сосредоточенных
на конечном множестве  точек вида $t_jw\in K_j$,
$w\in \Z^n$. Вес точки $t_jw$ определяется равенством
$$\sigma_{N,j}(\{t_jw\})=
\frac{1}{d(N)}\sum_{w}\rho(C_{w,j}).
$$
Но меры $\rho(\ \ |\, U_j\times U_j)$ также аппроксимируются суммами
сдвигов  диагональной меры,  их пределом будет мера
$$
a\int_{\R^n}(Id\times T_{v})\Delta \, d\sigma(v) +
(1-a)\mu\otimes\mu,
\ \ (\ a = \frac{1}{c}\sup_{N}d(N)\ ), \eqno (3.6) $$
где второе слагаемое возникает из-за свойства
перемешивания потока. Напомним, что оно  эквивалентно сходимости
мер
$(Id\times T_{v})\Delta$ к мере $\mu\otimes\mu$ при  $v\to\infty$.
Теперь с учетом (3.4) и (3.6) получаем представление (3.5).
 
\medskip

\bf ТЕОРЕМА 3.3.3. \it а) Эргодический джойнинг $\nu\neq\mu\otimes \mu$ двух
копий перемешивающего потока $\{T_r\}$, $r\in\R^n$, $n \geq 1,$
при $\beta(\{T_r\}) > 0$ является мерой, сосредоточенной на  графике
конечнозначного отображения.

б)
Для набора таких джойнингов  $\nu_1,\dots,\nu_p$ при
$p\beta(\{T_r\}) > 1$ для некоторых
$i,k$   выполнено  $\nu_i=(I\times T_v)\nu_k$, $1\leq i< j\leq p$.
\rm
\medskip

Доказательство а).   Рассмотрим
оператор $P\,:\,L_2(X,\mu)\to L_2(X,\mu)$, соответствующий джойнингу
$\nu$.
Обозначая  $X'=X''=X$, рассмотрим меру $\eta$ на $X\times X'\times X''$,
   определенную следующим образом:
$$
\eta(A\times A'\times A'') = \int_X \chi_A\,(P\chi_{A'})\,(P\chi_{A''}) d\mu.
$$
Покажем следующее:
если проекция меры  $\eta$ на  $X'\times X''$   имеет компоненту
вида
$$
 \int_{\R^n}(I\times T_v)\Delta d\sigma(v),  \ \ \ \ (\sigma(\R^n)=1),
$$
то $\sigma(\R^n\setminus\{0\})=0$, что равносильно $\sigma(\{0\})=1$.

Пусть $\sigma(\R^n\setminus\{0\})>0$, тогда для некоторого
$v\neq 0$ найдется  эргодическая
компонента $\zeta$(меры   $\eta$), которая   обладает следующими свойствами:

(1) проекции $\zeta$ на $X\times X'$   и  $X\times X''$ совпадают
(рассматривая их как меры  на  $X\times X$, замечаем, что
они равны мере $\nu$);

(2) проекция $\zeta$ на  $X'\times X''$ равна  $(I\times T_v)\Delta$.

Из (2) следует, что  проекция  $\zeta$ на
$X'\times X''$ есть сдвиг на $(Id\times T_v)$ проекции  $\zeta$
на $X\times X''$. С учетом (2) получим       $\nu=(Id\times T_v)\nu$,
что эквивалентно
$P=PT_{v}$.  Если ${v}\neq 0$, то   в силу свойства перемешивания
автоморфизма $T_{v}$ имеем
$$
P=P T_{nv}=P\lim_{n\to\infty} T_{nv}= P\Theta=\Theta.
$$
Мы получили $\nu=\mu\otimes\mu$, что противоречит
предположению о взаимной сингулярности этих мер. Таким образом,
получили, что $\sigma(\R^n\setminus\{0\})=0$. Учитывая также
лемму 2.1   и неравенство
$\int_X P\chi_{U_j} \, \, P\chi_{U_j}\ d\mu   \ \geq \mu(U_j)^2$,
получим
$$
\eta(A\times B)=\int_X  P\chi_A \, \, P\chi_B\ d\mu \geq
$$
$$\geq \lim_j \mu(U_j)^2\Delta(A\times B) =
\beta(\{T_r\}))^2 \int_X \chi_A \ \chi_B\ d\mu .
$$
Следовательно,
$P^\ast P\geq \beta(\{T_r\}))^2 I$, т.е.
мера  $\nu$ сосредоточена   на графике
конечнозначного отображения.

\medskip
Доказательство  б).
Пусть мерам  $\nu_1,\nu_2,\dots,\nu_p$ соответствует набор
марковских операторов  $P_1,P_2,\dots,P_p$. Заметим,
что  неотрицательные функции  $P_i\chi_{U_j}$  не превосходят 1,
а интеграл от каждой функции  $P_i\chi_{U_j}$ равен   $\beta(\{T_r\})$.
В силу  $p\beta(\{T_r\}) > 1$ для некоторого $c>0$ и
некоторых $i,k$  будет выполнено
$$   \limsup_j     \int_X  P_i\chi_{U_j}\,\,P_k\chi_{U_j}\ d\mu  \geq \
     c, \ \ \ 1\leq i< k\leq p.
$$
Рассматривая меру $\eta_{ik}$, отвечающую оператору $P_k^\ast P_i$,
с учетом леммы 2.1 получим
$$
\eta_{ik}(A\times B) = \int_X P_k^\ast P_i\chi_A\,\,\chi_B \ d\mu  =
\int_X P_i\chi_A\,\,P_k\chi_B \ d\mu \geq c.
$$
Мера $\nu_{ik}$ сосредоточена на графике композиции полиморфизмов
$\nu_k^\ast$ и $\nu_i$, следовательно, она сингулярна относительно
$\mu\otimes\mu$.
В силу сказанного из леммы 3.3.2 вытекает, что мера $\nu_{ik}$ имеет
компоненту  вида $(I\times T_v)\Delta$.
Аналогично тому, как в доказательстве п. a)
мы установили $\nu=(I\times T_v)\nu$, теперь мы получим
равенство  $\nu_i=(I\times T_v)\nu_k$.

\medskip
Из теоремы о тензорной простоте $\omega$-простого потока вытекает
теорема 3.3.1.

Приведем набросок  прямого доказательства
тензорной простоты перемешивающего потока положительного локального
ранга.
Пусть  $\lambda$ -- джойнинг третьего порядка с парной независимостью,
т.е.
$$
\forall r\in \R^n \ \ \  \lambda =(T_r\times T_r\times T_r) \lambda\,;
$$
$$
\int_{X\times X\times X}  f\otimes g\otimes h \ d\lambda =
\int f \ d\mu \ \int g \ d\mu \ \int h \ d\mu,
$$
если одна из функций $f,g,h$ является константой.
Нам нужно доказать, что   $\lambda= \mu\otimes\mu\otimes\mu$.

Определим  оператор $J:L_2(\mu)\to L_2(\mu\otimes\mu)$,
отвечающий мере $\lambda$:
$$
\int_{X\times X}  \chi_{C_1}\otimes\chi_{C_2} \ J\chi_B \ d\mu\otimes\mu
= \lambda  (C_1\times C_2\times B).
$$
Для достаточно большого $N$ рассмотрим меру
$\eta$ на $X_{(1)}\times X_{(2)}\times\dots\times X_{(N+2)}$:

\medskip
$
\eta({C_1\times C_2}\times B_1\times\dots\times B_N) \ =
$
$$
\int_{X\times X} \chi_{C_1\times C_2}
(I\otimes T_{\eps})J\chi_{B_1}\,\dots\,(I\otimes T_{k\eps})J\chi_{B_k}
\,\dots\,(I\otimes T_{N\eps})J\chi_{B_N}\, d\mu\otimes\mu .
$$
Если    $\lambda \ne  \mu\otimes\mu\otimes\mu$, то, как  известно,
$J^\ast J\neq\Theta$. Следовательно, обозначая $ J_k=(I\otimes T_{k\eps})J$,
для достаточно малого $\eps\neq 0$  получаем
$J_h^\ast J_k\ne\Theta$ при $h\ne k$. Используя это, лемму 3.3.2 и рассуждения
 из доказательства  теоремы 3.3.3,  показываем, что
для некоторых $h\neq k$ проекция
эргодической компоненты меры  $\eta$  на
$X_{h}\times X_{k}$   является сдвигом диагональной меры
$(I\times T_v)\Delta$, $v\in\R^n$. Это приводит к равенству
$$
(T_{h\eps}\otimes I\otimes I)\lambda =
(T_{k\eps}\otimes I\otimes T_v)\lambda .      \eqno     (3.7)
$$
Следовательно, при  $v\neq 0$ получим
$$\lambda = (T_{(h-k)\eps}\otimes I\otimes T_v)\lambda =
\lim_{p\to\infty}(T_{(h-k)p\eps}\otimes I\otimes T_{pv})\lambda  =
$$
$$
=(\Theta\otimes I\otimes \Theta)\lambda  =  \mu\otimes\mu\otimes\mu
$$
(здесь мы пользуемся тем, что  $T_u\to\Theta$ при $u\to\infty$, так
как поток является перемешивающим).

При  $v=0$ результат не меняется: для любого полиморфизма
$\lambda\in M(2,3)$ выполнено:
 равенство $\lambda =(\Theta\otimes I\otimes I)\lambda$
влечет за собой $\lambda =\mu\otimes\mu\otimes\mu$.
Таким образом, поток обладает свойством тензорной простоты.
\medskip

ЗАМЕЧАНИЕ.
В рассмотренном случае джойнинг  $\lambda$  допускает малые возмущения,
что приводит  к дополнительной симметрии (3.7),
которая решает задачу для потоков.
  В случае $\Z^n$-действий  проекция
меры, аналогичной  мере $\eta$,
не обязана иметь компоненту вида
$(I\times T_v)\Delta$, так как
все  проекции на двумерные грани
могут  совпадать с мерой $\mu\otimes\mu$.
Это  препятствует перенесению  теоремы 3.3.1.  на  $\Z^n$-действия.

\newpage

\begin{center}
{ \bf {\Large \bf  ГЛАВА 4} \\  ДЖОЙНИНГИ ДЕЙСТВИЙ КОНЕЧНОГО \\
    \ \  И ПОЛОЖИТЕЛЬНОГО ЛОКАЛЬНОГО РАНГА.   }
\end{center}
\medskip

  В этой главе вводится понятие  D-свойства для   $\Z^n$-действий и
доказывается, что
перемешивающий автоморфизм конечного ранга
обладает D-свойством.
   Если же локальный ранг перемешивающего $\Z^n$-действия
$\{T_z\}$ больше, чем
$\frac{1}{2^{n}}$, то действие обладает  \it D-свойством. \rm
Эти технические утверждения вместе с теоремой
о тензорной простоте перемешивающих систем с D-свойством приводят
к положительному решению проблемы
Рохлина для автоморфизмов конечного ранга и $\Z^n$-действий
локального ранга, превосходящего значение $\frac{1}{2^{n}}$.
 Показано, что  перемешивающий с кратностью 2
автоморфизм положительного локального ранга является тензорно простым и,
следовательно, перемешивает с любой кратностью.

Также установлено: ранг эргодического автоморфизма
$T\times T$ равен бесконечности,
 локальный ранг  эргодического автоморфизма $T\times T$ не превосходит
$\frac{1}{4}$, а если равенство достигается,
 то $T$ обладает свойством
$ \kappa$-перемешивания при $ \kappa=\frac{1}{2}$.
\medskip

\begin{center}
{\bf  4.1. D-свойство перемешивающих автоморфизмов  конечного  ранга}
\end{center}
\medskip

Говорят, что автоморфизм $S$ имеет {\it ранг}  $r=Rank(S)$,
если $r$ есть минимальное число такое, что найдется
последовательность разбиений $\xi_j\to\eps$  вида
$$ \xi_j=\{  B^1_j,\,    SB^1_j,\,    S^2B^1_j,\,    \dots
,\,    S^{h_j^1}B^1_j,$$
$$ \dots\ \dots\ \dots\  $$
$$  B^r_j,\,  SB^r_j,\,   S^2B^r_j,\,   \dots ,
,\,  S^{h_j^r}B^r_j,\   Y_j  \}$$
( условие $\xi_j\to\eps$
влечет за собой   ${h_j^1}$, ${h_j^2}\dots$, ${h_j^r}\to \infty$
и  $\mu(Y_j)\to 0$).

Будем говорить, что автоморфизм $T$ обладает D-свойством, если
найдутся последовательности аппроксимирующих башен
$ (U_j,\xi_j),\ \  (U'_j,\xi'_j), \ \ (U''_j,\xi''_j),   $
где
$$ \xi_j=\{ E_j, TE_j,\dots T^{h_j}E_j\},$$
$$ \xi'_j=\{ E'_j, TE'_j,\dots T^{h_j}E'_j\},$$
$$ \xi''_j=\{ E''_j, TE''_j,\dots T^{h_j}E''_j\},$$
причем  для некоторой последовательности $\{m_j\}$, $m_j>h_j$,
выполняются следующие условия:
$$\lim_j \mu(U_j) =a>0, \ \  \mu(E_j)=\mu(E'_j)=\mu(E''_j), $$
$$ E'_j= T^{m_j}E_j,\ \  \mu(T^{m_j}U'_j\Delta U''_j)\to 0, $$
$$ \max_{m>h_j}\mu(T^{m}E'_j\ |\ E''_j)\ \to\ 0.$$
\medskip

{\bf ТЕОРЕМА 4.1.1. }\it   Перемешивающий автоморфизм конечного ранга
обладает D-свойством.
         \rm
\medskip

Доказательство. Пусть ранг $T$ равен $n$ и
$$ (U^1_j,\xi^1_j),\ \  (U^2_j,\xi^2_j), \dots, (U^n_j,\xi^n_j)   $$
-- соответствующие последовательности башен при $j\to\infty$.
Будем считать, что
$$    h_j=|\xi^1_j|\leq |\xi^2_j|\leq\dots\leq |\xi^n_j|.$$
Найдется  подпоследовательность $j'$ такая, что существуют пределы
$$   \lim_{j'\to\infty} \frac {|\xi_{j'}^m|} {h_{j'}}, \ \ m=2,\dots,n.$$
(Далее вместо $j'$ пишем $j$.) Положим
$$    r = \max \{ m : \ \frac{|\xi_{j'}^m|}{h_{j'}}<\infty \}.$$
Покажем, что в наборе $\xi^1_j, \xi^2_j, \dots, \xi^r_j$ найдеся
последовательность пар башен $\xi^k_j , \xi^l_j$, из которых
можно  вырезать   требуемые башни $\xi_j, {\xi'}_j,{\xi''}_j$.
Фиксируем $j,k,l$, $1\leq k<l\leq r$, и обозначим через
$E^{kl}_j$   подмножество основания   $E^{k}_j$ башни $U^k_j$, состоящее
из тех точек $x\in E^{k}_j$, которые непосредственно переходят
  из башни $U^k_j$ в
башню $U^l_j$ (не проходя через другие башни):
$$ \exists \ q(x) \ \ T^{q(x)} x \in  E^l_j, \ \  x \in  E^k_j,
$$
и
$$T^i x  \notin  E^m_j, \ \ \forall i \ \ 1\leq i\leq q(x).
$$
На множестве   $E^{kl}_j$ рассматриваем функцию перехода
$q^{kl}_j(x)=q(x)$ и разбиение $\theta^{kl}_j$ на множестве $E^{kl}_j$,
порожденное прообразами функции
$q^{kl}_j(x)$.   Обозначим через
             $D_\e E^{kl}_j$     объединение атомов
 разбиения $\theta^{kl}_j$, меры которых не превосходят
   значения          $\e h^{-1}$.
Покажем, что для некоторых $k,l$ (возможен случай $k=l$)
выполняется
$$         \inf_{\e >0}\{ \limsup_j  \mu(D_\e E^{kl}_j\
|\ E^{k}_j \} =\ d\ >\ 0
$$
( в этом случае говорим, что последовательность  $q^{kl}_j(x)$
обладает D-свойством).

Рассмотрим последовательность  функций возвращения $p_j(x)$ на
множестве $ E^{1}_j$, которые задаются следующим образом:
$$  T^{p_j(x)} x \in  E^1_j, \ \  x \in  E^1_j $$
и
$$T^i x \notin  E^1_j, \ \ \forall i \ \ 1\leq i\leq p_j(x). $$
Отметим, что  "длинные и тонкие" \ башни $U^m_j$ при
$m>r$ асимтотически не   влияют на процесс возвращения множества
$E^1_j$. Поэтому существенным образом башня $U^1_j$ взаимодействует
только с башнями $U^k_j$ при $k\leq r$.
Предположим, что все  последовательности функций перехода (для
$k,l\leq r$)
не  обладают D-свойством.
Возвращение есть композиция переходов, отсутствие D-свойства
сохраняется при композиции, поэтому последовательность
функций возвращения $p_j$ также не обладает D-свойством:
$$ \limsup_j \ \max_{p>h_j}\{ \mu(T^pE^1_j\ |\ E^1_j)\} \ = c\ >0,
$$
т.е.  для некоторой последовательности $p_j>h_j$
$$  \mu(T^{p_j}E^1_j\ |\ E^1_j) \ \to c\ >0.$$
Последнее невозможно, так как противоречит свойству перемешивания.
Действительно, с учетом  $\mu(U^1_j)\to a>0$  получим
для любого измеримого  $A$, $0< \mu(A)< ac$,
$$  \mu(T^{p_j}A \cap A) \ \geq\ ac \mu(A) > \mu(A)\mu(A).$$
Но из перемешивания вытекает  $  \mu(T^{p_j}A\cap A)\to \mu(A)\mu(A).$

Таким образом, для некоторых $k,l$ переход из башни $U^k_j$ в $U^l_j$
обладает D-свойством.
Имеем
$$          \limsup_j  \mu(D_{\e_j} E^{kl}_j\
|\ E^{k}_j)  =\ d\ >\ 0
$$
для   некоторой медленно стремящейся к нулю последовательности
$\e_j$.

Пусть последовательность $h_j$ определена  условием:
$$  h_j = \min \{ \frac{h^k_j}{2},  h^l_j \}.
$$
Положим
$$ E_j =T^{h^k_j -2h_j}D_{\e_j} E^{kl}_j, \
   E'_j =T^{h_j}E_j, $$
а множество $E''_j$ определим    так:
для некоторого малого $\delta > 0$

$$
  E''_j = E^l_j \bigcap \left(\bigcup_{h_j<i<(1+\delta)h_j} T^iE_j\right).
$$
Последовательности башен с основаниями $E_j, E'_j,E''_j$ высоты $h_j$
являются искомыми.
 
\newpage

\begin{center}   \bf
{\bf   4.2.  D-свойство  перемешивающих $\Z^n$-действий и локальный ранг
}
\end{center}
\medskip

\bf Локальный ранг автоморфизмов. \rm
Говорят, что автоморфизм $T$ обладает \it локальным рангом \rm
$\beta > 0$, если
для некоторой последовательности $U_j=\bigsqcup_{k\in Q_j}T^zB_j$ башен
Рохлина-Халмоша, где $Q_j=\{0,1,\dots,h_j\}$, выполнены условия:

$\mu(U_j)\to \beta$, и для каждого $A\in\B$ пересечение  $U_j\cap A$
асимптотически близко к
объединению некоторого набора этажей в башне. Это означает, что
для некоторой последовательности множеств $S_j\subset   \{0,1,\dots, h(j)-1\}$
(последовательность зависит от $A$) выполнено
$$\mu((U_j\cap A)\Delta\bigsqcup_{z\in S_j}T^zB_j)\ \to \ 0,  \ j\to\infty .
$$
\medskip
\bf Локальный ранг $\Z^n$-действий. \rm
Пусть $[0,h]$ обозначает множество $\lbrace 0,1,\dots ,h\rbrace $, а $Q$ есть
куб $[0,h]^n$. Пусть $\xi =\lbrace \xi^z \rbrace_{z\epsilon Q}$ - измеримое
разбиение множества $U\subset X$, т.е. $U = \bigcup_{z\epsilon Q} \xi^z $ и
$\xi^v\bigcap \xi^w =\phi $ при $v\not= w$.
Если такое разбиение $\xi $ с конфигурацией $Q$ удовлетворяет условию
$$
  \forall z\epsilon Q \quad T^z\xi^{0} = \xi^{z}
$$
(здесь $0$ обозначает нулевой вектор в $Z^n$) говорим, что $\xi$ является башней.

Последовательность башен $\xi_j$ называется $ аппроксимирующей $,  если
для любого $A\epsilon \cal B$ выполнено
$$
   |Q_j|^{-1}\sum_{z\epsilon Q_j}  (\mu (A|\xi^{z}_j)^2 -
  \mu (A|\xi^{z}_j)) \quad \rightarrow \quad 0.
$$
(Другими словами: $\mu (A|\xi^{z}_j)\approx 1$ или $\mu (A|\xi^{z}_j)\approx 0$
для большинства $z\epsilon \, Q_j$ при больших $j$).
Последовательность таких троек $(Q_j,\xi_j,U_j)$ называется  аппроксимирующей
последовательностью. Если для действия $\Psi$ найдется аппроксимирующая
последовательность такая, что $\mu (U_j)=b$ для всех $j$, говорим, что действие
$\Psi$ имеет $ локальный\quad ранг\quad b$.
\medskip

\bf D-свойство. \rm
По определению, действие $\{T_z\: z\in \Z^n\}$ положительного локального ранга с
аппроксимрирующей последовательностью $(Q_j,\xi_j,U_j)$, где $\xi_j\to\eps$
и $\mu (U_j)\to \beta >0$,
обладает \it D-свойством, \rm если найдется другая аппроксимирующая
последовательность
$(\tilde{Q}_j,  \tilde{\xi}_j, V_j)$ и последовательность
$\{v(j)\}$, $v(j)\in Z^n,$ такая, что  выполнены
следующие условия:
\medskip

 (i)       $\ \ V_j =
           \bigsqcup_{z\in \tilde{Q}_j}  T_{z}\tilde{B}_j,\ \ \
            \tilde{{Q}_j}\subset Q_j ; $
 \vspace{5mm}

 (ii)      $ \ \ V_j\subset U_j , \ T_{v(j)}V_j  \subset U_j ,  \ \ \
              lim_{j\to\infty} \mu(V_j \, |\,  U_j) > 0 ; $
 \vspace{5mm}

 (iii)      $ \ \  \forall z\in  \tilde{{Q}_j} \ \exists u\in Q_j \ \
     \ (T_{v(j)}T_{z}\tilde{B}_j  \subset T_{u}{B}_j) ;    $
 \vspace{5mm}

 (iv)$ \ \   \max_{z\in  \tilde{{Q}_j}, u\in Q_j} \
\{\mu(T_{2v(j)}T_{z}\tilde{B}_j \ | \  T_{u}{B}_j)\} \ \to \ 0.$
 \vspace{5mm}

Прокомментируем это определение следующими примерами.
Пусть $T$ -- перемешивающий автоморфизм ранга 1. Как
показать, что он обладает D-свойством\,?
Пусть  $${U}_j = \bigsqcup_{0\leq z\leq  h(j)}  T_{z}\tilde{B}_j
$$ --
соответствующая последовательность башен для $T$.
Тогда в качестве  ${V}_j$ следует взять
башню $\bigsqcup_{0\leq z\leq 0.5h(j)}  T_{z}\tilde{B}_j,$
а в качестве $v(j)$ -- целую часть числа $0.5h(j)$.
Условия (i)--(iii) легко проверяются.  Покажем, что условие (iv) вытекает из
свойства перемешивания. Пусть
$$\mu(T_{2v(j)+z(j)}{B}_j \ | \  {B}_j)\ \approx d > 0, \ \ z(j) > 0. $$
Так как измеримое множество  $A$ аппроксимируется
для  достаточно большого $j$  $\xi_j$-измеримым множеством $A_j$,
 получим
$$\mu(T_{2v(j)+z(j)}A \ | \  A)\ \approx \
\mu(T_{2v(j)+z(j)}A_j \ | \  A_j)\ \geq d - \eps,
$$
что несовместимо со свойством перемешивания:
$$\mu(T_{2v(j)+z(j)}A \ | \  A)\to\mu(A).$$

Ситуация меняется, когда локальный ранг $\beta(T)$  близок к $\frac{1}{2}$.
Пусть  $\beta(T)=\mu(U_j)= 0.5+2\delta$.  Тогда   найдутся
$\tilde{B}_j, V_j$:
 $$
    \mu(\tilde{B}_j | B_j) > c(\delta)>0, \ \ \
V_j = \bigsqcup_{0\leq z\leq  \delta h(j)}  T_{z}\tilde{B}_j,
$$
причем для некоторой последовательности
  $v(j)$, $0.5h(j)<v(j)< (1-\delta)h(j)$, будут выполнены условия
(i)--(iii).  Условие (iv) опять будет следствием свойства перемешивания.
\medskip
Пусть $\Psi = \{ T^z : z \epsilon Z^n , \forall v,w\quad T^vT^w=T^{v+w} \}$
-- сохраняющее меру $Z^n$-действие на $(X, \cal {B},\mu )$.
\medskip

{\bf ТЕОРЕМА 4.2.1.}
\it   Если $\Z^n$-действие $\{T_z\}$ обладает свойством перемешивания и
$\beta\{T_z\}>\frac{1}{2^{n}}$, то действие обладает  \it D-свойством. \rm
\medskip

Доказательство.    Пусть
$(\tilde{Q}_j,  \tilde{\xi}_j, U_j)$ -- последовательность башен, которая
фигурирует в определении   $\beta(\{T_z\})$.
На  шаге с номером $j$ рассмотрим $2^n$ множеств вида
   $ T_{({w_1},{w_1},\dots,{w_n})}U_j$, где $w_k\in \{0,h_k(j)\}$.
Так как $\beta(\{T_z\})>\frac{1}{2^{n}}$,  то для достаточно больших $j$
для некоторой константы $c>0$
для различных $({w_1},\dots,{w_n})$ и $({w'_1},\dots,{w'_n})$
(которые зависят от $j$) имеем
$$
\mu(T_{({w_1},{w_1},\dots,{w_n})}U_j\cap
     T_{({w'_1},{w'_2},\dots,{w'_n})}U_j ) =  $$ $$
= \mu(U_j \cap    T_{({w'_1 -w_1},\dots,{w'_n}-w_n)}U_j ) > c. \eqno (4.1)
$$
Не ограничивая общности рассуждений, можно предположить, что  все
величины $w'_k-w_k$ неотрицательны. Общий случай сводится к этому при
подходящей  замене времени в рассматриваемом действии группы $\Z^n$.
При этом роль $B_j$ (основания башни) будет играть один
из ``угловых'' этажей башни.

Фиксируем маленькое число $\delta >0$ и
рассмотрим немного уменьшенную по сравнению с $U_j$ башню
$$ \bar{U}_j = \bigsqcup_{z\in \bar{Q}_j}  T_{z}{B}_j,\ \ \
$$ где $$\bar{Q_j} =   \{0,1,\dots,[(1-\delta)h_1(j)]\}\times
\dots\times\{0,1,\dots, [(1-\delta)h_n(j)]\},$$
и  маленькую башню
$$ \bar{V}_j = \bigsqcup_{z\in \delta{Q}_j}  T_{z}{B}_j,\ \ \
$$ $$\delta{Q_j} =   \{0,1,\dots,[\delta h_1(j)]\}\times
\dots\times\{0,1,\dots, [\delta h_n(j)]\}$$
($[\delta h]$ обозначает целую часть числа     $\delta h$).
Из (4.1) вытекает, что
$$
\mu(\bar{U}_j \cap    T_{2v(j)}V_j ) > a>0
$$
для некоторой последовательности  $\{v(j)\}$, удовлетворяющей условиям
$ v(j)\in Q_j$, $\ 2v(j) \notin Q_j$.
Осюда мы получим, что
$  \mu( {U}_j \ | \ T_{2v(j)}{B}_j )\ > \ a$, следовательно,
положив
$\tilde{B}_j =T_{-2v(j)}{U}_j\cap{B}_j,$
мы определим  башню
$V_j=\bigsqcup_{z\in \tilde{Q}_j}  T_{z}\tilde{B}_j$.
Условия (i),(ii),(iii) выполнены по построению.

Покажем, что условие (iv) вытекает из  свойства перемешивания
нашего действия с учетом условия $2v(j) \notin Q_j$ и
свойств последовательности   $\xi_j$.
Пусть (iv) не выполнено, т.е.
для некоторого $d>0$ и последовательности   $\{z(j)\}$ ($z(j)\in Q_j$)
неравенство
$  \mu(T_{2v(j)-z(j)}\tilde{B}_j \ | \ {B}_j)\ > d
$
выполнено  для бесконечного множества индексов $j$, причем
$|2v(j)-z(j)|\to\infty$.    Это приводит к  неравенствам

$
\mu(T_{2v(j)-z(j)}A \cap A)\ >  $
$$ > 0.9\mu(T_{2v(j)-z(j)}\tilde{B}_j \ | \ {B}_j)\mu(V_j)\mu(A) >
c \mu(A), \ \ (c>0),
$$
что для больших $j$ противоречит свойству перемешивания
$$\mu(T_{2v(j)-z(j)}A \cap A)\ \to  \mu(A)\mu(A),$$
когда мера множества     $A$  мала ($0<\mu(A)<c$).   
\medskip

\begin{center}
{\bf  4.3.  Тензорная простота перемешивающих систем с D-свойством}
\end{center}
\medskip

\bf ТЕОРЕМА 4.3.1. \it Перемешивающее $Z^n$-действие, обладающее
 D-свойством, является тензорно простым и, следовательно, обладает
 перемешиванием всех порядков.  \rm  \\
\medskip

Для доказательства теоремы достаточно показать, что мера $\nu \in M(2,3)$,
инвариантная  относительно  $T_z\times T_z\times T_z$,
есть $\mu\otimes \mu\otimes\mu$.
Предположим, что найдется эргодический джойнинг  $\nu\in$ M(2,3),
сингулярный относительно меры
$\mu\otimes\mu\otimes\mu$. Мере $\nu$  отвечает семейство
операторов  $\{P_x\}$, которое определяется формулой
$$
    \int \, \chi_A(x)\langle P_x\chi_B  \, , \,\chi_C \rangle d\mu(x)=
\nu(A\times B\times C),
$$
где $\langle\ ,\ \rangle$ -- скалярное произведение в $L_2(X,\mu)$.
Определим \it индуцированные джойнинги \rm  $\eta_z$:
$$
    \eta_z (A\times B\times B') =
\int \, \chi_A(x)\langle P_{T_{z}x}^\ast P_x\chi_B \,| \,
 \chi_{B'} \rangle d\mu(x).
$$
\medskip

{\bf ЛЕММА 4.3.2.}
\it
a) Если для некоторого
$z\in \Z^n\setminus\{0\}$ мера $\eta_z$ имеет
в качестве  компоненты меру $\mu\otimes\mu\otimes\mu$,
то $\nu=\mu\otimes \mu\otimes\mu$.

b)
Пусть последовательность множеств $V_j$  удовлетворяет условию
$\mu(V_j)>a>0$.  Пусть $v(j)$ -- некоторая последовательность,
$v(j)\in \Z^n$.
Тогда для некоторого $z\ne 0$ для бесконечного числа индексов $j$
имеет место неравенство
$$
      \eta_z (T_{v(j)}V_j\times V_j\times V_j) \ > \ 0.5a^3,
$$
причем мера $\eta_z$ принадлежит  классу M(2,3). \rm
\medskip

Доказательство.
a) Если   эргодический джойнинг $\nu$ сингулярен относительно
меры $\mu\otimes \mu\otimes\mu$, то
условные меры $\nu_x$, возникающие  при разложении
$$
      \nu(A\times A'\times A'') = \int_A \, \nu_x(A'\times A'')d\mu(x),
$$
сосредоточены на графиках конечнозначных отображений. (Это равносильно
выполнению для  некоторого $a>0$ и почти всех $x$
неравенств $P_x^\ast P_x\geq aI$.)
Действительно, проекция $\pi\eta_0$ меры $\eta_0$  на $X'\times X''$ не
равна  $\mu\otimes\mu$ (иначе, как известно,
$\nu=\mu\otimes \mu\otimes\mu$). В силу леммы 3.3.2  и неравенства
$\eta(X\times U_j\times U_j) \geq \mu(U_j)^2$,
получим, что  мера $\pi\eta_0$ имеет компоненту вида
$(I\times T_v)\Delta$. Для нетривиального джойнинга $\nu$
это возможно только в случае $v=0$.
 Таким образом,
выполнено $\int_X P_x^\ast P_x \, d\mu \geq aI$, что влечет
(ввиду эргодичности нашего действия) неравенство
$P_x^\ast P_x \, \geq aI$ для почти всех $x\in X$.

Таким образом, носитель   условной меры $(\eta_z)_x$ лежит на графике
композиции конечнозначных  отображений,
отвечающих операторам    $P^{\ast}_{T_zx}$  и   $P_x$.
Отсюда вытекает, что мера $(\eta_z)_x$  сингулярна относительно
$\mu\otimes\mu$ (для почти всех $x$).  Это влечет за собой сингулярность
меры $\eta_z$ относительно $\mu^3$.

\medskip
b)  Следуя  2.3,
для всех $z\in\Z^n$ получим представление
     $$
       P_{T_z(x)}^\ast P_x  = \frac{1}{p}
          ( J_x^{i(z,1)}+J_x^{i(z,2)}+\dots +J_x^{i(z,p)}),   \eqno (4.2)
     $$
где семейство операторов  $J_x^{i(z,q)}$ отвечает некоторому
эргодическому   джойнингу.

Несложно доказать (см. доказательство п. б) теоремы 3.3.3),  что
для некоторой последовательности  $z(k)\to\infty$  выполнено
$$
      \eta_{z(k)} (T_{v(j)}V_j\times V_j\times V_j) >c_k a^3,
$$
$$
\int_X P_{T_{z(k)}(x)}^\ast P_x d\mu(x) =
c_k\Theta + (1-c_k)P_k, \ c_k\to 1.
$$
Из (4.2) видно, что для больших $k$ правая часть последнего  равенства
есть  $\Theta$.  Это равносильно условию $\eta_{z(k)}\in$ M(2,3).

\medskip
Наша дальнейшая  цель -- доказать, что  некоторый индуцированный
джойнинг
$\eta_z$ имеет  компоненту $\mu\otimes\mu\otimes\mu$.
 Пусть $\eta\in M(2,3)$ является джойнингом трех копий действия
$\{T_z\}$.
В зависимости от  последовательности $\beta$-башен Рохлина-Халмоша
$\bigsqcup_{k\in Q_j}T^zB_j$,  определим  число  $Di(\eta)$:
$$Di(\eta)
    = \lim_{\eps\to +0} \limsup_{j\to\infty}\sum_{w1,w2,w3\in Q_j}
 a(\eps,w)\eta(T_{w1}B_j\times T_{w2}B_j\times T_{w3}B_j),
$$
где
 $$
    a(\eps,w)= \left\{
       \begin{array}{lcl}
          1, \ \ \ если \ \ \  \eta(T_{w1}B_j\times T_{w2}B_j\times T_{w3}B_j)
\leq \eps \mu(B_j)^2   \\
          0  \ \ \ иначе \ .
       \end{array}
           \right.
 $$
(Неформально говоря, $Di(\eta)$ при бесконечно  большом $j$ есть
сумма таких чисел $\eta(T_{w1}B_j\times T_{w2}B_j\times T_{w3}B_j)$, которые
бесконечно малы по сравнению  с величиной  $\mu(B_j)^2$.)

\medskip
{\bf ЛЕММА 4.3.3.} \it    Пусть перемешивающее действие $\{T_z\}$
обладает D-свойством. Для самоприсоединения $\nu \in M(2,3)$ найдется
индуцированный мерой $\nu$ джойнинг  $\eta_z$ такой, что
 $Di(\eta_z)>0$.
\rm

\medskip
Мы установим  неравенство $Di(\eta_z)>0$ для некоторого
$z\ne 0$, воспользовавшись
 D-свойством.
Пусть мера $\eta=\eta_z$ удовлетворяет свойствам пункта b) леммы 4.3.2,
в частности, выполнено условие
$$
      \eta_z (T_{v(j)}V_j\times V_j\times V_j) > c > 0.     \eqno (4.3)
$$
Предположим, что $Di(\eta)=0$. Тогда большинство
(относительно меры $\eta$)
значений
$\eta(T_{w}B_j\times T_{w'}\tilde{B}_j\times T_{w''}\tilde{B}_j)$
сравнимо с величиной  $\mu(B_j)^2$. Это влечет за собой, что для некоторого
натурального $N$ большинство
(относительно $\eta$) блоков вида
$T_{w1+v(j)}\tilde{B}_j\times
T_{w2}{B}_j\times T_{w3}{B}_j$  имеют $\eta$-меру большую, чем
$ \frac{1}{N}\mu(B_j)^2,$
где $w1,w2,w3\in \delta Q_j$ (обозначения из предыдущего параграфа).
 Рассмотрим образ  множества
$T_{w1+v(j)}\tilde{B}_j\times T_{w2}{B}_j\times T_{w3}{B}_j$
под действием преобразования
$T_{v(j)}\times T_{v(j)}\times T_{v(j)}.$
Учитывая (i),(ii),(iii) из определения D-свойства,  получим
$$T_{w1+2v(j)}\tilde{B}_j\times T_{w2+v(j)}{B}_j\times T_{w3+v(j)}{B}_j
= $$
$$
\bigsqcup_{u\in Q_j}\left(T_{u}{B}_j\cap T_{w1+2v(j)}\tilde{B}_j\right)
\times  T_{u2}\tilde{B}_j\times T_{u3}{B}_j .
$$
В силу свойства (iv) получим, что
числа  $\mu( T_{u}{B}_j\ |\ T_{w1+2v(j)}\tilde{B}_j)$
бесконечно малы (но их сумма по всем $u\in Q_j$ равна 1).
    Отсюда вытекает, что
под действием  $T_{v(j)}\times T_{v(j)}\times T_{v(j)}$
почти вся $\eta$-мера блока
$T_{w1}\tilde{B}_j\times T_{w2}{B}_j\times T_{w3}{B}_j$
распределилась среди блоков
$T_{u}B_j\times T_{u2}{B}_j\times T_{u3}{B}_j$, имеющих
$\eta$-меру, бесконечно малую по сравнению с $\mu(B_j)^2$.
Действительно, так как
мера $\eta$ принадлежит классу  M(2,3), то при
фиксированных  $u2,u3$  не более, чем $N$  таких блоков удовлетворяют
неравенству
$$\mu(T_{u}{B}_j\times T_{u2}B_j\times T_{u3}{B}_j)  \ \geq
 \ \frac{1}{N}\mu(B_j)^2 .
$$
Таким   образом,  предположение о том , что $Di(\eta)=0$
ввиду (4.3) приводит к  $Di(\eta)\geq c$.
 
\medskip

{\bf ЛЕММА 4.3.4.} \it    Пусть для эргодического джойнинга
$\eta\in M(2,3)$  выполнено
 $Di(\eta)>0$. Тогда $\eta =\mu\otimes \mu\otimes\mu$.
\rm
\medskip

Доказательствo.
Для некоторой  последовательности множеств  вида
$$
F_j =\bigsqcup_{z\in S_j} С^z_j =
\bigsqcup_{z\in S_j}\left(\bigsqcup_{v\in \delta Q_j}
T_{w_j+v}B_j\times T_{v}B_j\times T_{z+v}B_j\right)  \eqno   (4.4)
$$
выполнено $\eta(\ | F_j)\to \eta$ ( следствие того, что множество $F_j$
почти инвариантно относительно $T_z$ для фиксированных $z$, а мера
$\eta$ эргодична относительного нашего действия).
Множество индексов $S_j$  в формуле (4.4) выбирается таким образом,
чтобы для чисел $a_z^j=\eta(С^z_j   \,|\, F_j)$ выполнялось
$$\max_{z\in S_j}\{a_z^j\} \to 0, \  j\to\infty, \ \ \
 \sum_{z\in S_j}a_z^j \to a >  0.
$$
Условие  $Di(\eta)>0$  обеспечивает такой выбор  для  любой константы
$a$, удовлетворяющей условию $a < Di(\eta)$.

Обозначим через  $Y_j$  проекцию множества $С^z_j$ (для некоторого $z$)
 на второй сомножитель  в $X\times X\times X$.
Техника аппроксимаций (\cite{R0})
эргодических джойнингов сдвигами диагональных мер дает следующее:
 мера $\eta(\ | F_j)$   близка к мере $\lambda_j$, где $\lambda_j$ --
часть сдвига диагональной меры в $X\times X\times X$,
определенная равенством
$$ \lambda_j(A\times B \times C)\ =
\ \frac{1}{a\mu(Y_{z,j})}
\sum_z a_z^j \int_X \chi_{Y_{z,j}}\, \chi_{T_{w_j}A}
\,\chi_B\, \chi_{T_z C}\, d\mu .
$$
Из  теоремы
Блюма-Хансона для перемешивающих $\Z^n$-действий  имеем
$$\sum_{z\in S_j} a_z^j \chi_{T_z C} \ \to_{L_2}  \ Const \
\equiv \ a\mu(C), \   (j\to\infty)
$$
(короткое доказательство этой теоремы для $\Z$-действий см. в  1.1,
теорема 1.1.4;
в  случае $\Z^n$-действий  рассуждения аналогичны).
Теперь с учетом $\eta\in$ M(2,3) получим
$$\eta(A\times B \times C) =
\lim_{j\to\infty} \eta(A\times B \times C\, |\, F_j) =
 \lim_{j\to\infty}  \lambda_j(A\times B \times C) =
$$
$$
= \lim_{j\to\infty} \ \frac{1}{a\mu(Y_j)}
\int_X \chi_{Y_j}\, \chi_{T_{w_j}A}   \,\chi_B\,
\left( \sum_{z\in S_j} a_z^j   \chi_{T_zC}\right)\, d\mu =
$$
$$
=\mu(C)\,\eta(A\times B \times X) =
\mu(A)\mu(B)\mu(C).
$$

Мы доказали, что некоторый индуцированный
джойнинг  имеет компоненту $\mu\times\mu\times\mu$. Как было объяснено
перед формулой (4.2), это возможно лишь в случае
$\nu=\mu\times\mu\times\mu$.  Таким образом, свойство S(2,3) нашего
действия установлено, что завершает доказательство теоремы 4.3.1.
\newpage
{\bf  4.4.  Кратное перемешивание и локальный ранг}
\medskip


Ниже мы определим $(1+\varepsilon)$- перемешивание -- свойство,
занимающее промежуточное положение между свойством
перемешивания кратности 1 и кратности 2.
Положим

$
  Der(\varepsilon ,A,B,C) = $
$$ \lbrace (z,w)\in Q(\varepsilon ,h): | \mu (A\bigcap T^z B\bigcap
T^w C) - \mu (A)\mu (B)\mu (C)| > \varepsilon \rbrace ,
$$
где
$$
  Q(\varepsilon ,h) = \lbrace (z,w) \epsilon [0,h]^n :
|z|,|w|,|z-w| > \varepsilon h \rbrace ,
$$
$$
  d(h)=\sharp  Der(\varepsilon ,A,B,C)/ h .
$$
Если действие $\lbrace T^z \rbrace $ перемешивает двукратно, то $d(h)=0$
для больших значений $h$, если действие перемешивает, то последовательность
$d(h)$ ограничена.

Говорим, что $\bf \Psi $ $=\lbrace T^z : z \epsilon Z^n\rbrace$  является
$(1+\varepsilon)$- перемешивающим, если $d(h) \rightarrow \, 0 $ если
для любых измеримых множеств $A,B,C$ и $\varepsilon >0.$   \\ \bf
\medskip

\bf ТЕОРЕМА 4.4.1. \it $(1+\varepsilon)$-Перемешивающее $Z^n$-действие $\Psi $
положительного локального ранга является  тензорно  простым.  \rm

 Доказательство. Пусть для действия $\Psi $ найдется
нетривиальное самоприсоединение
$\nu \not= \mu \otimes \mu \otimes \mu $.
 Определим семейство операторов
$P_x : L_2(\mu ) \rightarrow L_2(\mu )$ формулой
$$
  \int \! \chi_A(x)\left<P_x\chi_B ,\chi_C \right>dm(x) = \nu (A\times B\times C),
$$
где $\chi_A $ является индикатором множества $A$. Рассмотрим новую меру
$\eta_m $ (здесь $m\epsilon Z^n$), определенную равенством
$$
    \eta_m (A\times B\times C) = \int \! \chi_A(x)\left<
P_x\chi_B ,P_{T^{m}x}\chi_C \right> dm(x)
$$
($\left<\,\ , \right>$ обозначает скалярное произведение в $L_2(\mu)$).

 Как и в предыдущем параграфе, выбираем такое $m$,
чтобы для некоторой эргодической компоненты
$\eta $ индуцированного джойнинга  $\eta_m $ выполнялось
$$
    lim_j   \eta (U_j\times U_j\times U_j ) > b^3 /2 ,
$$
причем для джойнинг $\eta $  обладает свойством парной независимости.
Поскольку мера $\eta $ сингулярна относительно $\mu \otimes \mu \otimes \mu $,
для любой (большой) константы $K$ найдутся некоторые множества $A,B,C$ такие,
что
$$
     \eta (A\times B\times C) > K \mu (A)\mu (B)\mu (C)
$$
Но $\eta $ аппроксимируется суммами
$\sum_{\gamma }a^{j}_\gamma \Delta^{j}_\gamma $ "почти инвариантных" мер
$\Delta^{j}_\gamma $, где $\gamma =(z,w) \epsilon Q(\varepsilon ,h_j)$ и
$$
  \Delta^{j}_\gamma (A\times B\times C) =
\mu (U_{j}^\gamma \bigcap A\bigcap T^zB\bigcap T^wC)/\mu(U_{j}^\gamma ),
$$
причем $\mu (U_{j}^\gamma )$ больше некоторого фиксированного положительного
числа. Подобнаяя техника аппроксимаций изложена в  \cite{MS}.  Используя
эргодичность меры $\eta $ и проекционные свойства $\eta $, можно доказать
следующее утверждение: число различных $\gamma $ таких, что
$\Delta^{j}_{\gamma }(A\times B\times C)$ близко к $\eta (A\times B\times C)$
и, следовательно, число
$$
  \sharp \lbrace (z,w)\epsilon  Q(\varepsilon ,h_j) :
\mu A\bigcap T^zB\bigcap T^wC) > 2 \mu (A)\mu (B)\mu(C) \rbrace
$$
сравнимо с $\sharp Q_j$. Это противоречит $(1+\e)$-перемешиванию.
\\  
\medskip

\begin{center}
{\bf  4.5. Ранг и джойнинги $T\times T$}
\end{center}
\medskip

Примеры автоморфизмов $T$
с   локальным рангом $\beta(T\times T)\geq\frac{1}{4}$
были предъявлены Катком в связи
с изучением спектральной кратности автоморфизмов пространства Лебега.
 Как мы покажем,
ранг  $T\times T$  равен бесконечности.
 Будет установлено, что
$\beta(T\times T)\leq\frac{1}{4}$,
а в случае равенства автоморфизм  $T$
обладает $ \kappa$-перемешиванием с показателем $ \kappa=\frac{1}{2}$
(как известно,
 это влечет за собой взаимную сингулярность
сверточных степеней спектральной меры автоморфизма $T$ [5]).
А. Каток сообщил автору
 пример перекладывания $T$ трех отрезков со свойством
$\beta(T\times T)=\frac{1}{4}$.
Эти результаты подтверждают давнюю гипотезу Оселедца  \cite{Os}
 о существовании
перекладываний со свойством  $ \kappa$-перемешивания.

Пусть автоморфизм $S$ пространства  $(X,\mu)$, $\mu(X)=1$,
обладает положительным  локальным рангом. Пусть
 $\beta(S)>0$ есть максимум чисел $\beta$ таких,
что  для последовательности конечных разбиений
  вида
$$ \xi_j=\{ B_j,  SB_j,    S^2B_j,    \dots,  S^{h_j-1}B_j,\dots\}$$
выполнено:
любое фиксированное измеримое множество аппроксимируется
$\xi_j$-измеримыми множествами при $j\to\infty$,
причем $\mu(U_j)\to \beta$, где  башня
$U_j=\bigsqcup_{0\leq k<h_j}S^kB_j$.

Мы рассмотрим случай, когда  $T\times T$ эргодичен
относительно  меры $\mu\otimes \mu$. Последнее равносильно
тому, что оператор $\Theta$ ортопроекции
на пространство констант в $L_2(X,\mu)$ принадлежит    слабому замыканию
множества $\{\T^n : n\in \Z\}$, где $\T$ обозначает
 унитарный оператор в $L_2(X,\mu)$,  индуцированный автоморфизмом $T$.
\medskip

{\bf ТЕОРЕМА 4.5.1. } \it   a) $Rank(T\times T)=\infty$.

b) $\beta(T\times T)\leq\frac{1}{4}$.
\rm  \medskip

Доказательство  проведем в предположении
эргодичности $T\times T$ (общий случай несложно редуцировать
 к рассматриваемому).
\rm  \medskip

{\bf ЛЕММА 4.5.2. } \it Пусть  эргодический автоморфизм $S$
пространства  $(\bar{X},\bar{\mu})$, коммутирует с автоморфизмом
$R$, и   $\beta(S)>0$.
Тогда  для любого    $\delta >0$
найдется $m>0$,  и последовательность
$n_{j}$ такая, что  имеет место слабая операторная  сходимость
$$  \hat{Y}_j\circ S^{n_{j}}\to  (1-\delta')\beta(S) \hat{R}^m,
$$
где  $ \hat{Y}_j$ -- оператор умножения на индикатор
множества ${Y_j}$,  причем ${Y_j}$ являются подмножествами
$\beta$-башен и
 $\bar{\mu}({Y}_j)\to (1-\delta')\beta(S)$ для некоторого $\delta'\in [0,\delta]$.
\rm
\medskip

{ Доказательство.}
Пусть $U_j$ -- последовательность $\beta$-башен для $S$.
Определим последовательность  маленьких башен
$ U_j^{\delta}=\bigsqcup_{0\leq k \leq \delta h_j}S^kB_j$.
Для некоторого числа $m>0$  выполнено
$$     \limsup_{j} \bar{\mu}(S^m U_{j}^{\delta} \cap  U_j^{\delta}) > 0.$$
Действительно, меры множеств  $ U_j^{\delta}$ стремятся
 к ${\delta}\beta(S)$,  отсюда вытекает, что множества
$$U_j^{\delta},  R U_j^{\delta},\dots,   R^N U_j^{\delta}$$
при ${\delta}\beta(S) N > 1$   не могут быть дизъюнктными.
Отсюда получаем, что для некоторого $m\leq N$ и некоторой
константы $c>0$   неравенство
$$\bar{\mu} (U_{j}^{\delta} \cap R^m U_{j}^{\delta}) >c $$
 выполнено для бесконечного числа номеров $j$.
В этом случае, как показано в  \S 3 работы \cite{R0},  джойнинг
$\Delta_{R^m}$,
отвечающий   оператору $R^m$, аппроксимируется  частью
сдвига $\Delta_{S^{n_j}}$ диагональной меры, расположенного
в $ Y_j\times \bar{X}$, где
$$Y_j=\bigsqcup_{\delta' h_j\leq k \leq h_j}S^kB_j$$ для некоторого
$\delta'\leq\delta$.
В операторной формулировке сказанное
записано в утверждении леммы.     
\medskip

{\bf ЛЕММА 4.5.3.}   \it  Пусть для эргодического автоморфизма  $T\times T$
выполнено $\beta(T\times T)>0$, и $U_j$ --
соответствующая последовательность $\beta$-башен.
Тогда найдется последовательность $k_j\to\infty$ такая, что
$$  \hat{U}_j\circ (\T\times\T)^{k_{j}}\to \beta(T\times T)(I\otimes\Theta)$$
($I$ --тождественный оператор).
\medskip                \rm

Доказательство. Применим предыдущую лемму  для случая
$S= T\times T$, $R=I\times T$.
Из доказательства   леммы 4.5.2 видно, что
при фиксированном $\delta >0$    множество $M_{\delta}$
чисел $m$,   для  которых
выполняется утверждение этой леммы, имеет положительную плотность.
Так  как $T$ -- слабо  перемешивающий автоморфизм,
слабое замыкание множества
$\{I\otimes\T^m\ : \ m\in M_{\delta}$\} содержит
оператор $I\otimes\Theta$.
Устремляя $\delta$ к нулю, получим утверждение леммы 4.5.3.
\medskip

{\bf ТЕОРЕМА 4.5.4. } \it
Если  локальный ранг эргодического автоморфизма $T\times T$ равен
$\frac{1}{4}$, то $T$ обладает $ \kappa$-перемешиванием при $ \kappa=\frac{1}{2}$. \rm
\medskip

Доказательство. Пусть
$$ (\T\otimes\T)^{k_j}\to\ P\otimes P  \ \geq \frac{1}{4}(I\otimes \Theta). $$
Представим $P$ в виде $P=aI+b\Theta +cQ$, где полиморфизм,
отвечающий марковскому оператору $Q$, сингулярен относительно
полиморфизмов $\Delta$ и $\mu\otimes \mu$, отвечающих операторам $I$ и $\Theta$.
Отметим также, что выполнено
$a,b,c\geq 0$, $a+b+c=1$. Тогда  в произведении  $$(aI+b\Theta +cQ)\otimes
(aI+b\Theta +cQ)$$
компоненту $(I\otimes \Theta)$  дает только  слагаемое $ab(I\otimes \Theta)$,
следовательно, $ab\geq \frac{1}{4}$, что  возможно только при
$a=b=\frac{1}{2}$.

Непосредственным следствием  леммы 4.5.3 являются сходимости
$$ \pi  \chi_{U_j} \circ \T^{k_{j}}\to\beta I, \ \ \
\pi  \chi_{ U_j^\ast} \circ \T^{k_{j}}\to\beta\Theta ,
$$
где $\pi F(x,y)$ обозначает $\int_X  F(x,y)d\mu(y)$,
а $U_j^\ast$  определено формулой
$$\chi_{ U_j^\ast}(x,y)=\chi_{ U_j} (y,x).$$
\medskip

{\bf ЛЕММА 4.5.5. }  \it
Пусть  $U_j$ -- последовательность $\beta$-башен
эргодического автоморфизма  $T\times T$.
Тогда имеет место сходимость
$$\int_X \pi   \chi_{U_j}  \pi  \chi_{U_j^\ast} d\mu  \ \to\ 0.$$
\rm
\medskip

Доказательство. Если
$$\int_X \pi   \chi_{U_j}  \pi  \chi_{U_j^\ast} d\mu  \ > c \ >\ 0,$$
то операторы $I$ и $\Theta$ имеют общую
компоненту, являющуюся предельной точкой последовательности
операторов
 $$\pi\chi_{U_j}\pi\chi_{U_j^\ast}\circ \T^{k_{j}}.$$
Но это невозможно, так как отвечающие им полиморфизмы
$\Delta$ и $\mu\otimes\mu$ взаимно сингулярны.
Таким образом, имеет место  асимптотическая
дизъюнктность проекций $\beta$-башни на сомножители в $X\times X$.

Теперь докажем пункт b) теоремы 4.5.1.
Из леммы 4.5.5 получаем:
для  некоторых непересекающихся множеств   $A_j$ и $B_j$ таких, что
$$\mu(A_j \Delta \{x : \pi \chi_{U_j} (x)> \eps>0\})\to 0$$ и
$$\mu(B_j \Delta \{x : \pi \chi_{U_j^\ast} (x)> \eps>0\}\to 0$$ выполнено
$$\mu\otimes\mu\left( U_j \setminus (A_j\times B_j)\right) < 2\eps.$$
Так как $\mu(A_j)\mu(B_j)\leq \frac{1}{4}$,
верхний предел последовательности $\mu\otimes\mu(U_j)$
 не превосходит $\frac{1}{4}$.
\medskip

Доказательство пункта a) теоремы 4.5.1.   Предположим, что
$Rank(T\times T)=r$, т.е. имеется некоторая
последовательность разбиений $\xi_j\to\eps$  вида
$$ \xi_j=\{  B^1_j,\,    SB^1_j,\,    S^2B^1_j,\,    \dots
,\,    S^{h_j^1}B^1_j,\ $$ $$ \dots\ \dots\ \dots\  $$
$$  B^r_j,\,  SB^r_j,\,   S^2B^r_j,\,   \dots ,
,\,  S^{h_j^r}B^r_j,\   Y_j  \}$$
(при этом
 ${h_j^1}$, ${h_j^2}\dots$, ${h_j^r}\to \infty$
и  $\mu(Y_j)\to 0$).
Тогда  меры объединений  соответствующих башен
$U_j^{1}, U_j^{2},\dots, U_j^{r}$  стремятся к $1$  при $j\to\infty$
(число $r$ фиксировано !).
Как мы пояснили выше,  башни из этого набора с
асимптотически исчезающей погрешностью  содержатся в
декартовых произведениях  непересекающихся множеств.
Но это противоречит следующей лемме, из которой следует,
что верхний предел упомянутой последовательности мер не
больше, чем  число $1 - 2^{-4r}$.
\medskip

{\bf ЛЕММА 4.5.6. }  \it
      Пусть  $A_1, A_2, \dots, A_r$  и  $B_1, B_2, \dots, B_r$
суть измеримые подмножества $X$, $\mu(X)=1$,
и для всех $i=1,2,\dots, r$
выполнено $A_i\cap B_i=\emptyset$.
Тогда
$$   \mu\otimes\mu\left(\bigcup_{i=1}^{r} (A_i\times B_i) \right) \
\leq \ 1 - 2^{-4r}.
$$
\rm
\medskip

{ Доказательство}.   Рассмотрим минимальное разбиение множества
$X$, порожденное множествами
$A_1,  \dots, A_r, B_1,  \dots, B_r$.
Пусть $C$ -- атом этого разбиения.
Замечаем, что
$$   \mu\otimes\mu\left( (A_i\times B_i) \bigcap (C\times C)\right) =0. $$
Действительно, если
$A_i\cap C$ и $B_i\cap C$ -- непустые множества,
то  $A_i\cap C=C=B_i\cap C$, так как $C$ является  атомом.
Но  это противоречит условию  $A_i\cap B_i=\emptyset$.
Таким образом, $C\times C$ не пересекается с объединением
множеств $A_i\times B_i$.
Всегда найдется такой атом $C$, что $\mu\otimes\mu(C\times C)\geq 2^{-4r}$.
\newpage 

\begin{center}
{\Large \bf  ГЛАВА 5} \\  { \bf
   НЕКОТОРЫЕ СПЕКТРАЛЬНЫЕ, \\ АЛГЕБРАИЧЕСКИЕ И АСИМПТОТИЧЕСКИЕ \\ СВОЙСТВА
           ДИНАМИЧЕСКИХ СИСТЕМ
}
\end{center}
\medskip

В главе  дается  положительное решение
проблемы Рохлина о непростом однородном спектре  для неперемешивающих и
перемешивающих автоморфизмов вида $T\times T$.
Для  $ \kappa$-перемешивающего автоморфизма
$T$ доказано, что изоморфизм $T\times T$ и $S\times S$ влечет за собой
изоморфизм $T$ и $S$.
Показано, что такое асимптотическое свойство
как частичное  кратное возвращение
может различать некоторые автоморфизмы с их
обратными. Рассмотрен  новый класс расширений,
 сохраняющих
 свойства тензорной простоты  и  кратного перемешивания.

\begin{center}
{\bf \S 5.1. Неперемешивающие автоморфизмы с однородным непростым спектром}
\end{center}
\medskip
В спектральной теории динамических
систем известна  задача В.А.Рохлина:  существует ли эргодический
автоморфизм $T$ с однородным спектром  кратности $m>1$?
Здесь подразумевается  кратность спектра унитарного оператора
$$\widehat{T}:L_2(X,\mu)\to L_2(X,\mu)\ \ \  \widehat{T} f(x)=f(Tx),$$
ограниченного на пространство функций с нулевым средним.
Напомним, что мы одинаково обозначаем автоморфизм пространства Лебега
и соответствующий оператор.

Впервые автоморфизмы с конечнократным непростым спектром были
построены В.И.Оселедцем \cite{Oseledets}.
В  лекциях А.Катка доказано, что
для типичного множества автоморфизмов
$T$ существенными значениями функции кратности  спектра
$T\times T$ будут наборы ${\cal M}_{T\times T}=\{ 2\}$
или ${\cal M}_{T\times T}=\{ 2, 4 \}$
(см. обзор \cite{GoodsonR}.)
Каток высказал гипотезу, что в типичном случае выполнено
${\cal M}_{T\times T}=\{ 2\}$.
Дж.Гудзон и М.Леманчик доказали в \cite{GoodsonL}, что
автоморфизм $T\times T$ не имеет  в спектре компоненты  нечетной кратности.

В работе \cite{G}  Дж.Гудзон рассмотрел
преобразование  $R:X\times X \to X\times X,$
опреденное  для фиксированного автоморфизма $T$ формулой
$$R(x,y)=(y,Tx),$$
и отметил, что из простоты спектра преобразования $R$ вытекает,
что спектр $T\times T$
является однородным кратности 2. Ниже показано,
что для класса автоморфизмов  соответствующее преобразование $R$ имеет
простой спектр.

А именно, рассматриваются автоморфизмы $T$ с простым спектром такие, что
при $a\in (0,1)$  оператор
$(aI+(1-a){\widehat{T}})$ (первый случай) или
$(1-a)(I+a{\widehat{T}}+a^2{\widehat{T}}^2+\dots)$ (второй случай)
 принадлежит слабому замыканию степеней оператора $\widehat{T}$.
Автоморфизмы  с этими свойствами строятся
в классе действий ранга 1.
Отметим, что преобразования $T$ со свойством
$\ \ {\widehat{T}}^{k_i}\to 0.5 (I+{\widehat{T}})$
рассматривались в статье А.Б.Катка и А.М.Степина \cite{KS1}.
Второй случай интересен тем, что при стремлении
$a$ к 1 автоморфизм $T$ становится все более перемешивающим.
 Системы с однородным спектром кратности 2  найдены автором
в классе перемешивающих систем
 (см. следующий параграф).

 Под спектром $T$ мы подразумеваем
спектральную меру  максимального типа $\sigma$
унитарного оператора $\widehat{T}$,
действующего в пространстве $H$ -- ортогональном дополнении
к постоянным функциям. Известно, что спектральный тип
$\widehat{T}\otimes\widehat{T}$
подчинен  $\sigma + \sigma\ast\sigma$.  Последнее слагаемое
(свертка) есть максимальный спектральный тип ограничения оператора
$\widehat{T}\otimes\widehat{T}$ на пространство $H\otimes H$.
В дальнейшем мы воспользуемся  тем, что взаимная сингулярность
спектральной  меры  $\sigma$ и свертки $\sigma\ast\sigma$
эквивалентна отсутствию ненулевых операторов, сплетающих
$\widehat{T} | H$ и  $\widehat{T}\otimes \widehat{T} | H\otimes H$.
(Пояснение: пусть $\widehat{T}\otimes \widehat{T} J=J\widehat{T}$, $J\neq 0$,
тогда найдется ненулевая комплексная мера $\lambda$ с коэффициентами Фурье
$$  \widehat{\lambda}(n)=\langle \widehat{T}^n f|J^\ast g\rangle_{L_2(\mu)}   =
\langle \widehat{T}^n\otimes \widehat{T}^n  Jf|g\rangle_{L_2(\mu\otimes\mu) },
$$
которая подчинена одновременно спектральной мере
$\sigma$ и свертке $\sigma\ast\sigma$.)

\medskip

{\bf ТЕОРЕМА 5.1.1.} \it Пусть $T$ -- эргодический автоморфизм,  и для
некоторой последовательности ${k_i}\to \infty$ и числа $a\in (0,1)$
выполнено
$\ \ {\widehat{T}}^{k_i}\to  (aI+(1-a){\widehat{T}})$.
Тогда

1) для спектральной меры $\sigma$ автоморфизма $T$
  выполнено $\sigma\ast\sigma\perp\sigma$;

2) если $T$ имеет простой спектр, то
автоморфизм $R$, $R(x,y)=(y,Tx),$  имеет простой спектр,
а автоморфизм
$(T\times T)$ имеет однородный  спектр кратности 2.
\vspace{5mm}
\rm
\medskip

(Утверждение, аналогичное теореме 5.1.1, независимо и иначе доказал
О.Н.Агеев \cite{Age}.)

Доказательство.
Предположим, что  ограниченный оператор $J$ сплетает ${\widehat{T}}$ с
оператором $({\widehat{T}}\otimes {\widehat{T}})$:
$$J{\widehat{T}}= (\widehat{T}\otimes \widehat{T})J.$$
Получим при   $b=1-a$
$$ J(aI + b{\widehat{T}})=
((aI + b{\widehat{T}})\otimes (aI + b{\widehat{T}}))J,$$

$$ (a(I\otimes I) + b({\widehat{T}}\otimes{\widehat{T}}))J =
 (a^2(I\otimes I) + ab({\widehat{T}}\otimes I) +ab (I\otimes \widehat{T}) +b^2(\widehat{T}\otimes \widehat{T}))J,$$

$$ J + (\widehat{T}\otimes \widehat{T})J =
 (I\otimes \widehat{T})J +   (\widehat{T}\otimes I)J. $$
Отсюда вытекает, что для любых  $i,j$ выполнено
$$(\widehat{T}^i\otimes \widehat{T}^j)J +
(\widehat{T}^{i+1}\otimes \widehat{T}^{j+1})J -
(\widehat{T}^{i}\otimes \widehat{T}^{j+1})J -
(\widehat{T}^{i+1}\otimes \widehat{T}^{j})J =0.
$$
Теперь получаем
$$  \sum_{0\leq i,j <n} (\widehat{T}^i\otimes \widehat{T}^j)J + (\widehat{T}^{i+1}\otimes \widehat{T}^{j+1})J -
 (\widehat{T}^{i}\otimes \widehat{T}^{j+1})J -   (\widehat{T}^{i+1}\otimes \widehat{T}^{j})J =  0,
$$
$$
  J + (\widehat{T}^n\otimes \widehat{T}^n)J -
 (I\otimes \widehat{T}^n)J -  (\widehat{T}^n\otimes I)J =0.       \eqno (5.1)
$$
Так как $T$ слабо перемешивающий ( нетрудно проверить, что
собственной функцией может быть только константа),
для некоторой последовательности
 ${n_i}\to \infty$
оператор $\Theta$  ( ортопрекция на пространство констант)
является слабым пределом  степеней:  $\widehat{T}^{n_i}\to\Theta$.

Учитывая (5.1),  получим
$$J + (\Theta\otimes\Theta) J = (I\otimes\Theta) J + (\Theta\otimes I)J.$$
Следовательно,  $Im(J)\perp H\otimes H$, что означает отсутствие
ненулевых операторов,  сплетающих $\widehat{T} | H$ с  $\widehat{T}\otimes\widehat{T} |H\otimes H$.
Последнее влечет за собой свойство $\sigma\ast\sigma\perp\sigma$.
Пункт 1) доказан.

 Доказательство пункта 2).   Достаточно показать, что $R$
имеет простой спектр.
 Пусть $f$ --  циклический вектор оператора $\widehat{T}$, действующего в
пространстве $H$.
Докажем, что
$V_{0,0}=f\otimes f$ является циклическим вектором для ограничения
оператора $\widehat{R}$ на
инвариантное пространство $H\otimes H$.
С этой целю установим, что все векторы $V_{m,n}= T^mf\otimes T^nf$
принадлежат пространству  $L$ -- замыканию линейной оболочки
множества векторов  $\{R^iV_{0,0}: i\in {\bf Z}\}$.
Заметим, что   $\widehat{R} V_{m,n}=V_{n,m+1}$.
Векторы
$V_{0,1},    V_{0,0},  V_{1,1} $ принадлежат $L$.
Но из    $ {\widehat{T}}^{k_i}\to  (aI+b{\widehat{T}})$  получаем
$$[(aI+b\widehat{T})\otimes (aI+b\widehat{T} )]V_{0,0}= a^2V_{0,0}+b^2V_{1,1}
+abV_{0,1}+abV_{1,0}\ \in L.$$
Следовательно, $V_{1,0}\in L,$ и $V_{0,2}= \widehat{R} V_{1,0}\in L.$

Докажем по индукции, что $V_{0,p+1}\in L$.
Предположим, что  $V_{0,i}\in L$ при $i=0,1, \dots, p$.
Так как слабое замыкание степеней автоморфизма является полугруппой,
то эта полугруппа   содежит  операторы вида  $(aI+b\widehat{T})^p$.
Следовательно, вектор
$U_p=[(aI+b\widehat{T})^p\otimes (aI+b\widehat{T} )^p]V_{0,0}$ принадлежит $L$.
Вектор $U_p$ является линейной комбинацией векторов  $V_{m,n}$,
$0\leq m,n\leq  p$, причем мы уже знаем, что все векторы $V_{m,n}$, кроме
$V_{p,0}$, принадлежат $L$. Но из $U_p\in L$ получим
$V_{p,0}\in L$.
Следовательно,
$V_{0,p+1}=\widehat{R} V_{p,0}\in L.$

Таким образом, все $V_{p,0}$ и, следовательно, все $V_{0,p}$ принадлежат
циклическому пространству $L$. Окончательно получаем
$V_{m,n}=R^{2m}V_{0,n-m}\in L,$ т.е. $L=H\otimes H.$
Таким образом, пространство $H\otimes H$ является циклическим пространством
оператора $\widehat{R}$, т.е.
 ограничение  $\widehat{R}$ на $H\otimes H$ имеет  простой  спектр.
 Очевидно, что
ограничение оператора $\widehat{R}$ на $(1\otimes H)+ (H\otimes 1)$
также имеет простой спектр.
Приходим к выводу, что автоморфизм $R$ имеет простой  спектр, так как
ограничение оператора $\widehat{R}$ на $(1\otimes H)+ (H\otimes 1)$ дизъюнктно
с ограничением $\widehat{R}$ на $(H\otimes H)$ : нет ненулевого сплетающего
оператора. Последнее вытекает  из того,
что ограничение оператора $\widehat{R}^2$ на $(1\otimes H)+ (H\otimes 1)$ дизъюнктно
с ограничением $\widehat{R}^2$ на $(H\otimes H)$.
  
\vspace{5mm}

{\bf ТЕОРЕМА 5.1.2} \it Пусть для эргодического автоморфизма $T$
 для некоторой
последовательности ${k_i}\to \infty$ и числа  $a\in (0,1)$  выполнено
$ {\widehat{T}}^{k_i}\to \
{(1-a)}(I+a\widehat{T} +a^2\widehat{T}^2+\dots).$
Тогда

1) для спектральной меры $\sigma$ автоморфизма $T$ выполнено
$\sigma\ast\sigma\perp\sigma$;

2) если $T$ имеет простой спектр, то
автоморфизм $R$  имеет простой спектр, а автоморфизм
$(T\times T)$ имеет однородный  спектр кратности 2.
\vspace{5mm}

\rm

Доказательство п. 1).  Обозначим $P= (1-a)(I-a\widehat{T})^{-1}$.
Пусть оператор $J:H\to H\otimes H$ удовлетворяет условию сплетения
$$J{\widehat{T}}= ({\widehat{T}}\otimes {\widehat{T}})J.$$
Так как $ {\widehat{T}}^{k_i}\to P$, получим
$$ JP=  (P\otimes P) J,  $$

$$ J{(1-a)}(I+a\widehat{T} +a^2\widehat{T}^2+\dots) =(P\otimes P)J,$$

$$ (1-a)[I\otimes I+ a(\widehat{T}\otimes\widehat{T})
 +a^2(\widehat{T}\otimes\widehat{T})^2+\dots]J
=(P\otimes P)J,$$

$$ (1-a)[I\otimes I - a(\widehat{T}\otimes\widehat{T})]^{-1} J =
(1-a)^2(I -a\widehat{T})^{-1}\otimes (I-a\widehat{T})^{-1}J. $$
Так как для коммутирующих операторов
$A=(I\otimes I - a(\widehat{T}\otimes\widehat{T}))$
и $B=(I -a\widehat{T})\otimes (I-a\widehat{T})$ равенство $A^{-1}J=B^{-1}J$
влечет за собой $AJ=BJ$, мы получим
$$ {(1-a)}[(I\otimes I) - a(\widehat{T}\otimes\widehat{T})] J
 =   (I-a\widehat{T})\otimes (I-a\widehat{T})J,$$
$$  [ I\otimes I +  \widehat{T}\otimes\widehat{T}]J
 = [ I\otimes\widehat{T} + \widehat{T}\otimes I]J.$$
Последнее равенство, как показано в доказательстве предыдущей теоремы,
приводит к $J=0$, что влечет за собой $\sigma\ast\sigma\perp\sigma$.

Доказательство п. 2). Наша цель  --  показать, что
пространство $H\otimes H$ является циклическим для
оператора $\widehat{R}$.  Для этого достаточно
предъявить систему вложенных циклических пространств
$$ C_0 \subseteq C_1\subseteq C_2 \subseteq C_3 \subseteq \dots,  $$
объединение которых плотно в $H\otimes H$.
Пусть $C_n$ -- циклическое пространство опрератора $\widehat{R}$, порожденное
вектором $W_n$, где
$$  W_n= [(I-a\widehat{T})^n \otimes (I-a\widehat{T})^n] f\otimes f,
$$
$f$ -- циклический вектор для $\widehat{T}$ на $H$,  $n$ принимает значения
$0,1,2,\dots$.
Так как $C_n$ инвариантно относительно действия
$\widehat{R}^2=\widehat{T}\otimes \widehat{T}$,
при $\widehat{T}^{k(i)}\to P$  получим
$(P\otimes P)C_{n}\subseteq C_{n}$.
Но    $(P\otimes P)W_n  =W_{n-1}$, откуда вытекает требуемое включение
$C_{n-1}\subseteq C_n$.

Обозначим через $L$  замыкание объединения пространств $C_n$.
Чтобы доказать, что $L=H\otimes H$ достаточно установить,
что все $V_{n,0}$ принадлежат $L$, так как
в этом случае все  $V_{0,n}=\widehat{R} V_{n-1,0}\in L$, следовательно, все
векторы вида $\widehat{T}^k\otimes \widehat{T}^kV_{n,0}$
или вида $\widehat{T}^k\otimes \widehat{T}^kV_{0,n}$
  ($k\in\Z, \ n=0,1,2, \dots$ )  принадлежат    пространству $L$.

Докажем $V_{1,0}\in L$.
Мы знаем, что  векторы
$$V_{0,1}, V_{0,0}, V_{1,1},\ W_1=
[(I -a\widehat{T} )\otimes (I -a\widehat{T})]V_{0,0}$$
принадлежат   пространству $C_1$.
Так как $$W_1= V_{0,0} + a^2V_{1,1} -a V_{0,1} - aV_{1,0},$$
получим, что   $V_{1,0}\in C_1$. Тогда имеем
$\widehat{R} V_{1,0} = V_{0,2}\in C_1$.

Докажем, что  $V_{2,0}\in C_2$. Вектор  $W_2$ является
линейной комбинацией векторов $V_{i,j}$, $0\leq i,j \leq 2$.
Мы установили, что все векторы, кроме  $V_{2,0}$, лежат в
$C_1$. Но  $W_2\in C_2$ и $C_1\subseteq C_2$, следовательно,
$V_{2,0}\in C_2$.

Рассуждая по индукции, для всех $n=0,1,2,\dots $ устанавливаем
$V_{n,0}\in C_n\subseteq L$. Как пояснялось в доказательстве теоремы 1,
это приводит к равенству $L=H\otimes H$ и простоте спектра
оператора $\widehat{R}$.  

\begin{center}
{\bf \S 5.2. Перемешивающие автоморфизмы с однородным непростым спектром}
\end{center}
\medskip

В этом параграфе мы докажем существование перемешивающего автоморфизма
 $T$, обладающего следующими свойствами:

1. Спектр симметрического произведения $T\odot T$ простой:
 ${\cal M}_{T\odot T}=\{ 1\}$,

2. ${\cal M}_{T\times T}=\{ 2\}$,

3. $\sigma_T\ast\sigma_T\perp\sigma_T$.

Напомним, что ${T\odot T}$ означает ограничение
${T\times T}$ на фактор $S$-неподвижных измеримых множеств в $X\times X$,
где $S$ -- симметрия $S(x,y)=(y,x)$.
Очевидно, что свойство 1  влечет за собой 2  и 3.

Как будет показано, такой автоморфизм можно найти в классе
так называемых лестничных конструкций.
Напомним ее определение. Пусть
автоморфизм $T$ допускает последовательность разбиений $\xi_n$
фазового пространства
$X$ следующего вида:
\medskip

$   B^1_n,\ \   TB^1_n,\ \    \dots
,\ \      T^{h_n-1}B^1_n,$
\medskip

$   B^2_n,\ \   TB^2_n,\ \     \dots
,\ \      T^{h_n-1}B^2_n,\ \   T^{h_n}B^2_n,$
\medskip

$   B^3_n,\ \   TB^3_n,\ \     \dots
,\ \      T^{h_n-1}B^3_n,\ \   T^{h_n}B^3_n,\ \ T^{h_n+1}B^3_n,$
\medskip

\dots \ \dots \ \dots \     \dots \ \dots \ \dots \     \dots \ \dots \ \dots \
 \dots \ \dots \ \dots \
\medskip

$ B^{r_n}_n,\   TB^{r_n}_n,\   \dots
,\   T^{h_n}B^{r_n}_n,\   T^{h_n+1}B^{r_n}_n,  \
\ \ \dots, \ T^{h_n+r_n-2}B^{r_n}_n, \ \ \ Y_n$
\medskip
такие, что для всех $n$ выполнено

$B^2_n = T^{h_n}B^1_n  ,   $

$ B^3_n = T^{h_n+1}B^2_n,     $

$\dots    $

$  B^{r_n}_n = T^{h_n+r_n-1}B^{r_n-1}_n .    $
\\
Если разбиения $\xi_n$ стремятся к разбиению на точки
($\xi_n\to\eps$) и
для всех  $n$ выполнено
$$B^1_{n-1}= B^1_{n}\cup B^2_{n} \cup \dots\cup B^{r_{n}}_{n}, $$
говорим, что  $T$ является лестничной кострукцией.
Можно заметить, что конструкция однозначно (с точность до изоморфизма)
 определена параметрами $h_1$ и
$\{r_n\}$.

Теорема Адамса \cite{Adams} утверждает, что из условий $r_n \to\infty$ и
$\frac{(r_n)^2}{h_n} \to 0$  вытекает свойство перемешивания
соответствующего $T$.

ЗАМЕЧАНИЕ. Как показал автор, для перемешивания достаточно
потребовать только условие
$r_n \to\infty$ с единственным ограничением на конечность меры $X$.
\medskip

{\bf ТЕОРЕМА 5.2.1.}
\it Существует перемешивающая лестничная конструкция $T$, обладающая свойством
 ${\cal M}_{T\odot T}=\{ 1\}$.

{\bf СЛЕДСТВИЕ.} Соответствующий перемешивающий
автоморфизм  $T\times T$
имеет однородный спектр   кратности 2.   \rm

\medskip

 ОПРЕДЕЛЕНИЕ.   Будем говорить, что лестничная конструкция
$T$ принадлежит  классу $St.C.(p, h)$,
если последовательность $r_n$ удовлетворяет
условиям:

1. $\liminf_{n\to\infty} r_n =p$;

2. для каждого $q, \ p\leq q\leq h,$ найдется $n_i\to \infty$
такая, что  $r_{n_i +1}\to \infty$ и
$\forall \ \ i \ \ r_{n_i} =q$.
\\
Из последнего условия при $p\leq q\leq h$ мы получим, что  операторы  вида
$$P_q=\frac{1}{q}(I+\T+\T^2+\dots + \T^{q-2} +\Theta)$$
принадлежат $WCl(T)$ -- слабому замыканию степеней $T$.
\medskip

{\bf ЛЕММА.} \it
  Пусть $T$ принадлежит классу $St.C.(p, h+2)$,  $3p<h+2$.
Пусть $B,TB,\dots, T^hB$ -- непересекающиеся измеримые множества,
обозначим через
$C_{B\times B}$
циклическое пространство, порожденное
$\chi_{B}\otimes \chi_{B}$  под действием оператора $\T\otimes \T$.
Тогда функции
$$F_{i,j}=\chi_{T^iB}\otimes \chi_{T^jB}+
\chi_{T^jB}\otimes \chi_{T^iB}, \ \ \ 0\leq i,j\leq h,
$$
принадлежат циклическому пространству $C_{B\times B}$.
\rm
\medskip

Доказательство. Так как
$$F_{q,0}\in C_{B\times B} \Longrightarrow  F_{q+i,i}\in C_{B\times B},
$$
достаточно установить $F_{q,0}\in C_{B\times B}$.
Функции
$$      P_{q+2}\chi_{B}\otimes  P_{q+2} \chi_{B},  \ \
P_{q+1} \chi_{B}\otimes  P_{q+1}\chi_{B},  \ \
         P_{q}\chi_{B} \otimes  P_{q}\chi_{ B}\
$$
принадлежат циклическому пространству $ C_{B\times B}$. Положим
$$
G_{m}= \frac{m^2}{(1+\mu(B))^2} P_{m} \chi_{B}\otimes  P_{m}\chi_{B}.
$$
Заметим, что
$$F_{q,0}= Const [G_{q+2} - G_{q+1}   -  (\T\otimes \T) G_{q+1} +
(\T\otimes \T) G_{q}].
$$
Таким образом, для всех  $q\geq p$  получили $F_{q,0}\in C_{B\times B}$,
следовательно,
 $F_{i,j}\in C_{B\times B}$.

Теперь покажем, что при
$q < p$  будет выполнено $F_{q,0}\in C_{B\times B}$.
Рассмотрим случай $q=1$ (общий случай доказывается аналогично).
Так как
$$F_{p,0}= [(\T^p\otimes I)+(I\otimes \T^p)]\chi_{B\times B}$$
и   $C_{B\times B}$  содержит   $F_{p+1,0}$, получим:
циклическое пространство, порожденное $F_{p,0}$ содержит
\vspace{5mm}

$[(\T^p\otimes I)+(I\otimes \T^p)]F_{p+1,0} = $ $$
[(\T^p\otimes I)+(I\otimes \T^p)]  [(\T^{p+1}\otimes I)+(I\otimes \T^{p+1})]
\chi_{B\times B}.$$
Последнее представим в виде суммы  $F_{2p+1,0} + F_{p+1, p}$, где
$F_{2p+1,0}\in C_{B\times B}$. Так как сумма принадлежит
$C_{B\times B}$, получим $F_{p+1, p}\in C_{B\times B}$, значит,  $F_{1, 0}$
принадлежит $C_{B\times B}$.
\\
 
\medskip

\bf ТЕОРЕМА 5.2.2. \it Если $T$ принадлежит классу \\ $St.C.(p,\infty)$,
то   ${\cal M}_{T\odot T}=\{ 1\}$.
\rm
\medskip

Доказательство.   Пусть $\xi_n$ -- соответствующая последовательность
разбиений.
Рассмотрим циклические пространства $C_{B_n\times B_n}$
для оператора $\T\otimes\T$.
Каждая симметрическая функция $F(x,y)$ может быть аппроксимирована
линейными комбинациями $F_{i,j}$, где $F_{i,j}$  ($0\leq i,j\leq h_n$)
зависят также от  $n$.   Из леммы получим, что
 пространство $L_2(\mu)\odot L_2(\mu)$ аппроксимируется
циклическими пространствами  $C_{B_n\times B_n}$.
Это означает, что
$L_2(\mu)\odot L_2(\mu)$ также является циклическим пространством.
Получаем ${\cal M}_{T\odot T}=\{ 1\}$.

Заметим, что автоморфизмы класса $St.C.(p,\infty)$ не являются
перемешивающими.  Однако с их помощью можно построить
перемешивающий  автоморфизм $T$  со свойством
 \\ ${\cal M}_{T\odot T}=\{ 1\}$.
Для этого мы будем рассматривать последовательность
$\{r_n \}$    такую, что
для всех $k$ выполнено $r_{2k+1} =  2k+1$, причем
последовательность $r_{2k}$ очень медленно стремится к бесконечности
при $k\to\infty$.
Оказывается, что лестничная конструкция, задаваемая такой последовательностью
 $\{r_n \}$, является искомой.

Теперь поясним, как выбирать $\{r_n \}$
(отметим, что выбор наш неконструктивен).
На шаге с номером $j$ рассматриваем автоморфизм
$T_j$ класса $St.C(p_j,\infty)$
на $X_j= \{x : T_j(x)\neq x \}$.
Автоморфизму $T_j$ соответствует последовательность $r_n^{(j)}$.
Выберем (очень большое) число $N_j$ и
изменим последовательность $r_n^{(j)}$
при $n>N_j$.  Получим новую последовательность
$r_n^{(j+1)}$, которой соответствует
автоморфизм $T_{j+1}$, отличающийся от $T_j$
только на очень маленьком множестве $Y_{N_j}$ -- крыше над башней.
В результате предельный автоморфизм $T$  будет
лестничной конструкцией с $r_n\to\infty$. Выбирая
$r_n \leq n$, мы получим $\frac{(r_n)^2}{h_n} \to 0$,
что в силу теоремы Адамса гарантирует свойство перемешивания.

Вкратце стратегия следующая.
Последовательность $\{T_{j}\}$  выбираем так,
чтобы для любой симметричной функции \\ $F(x,y)\in L_2(\mu)\otimes L_2(\mu)$
 расстояние между
$F$ и циклическим пространством $C_{B_j\times B_j}$  оператора
            $T_j\otimes T_j$
стремилось к нулю.  Если обеспечить быстрое стремление $T_j$ к $T$,
то получим, что последовательность
циклических пространств $C_{B_j\times B_j}$  оператора
    $T\odot  T$  аппроксимирует $L_2(\mu)\odot  L_2(\mu)$
(напомним, что $T\odot  T$ обозначает ограничение оператора
$T\otimes  T$ на пространство симметричных функций $L_2(\mu)\odot  L_2(\mu)$) .
Это, как хорошо известно, означает цикличность пространства
$L_2(\mu)\odot  L_2(\mu)$ относительно оператора
$T\odot   T$, что нам и требуется.

Теперь подробнее.   Последовательность
$\{T_{j}\}$  строим таким образом, чтобы  выполнялось
$$ supp (T_j)\subset supp (T_{j+1}) \subset \dots supp(T)=X.$$

Пусть на шаге $j$ задан набор  функций
$$\{ F_k\in L_2(\mu)\otimes L_2(\mu) : k=1,2,\dots, j\}, \ \
supp (F_{k+1})\subset  supp(T_k).$$
Добавляем  к набору произвольным образом симметричную функцию  $F_{j+1}$,
$supp (F_{j+1})\subset  supp(T_j)$, но при условии, что   семейство
$\{F_{j}\}$ будет плотным  в  $L_2(X,\mu)\odot L_2(X,\mu)$.
Пользуясь  цикличностью    пространства
$L_2(supp(T_{j+1}),\mu)\odot L_2(supp(T_{j+1}),\mu)$
для          $U_{j+1}  =\T_{j+1}\odot\T_{j+1}$,
находим такое $N$ и такие наборы коэффициентов $a_{ik}$, что
выполнено
$$ \|F_k - \sum_{i=-N}^{N}a_{ik} U_{j+1}^i \chi_{B_{j+1}\times B_{j+1}}\|
 \ < \ \frac{1}{j},
$$
где $U_{j+1}  = \T_{j+1}\odot\T_{j+1}$.
При достаточно быстром убывании последовательности
$\mu(\{x\ :\ T_j(x)\neq T_{j+1}(x)\})$ (здесь наши возможности не ограничены)
имеем неравенства
$$ \|F_k - \sum_{i=-N}^{N}a_{ik} U^i \chi_{B_{j+1}\times B_{j+1}}\|
 \ < \ \frac{2}{j},
$$
где $U  = \T\odot\T$.
С учетом сказанного  выше,  получаем, что $U$ имеет простой спектр.
 .

\begin{center}
{\bf \S 5.3.
        Изоморфизм  декартовых степеней преобразований и $ \kappa$-перемешивание}
\end{center}
\medskip
Пусть $S$ и $T$  --
сохраняющие меру обратимые преобразования вероятностного пространства Лебега
$(X,\mu)$ и известно, что их декартовы квадраты $S\times S$ и $T\times T$
метрически изоморфны. Последнее означает, что для некоторого
сохраняющего меру $\mu\otimes\mu$ обратимого преобразования
$\Phi: X\times X\to X\times X $  выполнено
$$
      \Phi (S\times S) \Phi^{-1} = T\times T \, .
$$
Будут ли изоморфны $S$ и $T$ ?

    В этом параграфе доказывается следующая
\medskip

 \bf ТЕОРЕМА 5.3.1. \it Если автоморфизм $T$ обладает свойством
 $ \kappa$-перемешивания, то изоморфизм $T\times T$ и $S\times S$
 влечет за собой изоморфизм $T$ и $S$.
\rm
\medskip

 Доказательство.
Пусть для оператора $\Phi$, отвечающего автоморфизму
$(X\times X,\mu\otimes \mu)$, выполнено равенство
$$
      \Phi (S\otimes S) \Phi^{-1} = T\otimes T \, .
$$
Преобразование $T\times T$ обладает двумя  координатными факторами,
которые под действием $\Phi$ переходят в некоторые факторы преобразования
$S\times S$. (Напомним, что фактором называется ограничение преобразования
на инвариантную  $\sigma$-подалгебру алгебры всех измеримых множеств.)
Координатные алгебры образованы множествами вида $X\times A$  и
множествами вида $A\times X$ -- им соответствуют пространства
${\bf 1}\otimes L_2(X,\mu)$ и $L_2(X,\mu)\otimes {\bf 1}$.

    Наша цель -- доказать, что $\Phi$ переводит координатную алгебру в
координатную. Тем самым, очевидно, будет осуществлен изоморфизм между
преобразованиями $S$ и $T$.
Для этого в нашем случае достаточно установить, что образ
пространства ${\bf 1}\otimes L_2(X,\mu)$ под действием $\Phi$
совпадает с пространством
${\bf 1}\otimes L_2(X,\mu)$ или $L_2(X,\mu)\otimes {\bf 1}$.

Пусть
$T^{n(i)} \to (1- \kappa )I +  \kappa\Theta$ для некоторого $ \kappa \in (0,1)$.
Докажем, что из (1.1) вытекает
$$
S^{n(i)} \to (1- \kappa )I +  \kappa\Theta .
$$
Имеем  сходимость
$$
\Phi (S \otimes S)^{n(i)} \Phi^{-1} \to \left((1- \kappa )I +  \kappa\Theta\right)
\otimes \left((1- \kappa )I +  \kappa\Theta\right)\, .
$$
Так как марковская полугруппа компактна ( в топологии слабой сходимости),
для некоторой подпоследовательности $n(i')$ и некоторого
марковского оператора $Q$ получим
$$
(S \otimes S)^{n(i')}\to Q\otimes Q \,.
$$
Учитывая, что  $I\otimes I$ и $\Theta\otimes \Theta$ коммутируют с $\Phi$,
а последовательность
$
 (S \otimes S)^{n(i)} \Phi^{-1}$
 сходится к
$$ \Phi\left((1- \kappa )I +  \kappa\Theta\right)
\otimes \left((1- \kappa )I +  \kappa\Theta\right) \Phi^{-1},
$$
получим
$$
  Q\otimes Q = (1- \kappa )^2 \left(I\otimes I\right)\ + \
 \kappa^2 \left(\Theta \otimes \Theta\right) \ + \
$$
$$
           +    \kappa (1- \kappa ) \Phi^{-1}\left((I\otimes\Theta )+
                                   (\Theta\otimes I)\right)\Phi.
$$
Таким образом, для любой неотрицательной функции $f \in L_2(X,\mu)$
выполнены неравенства
$$
  Qf\otimes Qf \geq  (1- \kappa )^2 f\otimes f\, ,
$$
$$
  Qf\otimes Qf \geq   \kappa^2 \Theta f\otimes \Theta f\, .
$$
Из этих неравенств вытекает, что
$$
            Qf \geq (1- \kappa ) f , \quad Qf \geq   \kappa \Theta f\, ,
$$
следовательно, в силу произвольности положительной функции $f$ получаем
$$
Q= (1- \kappa )I +  \kappa\Theta \, .
$$
Так как выбор сходящейся последовательности
$S^{n(i')}$ был произвольным, получим
$$
S^{n(i)}\ \to \ (1- \kappa )I +  \kappa\Theta \, .
$$

Пусть $H\subset L_2(X,\mu)$ обозначает пространство функций с нулевым средним.
Тогда имеет место разложение
$$
L_2(X\times X,\mu\otimes \mu) =\{Const\}\oplus
(H\otimes {\bf 1})\oplus ({\bf 1}\otimes H) \oplus (H\otimes H)\, ,
$$
где $\{Const\}$ обозначает пространство постоянных функций. Покажем, что
$$
\Phi(H\otimes {\bf 1})\subseteq (H\otimes {\bf 1})\oplus ({\bf 1}\otimes H)\, .
$$
Действительно, последовательность $(T\otimes T)^{n(i)}$ сходится на
линейном пространстве $(H\otimes {\bf 1})$ к оператору
$(1- \kappa)(I\otimes I)$, следовательно,
последовательность $(T\otimes T)^{n(i)}\, =\,
\Phi(S\otimes S)^{n(i)}\Phi^{-1}$   сходится на
$\Phi(H\otimes {\bf 1})$ также к оператору $(1- \kappa)(I\otimes I)$.
Это возможно лишь в случае, когда выполнено
$$\Phi(H\otimes {\1})\subseteq
(H\otimes {\1})\oplus ({\1}\otimes H),$$
так как на   пространстве
$(H\otimes H)$ эти последовательности операторов сходятся к
$(1- \kappa)^2(I\otimes I),$ но $(1- \kappa)^2 <(1- \kappa)\,.$
Теперь покажем, что в действительности имеет место
$$
\Phi(H\otimes {\bf 1})\subseteq (H\otimes {\bf 1})
\quad {или}
\quad \Phi(H\otimes {\bf 1})\subseteq({\bf 1}\otimes H).  \eqno (5.2)
$$
Для любой ограченной функции $f\in H$ существуют функции
 $f_1$ и $f_2$ такие, что
$$
\Phi(f\otimes {\bf 1}) = a(f_1\otimes {\bf 1}) + b({\bf 1}\otimes f_2)\, .
$$
Из того, что пространство
            $\Phi(L_2(X,\mu))\otimes {\bf 1})$
соответствует фактору и замкнуто относительно произведения
ограниченных функций, вытекает
$$
 (\Phi(f\otimes {\bf 1}))^2 \in (H\otimes {\bf 1})\oplus ({\bf 1}\otimes H)
\oplus \{Const\}\,.
$$
Тогда получим
$$
ab(f_1\otimes f_2)\in (H\otimes {\bf 1})\oplus ({\bf 1}\otimes H)
\oplus \{Const\}\,,
$$
поэтому
   $ab = 0$, или $f_1$ есть константа, или $f_2$ есть
константа. Таким образом, мы показали, что
$$
\Phi(H\otimes {\bf 1})\subset (H\otimes {\bf 1})\cup ({\bf 1}\otimes H)\,.
$$
Так как $H\otimes {\bf 1}$ и ${\bf 1}\otimes H$ имеют нулевое пересечение,
а $\Phi(H\otimes {\bf 1})$ есть линейное пространство, получим (5.2).
Пусть, например, имеет место
$\Phi(H\otimes {\bf 1})\subseteq (H\otimes {\bf 1})$. Повторим рассуждения,
поменяв местами $S$ и $T$. Тогда получим
$\Phi^{-1}(H\otimes {\bf 1})\subseteq (H\otimes {\bf 1})$\, . Эти включения
вместе дают равенство
$\Phi(H\otimes {\bf 1}) = (H\otimes {\bf 1})$ и требуемый изоморфизм
координатных факторов.

Теорема 5.3.1 допускает следующее естественное обобщение.
\medskip

 \bf ТЕОРЕМА 5.3.2. \it  Пусть $T$ --  $ \kappa$-перемешивающий автоморфизм,
$0< \kappa <1$, и для автоморфизма $S$ выполнено
$$
      \Phi S^{\otimes n} \Phi^{-1} = T^{\otimes n},
 $$
где $\Phi$ -- оператор, отвечающий автоморфизму пространства
$(X^{\times n},\mu^{\otimes n})$.
Тогда $S$ и $T$ изоморфны. \rm
\medskip

Рассуждения, аналогичные предыдущим,  при $f\in H$ приводят
последовательно к  утвеждениям:
$$
  \Phi(f\otimes {\bf 1}\otimes\dots \otimes {\bf 1}) \in
    (H\otimes {\bf 1}\otimes\dots \otimes {\bf 1})\oplus\dots
     \oplus ({\bf 1}\otimes\dots\otimes {\bf 1}\otimes H),
$$
функция $
 (\Phi(f\otimes {\bf 1}\otimes\dots \otimes {\bf 1}))^2$
 принадлежит пространству
$$
(H\otimes {\bf 1}\otimes\dots \otimes {\bf 1})\oplus\dots\oplus ({\bf 1}
\otimes\dots\otimes {\bf 1}\otimes H)
\oplus \{Const\}\,,
$$
следовательно,
$$
  \Phi(H\otimes {\bf 1}\otimes\dots \otimes {\bf 1}) \subset
    (H\otimes {\bf 1}\otimes\dots \otimes {\bf 1})\cup\dots
     \cup ({\bf 1}\otimes\dots\otimes {\bf 1}\otimes H).
$$
Последнее  возможно только в случае, когда
$$
  \Phi(H\otimes {\1}\otimes\dots \otimes {\bf 1}) =
    (\1\otimes\dots\otimes\1\otimes H\otimes
{\bf 1}\otimes\dots \otimes {\bf 1}),
$$
т.е. под действием $\Phi$   образ пространства координатного фактора системы
$T^{\otimes n}$
совпадает  с пространством некоторого координатнатного фактора  системы
$S^{\otimes n}$.
Это приводит к  изоморфизму координатных факторов и, следовательно,
к изоморфизму $T$ и $S$.
\medskip
\medskip
\begin{center}
 {\bf \S 5.4.
     Асимметрия прошлого и будущего динамической системы
            и кратная возвращаемость}
\end{center}
\medskip
Если автоморфизм $T$ пространства Лебега
$(X,\B,\mu)$, $\mu(X)=1$, метрически неизоморфен
автоморфизмому $ T^{-1}$,
будем называть его асимметричным.
Первый пример асимметричного действия
был опубликован в \cite{Anz}.
В \cite{Ose} предъявлены асимметричные каскады с
простым непрерывным спектром. В этих конструкциях  использовалось
косое произведение
(интересное обсуждение возникновения этого понятия
имеется в статье  \cite{An2}).
     О современном состоянии  проблемы симметрии прошлого и будущего
динамической системы  см., например,  работу \cite{Goodson1}.
В работе \cite{F96} отмечалось, что  $T$ и $T^{-1}$
не сопряжены
в полной группе эргодического автоморфизма  $T$.

Цель этого параграфа -- указать  новый асимптотический  инвариант
(кратное возвращение на подпоследовательностях),
который может различать  автоморфизмы  $T$  и  $T^{-1}$.
\medskip

{\bf ТЕОРЕМА 5.4.1.} \it
Существует автоморфизм $T$, обладающий
свойством:  для некоторой последовательности  \\
$n(i)\to\infty$ для любого множества  $A\in\B$  выполнено
$$
 \lim_{i\to \infty} \mu(A\cap T^{n(i)}A\cap T^{3n(i)}A) \geq
\frac{1}{5} \mu(A),                                             \eqno (5.3)
$$
при этом для некоторого множества $A', \ \mu(A') > 0,$ имеет место
$$
\lim_{i\to \infty} \mu(A'\cap T^{-n(i)}A'\cap T^{-3n(i)}A') =  0.   \eqno (5.4)
$$
В качестве такого множества $A'$ годится
любое множество, удовлетворяющее условию
$\mu(A'\cap TA')=\mu(A'\cap T^2A')=0.$   \rm
\medskip

{\bf СЛЕДСТВИЕ.} \it   Автоморфизм $T$ асимметричен.
\medskip       \rm

Свойство (5.3), которое можно назвать  \it
кратной возвращаемостью \rm с коэффициентом
$a=\frac{1}{5}$,  является инвариантом автоморфизма. Действительно, пусть
для некоторого автоморфизма $S$
$$S^{-1}TS =R,$$
тогда для любого  $A\in \B$   из (5.3) получим
$$
 \lim_{i\to \infty} \mu(SA\cap T^{n(i)}SA\cap T^{3n(i)}SA) \geq
a \mu(SA),
$$
следовательно,
$$
 \lim_{i\to \infty} \mu(A\cap S^{-1}T^{n(i)}SA\cap S^{-1}T^{3n(i)}SA) =
$$ $$
= \lim_{i\to \infty} \mu(A\cap R^{n(i)}A\cap R^{3n(i)}A)
\geq    a \mu(A).
$$
Из предположения, что $T$ и $T^{-1}$ сопряжены,
вытекает
 $$\lim_{i\to \infty}\mu(A'\cap T^{-n(i)}A'\cap T^{-3n(i)}A') \geq  a \mu(A), $$
но это противоречит (5.4).

{\bf Автоморфизмы, удовлетворяющие свойствам (5.3) и (5.4).}
В классе преобразований ранга 1 найдется такой автоморфизм $T$, что
для некоторой последовательности $h(i)\to\infty$
имеем следующее представление фазового пространства:
$$
  X= Y_i\sqcup Y^1_i \sqcup Y^2_i\sqcup  Y^3_i\sqcup  Y^4_i\sqcup  Y^5_i,
$$
где $\mu(Y_i)\to 0$,  а  множества
$Y^1_i,Y^2_i,Y^3_i,Y^4_i,Y^5_i$ имеют следующий вид:
\medskip
\medskip

$
  Y^1_i = B^1_i\cup TB^1_i\cup  \dots
\cup T^{h_i-2}B^1_i\cup T^{h_i-1}B^1_i,
$
\medskip

$
  Y^2_i = B^2_i\cup TB^2_i\cup \dots
\cup T^{h_i-2}B^2_i  \cup T^{h_i-1}B^2_i  \cup T^{h_i}B^2_i,
$
\medskip

$
  Y^3_i = B^3_i  \cup TB^3_i \cup \dots
  \cup T^{h_i-2}B^3_i  \cup T^{h_i-1}B^3_i  \cup T^{h_i}B^3_i,
$
\medskip

$
  Y^4_i = B^4_i  \cup TB^4_i \cup \dots
  \cup T^{h_i-2} _i \cup T^{h_i-1}B^4_i \cup T^{h_i}B^4_i\cup
          T^{h_i+1}B^4_i,
$
\medskip

$
  Y^5_i =  B^5_i  \cup TB^5_i \cup \dots
  \cup T^{h_i-2}B^5_i  \cup T^{h_i-1}B^5_i
                 \cup T^{h_i}B^5_i \cup  T^{h_i+1}B^5_i.
$
\\ Потребуем также, чтобы  выполнялось:
$$
  T^{h_i}B^1_i= B^2_i,
\  T^{h_i+1}B^2_i= B^3_i,
 \  T^{h_i+1}B^3_i= B^4_i,
  \ T^{h_i+2}B^4_i= B^5_i;
$$
$$
\mu(T^{h_i+2}B^5_i\Delta B^1_i)/\mu(B^1_i) \ \to \ 0;
$$
$$
  \xi_i = \{ C_i, TC^1_i, T^2C_i, \dots , T^{h_i-1}C_i\}\
       \to \ \eps,
$$
где
$
C_i= B^1_i\cup B^2_i\cup B^3_i\cup B^4_i\cup B^5_i,
$
а запись $\xi_i  \to \ \eps$ означает, что любое множество
$A\in \B$ аппроксимируется $\xi_i$-измеримыми множествами.

Положим $n(i)=h(i)+1$.
Мы утверждаем, что при $a=\frac{1}{5}$ для любых $A,B,C,D\in \B$ выполнено:
$$
 \mu(Y^1_i\cap A\cap T^{n(i)}B\cap T^{2n(i)}C\cap T^{3n(i)}D) \to
  a \mu(A\cap T^{-1}B\cap T^{-2}C\cap T^{-2} D)     ,
$$
$$
 \mu(Y^2_i\cap A\cap T^{n(i)}B\cap T^{2n(i)}C\cap T^{3n(i)}D) \to
  a \mu(A\cap T^{+1} B\cap T^{0}C\cap T^{-1} D)      ,
$$
$$
 \mu(Y^3_i\cap A\cap T^{n(i)}B\cap T^{2n(i)}C\cap T^{3n(i)}D) \to
  a \mu(A\cap T^{0} B\cap T^{+1}C\cap T^{0} D)      ,
$$
$$
 \mu(Y^4_i\cap A\cap T^{n(i)}B\cap T^{2n(i)}C\cap T^{3n(i)}D) \to
  a \mu(A\cap T^{0} B\cap T^{0}C\cap T^{+1} D)   ,
$$
$$
 \mu(Y^1_i\cap A\cap T^{n(i)}B\cap T^{2n(i)}C\cap T^{3n(i)}D) \to
  a \mu(A\cap T^{-1}B\cap T^{-1}C\cap T^{-1} D)  .
$$
Чтобы проверить эти утверждения,  вместо множеств $A,B,C,D$ следует
на каждом $i$-том шаге подставлять  $\xi_i$-измеримые множества,
которые аппроксимируют $A,B,C,D$, а  затем перейти к пределу.

Таким образом, при $B=A,\ D=A, \ C=X$  получим
$$
\lim_{i\to\infty} \mu(A\cap T^{n(i)}A\cap T^{3n(i)}A) \geq
   a\mu(A),
$$
так как
$$
 \mu(Y^3_i\cap A\cap T^{n(i)}A\cap T^{3n(i)}A) \to
  a \mu(A \cap T^{0} A \cap T^{0} A) = a\mu(A).
$$
Подстановка $B=X,\ C=A, \ D=A$ дает
$$
    5 \mu( A\cap  T^{2n(i)}A\cap T^{3n(i)}A)\ \to \
[\ \mu(A\cap T^{-2}A\cap T^{-1} A)\ +
$$
$$
 + \mu(A\cap T^{0}A\cap T^{-2} A) \ +\
 \mu(A\cap T^{+1}A\cap T^{0} A) \ + \
$$
$$  + \mu(A\cap  T^{+0}A\cap T^{+1} A)    \ + \
   \mu(A\cap  T^{-1}A\cap T^{-1} A)\ ].
$$
При условии $$\mu(A\cap TA)=\mu(A\cap T^2A)=0$$ имеем
$$
 \lim_{i\to\infty}\mu(A\cap T^{2n(i)}A\cap T^{3n(i)}A) = 0.
$$
Так как
$$
 \mu(A\cap T^{2n(i)}A\cap T^{3n(i)}A) =
\mu(T^{-3n(i)}(A\cap T^{2n(i)}A\cap T^{3n(i)}A)),
$$
теперь получим
$$
 \lim_{i\to\infty}\mu(A\cap T^{-n(i)}A\cap T^{-3n(i)}A) = 0.
$$

{\bf Асимметрия и кратное частичное перемешивание.}
Существует
автоморфизм $T$ такой, что для  некоторого числа  $b>0$
и некоторой последовательности $m(i)$ выполнено
$$\forall A,B,C\in\B \ \ \ \
 \lim_{i\to \infty} \mu(A\cap T^{m(i)}B\cap T^{3m(i)}C) \geq
 b\mu(A)\mu(B)\mu(C),
$$
причем  для некоторых множеств $A',B',C'$ положительной меры имеем
 $\mu(A'\cap T^{-m(i)}B'\cap T^{-3m(i)}C') \to  0.$
В этом случае $T$ и $T^{-1}$ не могут быть изоморфными .
Такие примеры предложены автором  в \cite{R03}.

Упомянем задачу, близкую по содержанию
к проблеме Рохлина о кратном перемешивании.
Неизвестно, существует ли такой автоморфизм  $T$, что
для некоторых последовательностей  $m(i)$, $k(i)$ выполнено  условие:
для любых множеств $A,B,C\in\B$
 $$
 \lim_{i\to \infty} \mu(A\cap T^{m(i)}B\cap T^{k(i)}C)\ = \
 \mu(A)\mu(B)\mu(C),                          \eqno (5.5)
$$
но существуют множества $A',B',C'$ такие, что
$$
 \lim_{i\to \infty} \mu(A'\cap T^{-m(i)}B'\cap T^{-k(i)}C') \neq
 \mu(A')\mu(B')\mu(C').                        \eqno (5.6)
$$
Отметим, что  в классе автоморфизмов с сингулярным спектром таких примеров нет.
Действительно, для некоторой последовательности $i'\to\infty$
найдется мера $\nu\neq \mu \otimes \mu\otimes \mu$ на
$(X\times X\times X,\, \B\times\B\times\B)$  такая, что
$$
 \mu(A\cap T^{-m(i')}B\cap T^{-k(i')}C) \ \to \
 \nu(A\times B\times C).
$$
Из (5.5) легко получить, что   проекции меры $\nu$  на грани
$X_{(i)}\times X_{(j)}$
совпадают с  $\mu \otimes \mu$. Результат \cite{Hos} гласит: если
спектр $T$ сингулярный, то  такая мера при условии
ее инвариантности относительно
$T\times T\times T$ есть  $\mu\otimes\mu\otimes\mu$.   Таким образом, если
спектр $T$ сингулярный, то из (5.5) вытекает
$$\mu(A\cap T^{-m(i)}B\cap T^{-k(i)}C) \to   \mu(A)\mu(B)\mu(C),$$
что противоречит (5.6).

\begin{center}
{\bf \S 5.5. Расширения, сохраняющие  кратное перемешивание и тензорную простоту}
\end{center}
\medskip
Напомним, что мы используем обозначение  $M(n,n+1)$ для класса мер
на декартовом $(n+1)$-мерном кубе $X^n$ с проекциями на $n$-мерные грани,
совпадающими с произведением  $\mu^n= \mu^{\otimes n}$.

\medskip
{\bf ТЕОРЕМА 5.5.1. } \it  Пусть автоморфизм $S$ перемешивает с кратностью $k$,
а косое произведение $R$,
$R(x,y) = (S(x),T_x(y))$, является перемешивающим. Если
все преобразования  $T_x$ коммутируют с некоторым тензорно простым
действием $\Psi$ (быть может, некоммутативным), то косое произведение
$R$    перемешивает с кратностью $k$.
\rm
\medskip

Доказательство. Рассмотрим случай $k=2$. Обозначим
через $\lambda$ меру $\mu\otimes\mu$, которая инвариантна относительно
 $R$. Пусть для некоторых последовательностей
$z_0,z_1,z_2\to\infty$ и некоторой
меры $\nu\in M(2,3)$ (класс $M(2,3)$ рассматривается здесь
относительно пространства $(X\times Y,\lambda)$)  выполнено
$$
\lambda (R^{z_0(j)}A_0 \cap R^{z_1(j)}A_1  \cap R^{z_2(j)}A_2) \quad
\to \quad \nu (A_0\times A_1 \times A_2) \ \ (j\to\infty).
$$
Мера   $\nu$ задана на  произведении
$(X\times Y)\times (X\times Y)\times (X\times Y)$,
которое нам удобно  представить   в виде
$(X\times X\times X)\times (Y\times Y\times Y)$.
Из того, что $S$ перемешивает кратно, вытекает, что
проекция меры $\nu$ на $(X\times X\times X)$  совпадает с
$\mu^3=\mu\otimes\mu\otimes\mu$.  Меру $\nu$ разложим в систему
условных мер $\{\nu_w \,:\, w\in X\times X\times X\}$:
$$
  \nu(A_1\times A_2\times A_3\times B_1\times B_2\times B_3)=
\int_{A_1\times A_2\times A_3} \nu_w(B_1\times B_2\times B_3)d\mu^3(w) .
$$
Для почти всех $w$ верно, что $\nu_w\in M(2,3)$.
Это вытекает из
того, что мера  $\mu\otimes\mu$ эргодична относительно
$\Psi\otimes\Psi$.
Так как косое произведение $R$ коммутирует с действием
$\Phi= I\otimes\Psi$, мера $\nu$
инвариантна относительно $\Phi\otimes\Phi\otimes\Phi$.
Действительно,  для любого $\varphi$ из $\Phi$ имеем
$$
\nu (\varphi A_0\times \varphi A_1 \times \varphi A_2) =
$$
$$
=\lim_{j\to\infty} \lambda (R^{z_0(j)}\varphi A_0
\cap R^{z_1(j)}\varphi A_1
 \cap R^{z_2(j)}\varphi A_2) =
$$
$$
=\lim_{j\to\infty} \lambda ( (R^{z_0(j)} A_0 \cap R^{z_1(j)} A_1
 \cap R^{z_2(j)} A_2)) =
\nu ( A_0\times  A_1 \times  A_2).
$$
Из установленной инвариантности меры $\nu$ получаем,
что условные меры  $\nu_w\in M(2,3)$ инвариантны
относительно  $\Psi\otimes\Psi\otimes\Psi$. Но $\Psi$ является
тензорно простой
системой, следовательно, почти все  $\nu_w$ равны
$\mu\otimes\mu\otimes\mu$. Таким образом, мы получили, что
$\nu = \lambda\otimes\lambda\otimes\lambda$. Тем самым
для косого произведения $R$ свойство кратного
перемешивания порядка $k=2$ установлено.
Случай $k>2$ рассматривается аналогично.
 
\medskip

{\bf ТЕОРЕМА 5.5.2. } \it Пусть $R,T$ -- перемешивающие преобразования,
где $R$ есть косое произведение над $S$ следующего вида:
$$
  R(x,y)=(S(x),T^{n(x)}(y)), \quad \int n(x)d\mu = 0.
$$
Если автоморфизм $S$ перемешивает с кратностью $k$, то  косое произведение
$R$ также обладает перемешиванием  кратности $k$.
Если автоморфизм $S$ является тензорно простым,
 то  $R$ также является тензорно простым.
\medskip

\rm

Теорема 5.5.2  вытекает, как мы покажем, из следующего ключевого
утверждения.
\medskip

{\bf ТЕОРЕМА 5.5.3.} \it  Пусть ${\cal C} : X^n \to M(n-1,n)$. Если $(S,T^{f(x)})$
-- перемешивающее  косое произведение,
где $T$ -- перемешивающее преобразование и $\int f d\mu =0$,  то уравнение
$$
   {\cal C} (S(x_1),\dots,S(x_n))\equiv
(T^{f(x_1)}\otimes \dots\otimes T^{f(x_n)}){\cal C}(x_1,\dots,x_n)
$$
имеет единственное решение: ${\cal C}(x_1,\dots,x_n) \equiv \mu^n$.
\rm
\medskip

Доказательство теоремы 5.5.3. проведем после вспомогательных утверждений.

Цилиндрическим каскадом  называется
отображение $(S,f): (X\times \Z)\to (X\times \Z)$,
определенное соотношением:
$$\forall x\in X,\,\forall a\in \Z\quad
(S,f)(x,a) = (S(x), a + f(x)),
$$
где  $S$ -- преобразование множества $X$,
а $f:X\to \Z$ -- некоторая функция.
\medskip

{\bf ТЕОРЕМА 5.5.4.} (Крыгин, Аткинсон) \it
Пусть $S$ есть эргодическое преобразование, и
для функции $f:X\to\Z$ выполнено  $\int f(x)d\mu = 0$, тогда

a) цилиндрический каскад $(S,f): (X\times \Z)\to (X\times \Z)$
является консервативной системой,

b) для почти всех $x\in X$ найдется бесконечная
последовательность натуральных
чисел $\{q_i(x)\}$ такая, что $q_0(x)=0,\,$
$\,q_i(x) < q_{i+1}(x)$
и для всех  $i$ выполнено
$$
    \sum_{n=0}^{q_i(x)-1}f(S^{n}(x)) = 0.
$$
\rm
\medskip

Доказательство теоремы см. в \cite{Kry}, \cite{Atk};
пункт b) вытекает непосредственно из a)
и означает, что множество
$\{N\,: \,\sum_{n=0}^{N}f(S^{n}(x)) = 0\}$ бесконечно для почти всех
$x\in X$.
(Отметим, что теорема, аналогичная п. b),  выполнена для
потоков \cite{Sch}).
\medskip

 {\bf ЛЕММА 5.5.5.} \it   Пусть $T$ перемешивает, и задана мера
$\nu\in M(n,n+1)$. Если хотя бы одна из последовательностей
$m_1(i),\dots,m_n(i)$ стремится к бесконечности при $i\to\infty$,
то для последовательности мер
$\nu_i=(I\times T^{m_1(i)}\times\dots\times T^{m_n(i)})\nu$
выполнено  $\nu_i\to\mu^{n+1}$.
\rm
\medskip

Доказательство проведем для случая $n=2$.
Рассмотрим оператор  $P:L_2^{\otimes 2} \to L_2$ такой, что
$$\langle P(\chi_A\otimes \chi_B)\,,\, \chi_C\rangle =
\nu(A\times B\times C).$$
Предположим для определенности, что $m_1(i)\to\infty$.
Нам нужно доказать, что
$$P_i = ( T^{m_1(i)}\otimes T^{m_2(i)})P \to P\Theta, \quad (i\to\infty) $$
( $\Theta$ -- ортопроектор на пространство  констант в $L_2(\mu)$).
Для некоторой последовательности $i'$ имеем
$$( T^{m_1(i')}\otimes T^{m_2(i')})P \to
 \left(\Theta\otimes Q\right)P = ( \Theta\otimes \Theta)P,
$$
где равенство возникает из-за того, что оператор $P$ внутренний.
Предельному оператору $(\Theta\otimes \Theta)P = P\Theta$
отвечает мера $\mu^3$.
Так как мере $\nu_i$ соответствует оператор $P_i$,
получим $\nu_{i'}\to\mu^3$.
Отсюда вытекает, что $\nu_{i}\to\mu^3$, так как выбор
сходящейся подпоследовательности $\nu_{i'}$ был произвольным.
 
\medskip

  {\bf ЛЕММА 5.5.6.} \it
Если $(S,T^{f(x)})$  перемешивает, то при $i\to\infty$  для любого $K$
меры множеств $ \{ x\,:\, |F(x,i)|<K \}$ стремятся к  0.
\rm
\medskip

Доказательство.
Рассмотрим множество  вида
$ U=X\times \bigcup_{-N<i<N}T^iY, $
где $T^iY$ -- непересекающиеся множества, а число $N$ выбрано
так,  чтобы выполнялось $a-2\frac{K}{N} > 2\tilde{\mu}(U)>0$.
Требуемое множество $Y$ всегда найдется (для построения
можно воспользоваться известной леммой Рохлина-Халмоша).

Предположим, что
для некоторой последовательности $i'\to\infty$   будет выполнено
$$ \mu(\{ x\,:\, |F(x,i')|<K \}) > a > 0.$$
Тогда получим
$$
\lim_{i\to\infty}\tilde{\mu}(U\cap R^iU) > a\tilde{\mu}(U)-2\frac{K}{N} >
2\tilde{\mu}(U)\tilde{\mu}(U),
$$
что противоречит свойству перемешивания автоморфизма $R$.  

\medskip
Доказательство теоремы 5.5.3. Для простоты рассмотрим случай $n=3$.
Зададим метрику на множестве $M(2,3)$. Пусть
некоторое семейство $\{B_i\}_{i\in\N}$    всюду плотно в $\B$, и задана
биекция $\sigma :\N^3 \to \N$. Определим расстояние между мерами
$\nu$ и $\eta$ следующим образом:
$$ \rho(\nu\,,\,\eta) = \sum_{i,j,k}
   {2^{-\sigma(i,j,k)}}{|\nu(B_i\times B_j\times B_k)\,-\,
\eta(B_i\times B_j\times B_k) |}.
$$
 Положим
$F(x,i) = \sum_{n=0}^{i}f(S^{n-1}(x)) .
$
Для всех $i$  выполнено
$$
   {\cal C}S^i(x),S^i(y),S^i(z))\equiv
(T^{F(x,i)}\otimes T^{F(y,i)}\otimes T^{F(z,i)}){\cal C}x,y,z).
$$

 Из эргодичности  $S\otimes S\otimes S$ вытекает, что почти всюду в
$X\times X\times X$ либо  ${\cal C}x,y,z)\neq \mu^3$,
либо  ${\cal C}x,y,z)\equiv \mu^3$.
Пусть выполнено первое, тогда для некоторого
$c>0$ мера тех $(x,y,z)$, для которых $\rho({\cal C}x,y,z), \mu^3)>c$
больше, чем $c$. Для некоторого
множества $A,\,\mu(A)>0$ и числа $b>0$  выполнено:
$$
  \forall x\in A \quad \mu\otimes\mu(\{(y,z):
  \rho({\cal C}x,y,z), \mu^3)>c \})>b .               \eqno (5.7)
$$
Заметим, что для почти всех точек из $A$
множество $\{i'\, : \,S^{q_{i'}(x)}(x)\in A\}$ бесконечно.
Это вытекает из того, что преобразование $U(x)=S^{q_1(x)}(x)$
(при условии, что число ${q_1(x)>0}$ выбирается минимальным для всех $x$)
будет обратимо почти всюду,
следовательно, оно сохраняет меру и  для
него выполняется свойство возвращаемости:
для почти всех $x$  множество
$\{i'\, : \,U^{i'}(x)\in A\}$ бесконечно
(ниже обозначим $i'$ через $i$.)
Приходим к противоречию с предположением (5.7).  Действительно,
в силу теоремы 5.5.4, леммы 5.5.5 и леммы 5.5.6,   при $x\in A$
$$
   {\cal C}S^{q_i(x)}(x),S^{q_i(x)}(y),S^{q_i(x)}(z))\equiv
(I\otimes T^{F(y,q_i(x))}\otimes T^{F(z,q_i(x))}){\cal C}x,y,z),
$$
$$
\mu\otimes\mu( \{(y,z)\,:\,
\rho({\cal C}S^{q_{i}(x)}(x),S^{q_{i}(x)}(y),S^{q_{i}(x)}(z)), \mu^3) > c
\})\,\to\,0.
$$
Таким образом, $${\cal C}x,y,z)\equiv \mu^3.$$

\medskip
Доказательство теоремы 5.5.2. Общий случай аналогичен случаю $n=3$,
который мы излагаем ниже.   Рассмотрим некоторую предельную
меру $\nu$ для сдвигов диагональной меры:
$$
\nu (A_0\times A_1 \times A_2) =
\lim_{j\to\infty}
 \lambda (A_0 \cap R^{z_1(j)}A_1  \cap R^{z_2(j)}A_2).
$$
Так как $R$ коммутирует с преобразованием $I\times T$, то
по причинам, указанным в доказательстве теоремы 5.5.1,
получим, что почти все условные меры   $\nu_w =\nu_{(x_1,x_2,x_3)}$
лежат в классе   $M(2,3)$.   В силу инвариантности меры
$\nu$ относительно $R\times R\times R$ получим,
что функция ${\cal C}: X^3 \to M(2,3)$, определенная равенством
$${\cal C}x_1,x_2,x_3) = \nu_{(x_1,x_2,x_3)},$$ удовлетворяет
условиям теоремы 5.5.3.  Следовательно,  $\nu=\mu\otimes\mu\otimes\mu$.
Отсюда вытекает, что  для любых последовательностей
$z_1(j),z_2(j)\to\infty$ таких, что  $|z_1(j)-z_2(j)|\to\infty$,
выполнено
$$
 \lambda (A_0 \cap R^{z_1(j)}A_1  \cap R^{z_2(j)}A_2)\to
\lambda (A_0)  \lambda (A_1) \lambda (A_2).
$$
Таким образом, косое произведение обладает перемешиванием
кратности 2.
Случай $k>2$ рассматривается аналогично.

Для  доказательства сохранения свойства тензорной простоты
опять применим теорему 5.5.3. Рассмотрим случай $n=3$.
Пусть $\nu$ -- джойнинг класса $M(2,3)$ трех копий $R$,
где косое произведение $R$ действует на $(X\times Y, \mu\otimes\mu')$.
Чтобы оказаться в условиях теоремы 5.5.3, нужно убедиться в том, что
условные меры $\nu_{(x_1,x_2,x_3)},$ лежат  в классе $M(2,3)$.
Последнее  вытекает из равенства
$$
  \nu\left((A_1 \times Y)\times  (A_2 \times B_2) \times (A_3 \times B_3)
  \right)  =  $$ $$ =\ \mu(A_1)\mu(A_2)\mu(A_3)\mu'(B_2)\mu'(B_3)
$$
для любых  наборов измеримых  множеств $A_1,A_2,A_3, B_2,B_3$,
которое является  следствием тензорной простоты $S$.
Действительно, фактор $S$  первой копии $R$ будет независим
с произведением других копий
$R$  (см. лемму 1.1.1).
Таким образом,
$$
  \nu_{(x_1,x_2,x_3)}(Y\times  B_2 \times B_3) =
 \mu'(B_2)\mu'(B_3),
$$
т.е. $\nu_{(x_1,x_2,x_3)}$ принадлежат   классу $M(2,3)$.
Случай $n>3$ рассматривается аналогично.
 



\begin{thebibliography}{99}
\bibitem{An1}
Аносов Д.В.
                  Геодезические потоки  на римановых многообразиях
                  отрицательной кривизны.  Труды МИАН.  Т.90. 1967.
\bibitem{An2}
Аносов Д.В.  О вкладе Н.Н. Боголюбова в теорию динамических систем.
         УМН. {\bf 49}(1994) No 5. 5--20.
\bibitem{An3}
Аносов Д.В.  О спектральных кратностях в эргодической теории.
Современные проблемы математики. Выпуск 3. М.: МИРАН. 2003.
\bibitem{Ver}
Вершик А.М. Многозначные отображения с инвариантной мерой (полиморфизмы)
и марковские процессы. Зап. науч. сем. ЛОМИ.   {\bf 72}(1977). 26--62.

\bibitem{V}
 Вершик А.М., Корнфельд И.П., Синай Я.Г.  Общая эргодическая теория
групп преобразований с инвариантной мерой. Совр. проблемы математики.
Фундаментальные направления. Итоги науки и техники. Т.2. М.: ВИНИТИ,
5--111 (1985).
\bibitem{Gur}
            Гуревич Б.М. Энтропия  потока орициклов.
            ДАН СССР.  {\bf 136}(1961). No 4. 768--770.
\bibitem{GeF}
Гельфанд И.М., Фомин С.В. Геодезические потоки на многообразиях
постоянной отрицательной кривизны.
УМН. {\bf 7 }(1952). No 1. 118--137.

\bibitem{KSa} Каток А.Б., Сатаев  E.A.  Стандартность
автоморфизмов перекладываний отрезков и потоков на поверхностях.
{ Матем. заметки.} {\bf 20} (1976). No 4. 479--488.

\bibitem{KS}
Каток А.Б., Степин А.М. Аппроксимации в эргодической теории.
УМН. {\bf 22 }(1967). No 5. 81--106.
\bibitem{KS1}
Каток А.Б., Степин А.М. Метрические свойства гомеоморфизмов, сохраняющих меру.
УМН. {\bf 25 }(1970). No 2. 193--220.

\bibitem{KSS}
Каток А.Б., Синай Я.Г., Степин А.М.
Теория динамических систем и общих групп преобразований с
инвариантной мерой.  Итоги науки. Математический анализ.
Т.3.(1985). М.: ВИНИТИ.   5--111

\bibitem{Kry}
            Крыгин А.Б. Пример цилиндрического каскада с аномальными
            метрическими свойствами. Вестник МГУ сер.1. 1975. No 5. 26--32.

\bibitem{Leo}
Леонов В.П. Применения характеристических функционалов и семиинвариантов
в эргодической теории стационарных процессов.
ДАН СССР. \bf 133 \rm (1960). No 3. 523--526.

\bibitem{Oseledets} Оселедец В.И. О спектре эргодических автоморфизмов.
ДАН СССР. {\bf 168}(1966). No 5. 1009--1011.

\bibitem{Os}
Оселедец В.И. Автоморфизмы с простым непрерывным спектром без группового
свойства.
Матем. заметки.  {\bf 5}(1969).  No 3. 323--326.

\bibitem{Ose}
            Оселедец В.И. Две неизоморфные динамические системы с
            одинаковым простым непрерывным спектром.
            Функц. анализ и  его прил.   {\bf 5}(1971). No 3. 75--79.

\bibitem{Para}
Парасюк О.С. Потоки гороциклов на поверхностях постоянной
отрицательной кривизны.
УМН. 1953. {\bf 8}(1953). No 3. 125--126.


\bibitem{Pi}  Пинскер М.С. Динамические системы с вполне положительной
и нулевой энтропией.
ДАН СССР. {\bf 133}(1960). No 5. 1025--1026.


\bibitem{Pri} Приходько А.А. Стохастические конструкции
перемешивающих систем положительного локального ранга.
Матем. заметки.   {\bf 69}(2000).  No 2. 316--319.

\bibitem{PR98} Приходько А.А., Рыжиков В.В. Простые некоммутативные
контрпримеры в теории джойнингов.
         Вестник МГУ сер.1. (1998). No 4. 16--19.
\bibitem{Rok}
Рохлин В.А. Эндоморфизмы компактных коммутативных групп.
Изв. АН СССР. сер. матем.   {\bf 13}(1949). 329--340.
\bibitem{UMN89}
Рыжиков В.В.  Замечание о кратном перемешивании. УМН \bf 44 \rm
(1989). No 1, 251--252.
\bibitem{flow}
Рыжиков В.В. Связь перемешивающих свойств потока
с изоморфизмом входящих в него преобразований.
Матем. заметки  \bf 49\rm(1991). No 6, 98--106.
\bibitem{R0}
Рыжиков В.В.   Перемешивание, ранг и минимальное самоприсоединение
    действий с инвариантной мерой.
Матем. сборник. \bf 183\rm(1992). No 3.  133--160.
\bibitem{F96}
Рыжиков В.В.
О когомологичности коциклов, отвечающих эргодическим косым
произведениям.
    Функц. анализ и его прил.  \bf 30\rm (1996), No 1,  86-88.
\bibitem{R03}
Рыжиков В.В.
Частичное кратное перемешивание на подпоследовательностях
может различать автоморфизмы $T$ и $T^{-1}$.
Матем. заметки  \bf 76\rm (2003). No 6, 889--895.
\bibitem{Si}
Синай Я.Г. О свойствах спектров эргодических динамических систем.
ДАН СССР.\bf 150\rm  (1963). No 6. 1235--1237.
\bibitem{Sin}
Синай Я.Г. О слабом изоморфизме преобразований с инвариантной
мерой. Матем. сборник. \bf 63\rm(1964). No 1.  23--42.

\bibitem{Star} Старков А.Н. О кратном перемешивании однородных потоков.
            ДАН. {\bf 333}(1993). No 1. 28--31.
\bibitem{Stara} Старков А.Н.
{ Динамические системы на однородных пространствах.}
М.: ФАЗИС. 1999.
\bibitem{St66}
Степин А.М.  О свойствах спектров эргодических
динамических систем с локально компактным временем.
ДАН СССР. \bf 169 \rm
(1966). No 4. 773--776.

\bibitem{St67}
Степин А.М.  Спектр и аппроксимация метрических автоморфизмов
 периодическими преобразованиями.
    Функц. анализ и его прил.  \bf 1\rm (1967), No 2,  77-80.

\bibitem{St73} Степин А.М. О связи аппроксимативных и
спектральных свойств метрических автоморфизмов.
Матем. заметки  \bf 13\rm(1973). No 3, 403--409.

\bibitem{St86}
Степин А.М. Спектральные свойства типичных динамических систем.
Изв. АН СССР. сер. матем.  {\bf 50}(1986).  801--834.
\bibitem{Sch}
             Шнейберг И.Я. Нули интегралов  вдоль траекторий
             эргодических систем.
             Функц. анализ и  его прил. {\bf 19}(1985). No 2. 92--93.
\bibitem{Abdalaoui1}  El Abdalaoui E. H. { \'Etude  Spectral des
Transformations D'Ornstein}. Ph.D Thesis. Universit\'e de Rouen. 1998.
\bibitem{Adams}  Adams T. Smorodinsky's conjecture on rank one systems.
{ Proc. Amer. Math. Soc.} {\bf 126} (1998). 739--744.

\bibitem{Age}  Ageev O.N.
On ergodic transformations with homogeneous spectrum.
{  J. Dynamical and Control Systems.} {\bf 5} (1999). No 1. 149--152.
\bibitem{Anz}
 Anzai H.  Ergodic skew product transformations on the
torus. { Osaka math. journal.} {\bf 3}(1951). No 1.  83--99 .
\bibitem{Atk}
            Atkinson G. Recurrence of co--cycles   and random walks.
            J.London Math. Soc.  {\bf 13}(1976). 486--488.
\bibitem{BH}
Blum J.R., Hanson D.L. On the mean ergodic  theorem for
   subsequences. Bull. Amer. Math. Soc.\bf 66 \rm (1960). 308--311.
\bibitem{Bou}
Bourgain J. On the spectral type of Ornstein's class one transformations
Israel J. Math.  {\bf 84}(1993).   53--63.

\bibitem{Chacon0} Chacon R.V. Transformations having continuous spectrum.
{ J. Math. Mech.} {\bf 16} (1966). 399--415.

\bibitem{Chacon1} Chacon R.V.  Approximation and spectral multiplicity.
Lecture Notes in Math {\bf 160} (1970). 18--27.
Academic Press. New York. San Fransisco. London. 1976.


\bibitem{F}
   Ferenczi S. Systems of finite rank. Colloquium mathematicum.
   {\bf 73}. No 1. 35--65 (1997).

\bibitem{Ferenczi1} Ferenczi S.  Systemes localement de rang un.
{ Ann. Inst. H. Poincar\'e.} {\bf 20} (1984). 33--51.


\bibitem{Friedman}  Friedman N. { Introduction to Ergodic Theory.}
Van Nostrand Reinhold. New York. 1970.

\bibitem{FTh}
  Friedman N.A., Thomas E.S. Higher order sweeping out, Illinois
  J.Math. \bf 29 \rm (1985). 401--417.




\bibitem{Fur}
Furstenberg H. Disjointness in ergodic theory, minimal sets, and a
problem in Diophantine approximation.  Math. Systems Theory.
{\bf 1}(1967).  1--49.

\bibitem{Furstenberg} Furstenberg H.  { Recurrence in Ergodic Theory
and Combinatorial Number Theory.}  Princeton University Press. Princeton.
1981.

\bibitem{FKO}
  Furstenberg H., Katznelson Y.,  Ornstein D. The ergodic  theoretical
  proof of Szemeredi's theorem.
  { Bull. Amer. Math. Soc.}(N.S.) {\bf 7} (1982). 527--552.


\bibitem{GHR}
Glasner E., Host B., Rudolph D. Simple systems and their higher order
self--joinings. Israel J. Math. \bf 78\rm (1992). 131--142.

\bibitem{GW}
Glasner E., Weiss B. A simple weakly mixing
transformation with non--unique prime factors.
{ Amer. J. Math.} {\bf 116} (1994). 361--375.
\bibitem{G} Goodson G.R.
A survey of recent results in the spectral theory of ergodic
dynamical systems.
{  J. Dynamical and Control Systems}. {\bf 5} (1999). No 2. 173--226.

\bibitem{GoodsonL} Goodson G.R.,    Lema\'{n}czyk M. Transformations
conjugate to their inverses have even essential values. {  Proc.
Amer. Math. Soc. }   {\bf 124} (1996). 2703--2710.


\bibitem{Goodson1}  Goodson G.R.,  del Junco  A.,  Lema\'{n}czyk M.,
 Rudolph D. Ergodic transformations
conjugate to their inverses by involutions.
{  Ergod. Th. Dynam. Sys.}
 {\bf 24} (1995). 95--124.

\bibitem{GoodsonR}  Goodson G.R.,   Ryzhikov V.V.  Conjugations,
joinings and direct products of locally rank--one dynamical systems.
{  J. Dynamical and Control Systems. } {\bf 3} (1997). No 3. 321--341.

\bibitem{Hal}
Halmos P.R. Lecture on ergodic  theory. Chelsea. New York. 1960.

\bibitem{Hos}
Host B. Mixing of all orders and pairwise independent joinings
   of systems with singular spectrum. Israel J. Math. \bf 76
\rm (1991). 289--298.

\bibitem{Jun}
del Junco A. A family of counter--examples in ergodic theory.
Israel J. Math. \bf 44 \rm (1983). 160--188.

\bibitem{Junco7}
del Junco  A., Lema\'{n}czyk M. Generic spectral
properties of measure--preserving maps and applications. {  Proc. Amer.
Math. Soc.}   {\bf 115} (1992). 725--736.

\bibitem{JL}
  del Junco A., Lemanczyk M. { Simple systems are disjoint from Gaussian
systems}. { Studia Math}.  {\bf 113}(1999). 249--256.

\bibitem{JuP}
del Junco A., Park K. An example of a measure--preserving flow with minimal
self--joining. J.d'Analyse Math. \bf 42 \rm (1983). 199--211.

\bibitem{JR}
del Junco A., Rudolph D. On ergodic  action whose self--joinings are graphs.
{  Ergod. Th. Dynam. Sys.} \bf 7 \rm (1987). 531--557.

\bibitem{JRS}
       del Junco A., Rahe A.M., Swanson L. Chacon's automorphism
      has minimal self--joinings.
      { J. d'Analyse Math.} \bf 37 \rm (1980). 276--284.


\bibitem{Kal}
  Kalikow S.A. Twofold mixing implies threefold mixing for rank
                one transformations.
               {  Ergod. Th. Dynam. Sys.} \bf4 \rm (1984). 237--259.


\bibitem{Kin}
King J.L. Joining--rank and the structure of finite--rank mixing transformation.
{  J. Analyse Math.}.  {\bf 51}.  182--227 (1988).

\bibitem{King}
King J.L. Ergodic properties where order 4 implies infinite order.
{  Israel J. Math.} \bf 80\rm (1992). 65--86.


\bibitem{KiT}
King J. L., Thouvenot J.--P. A canonical structure theorem for finite
 joining--rank maps.
{  J. Analyse. Math. }\bf 56\rm (1991). 211--230.



\bibitem{Led}
Ledrappier F. Un champ marcovien peut \^etre d'entropie null
et m\'elangeant. C. R. Acad. Sci.Paris Ser. A \bf 287 \rm (1978). 561--563.


\bibitem{LPT}  Lema\'{n}czyk M., Parreau F., Thouvenot  J.-P.
Gaussian automorphisms whose ergodic self-joinings are Gaussian.
Fundamenta Math. \bf 164\rm (2000). 253-- 293.

\bibitem{Mar}
Marcus B. The horocycle flow is mixing of all degrees.
          Inv. Math. \bf 46 \rm (1978). 201--209.

\bibitem{Moz}
Mozes S. Mixing of all orders of Lie group actions. Invent. Math.
\bf 107 \rm (1992). 235--241.

\bibitem{Ornstein} D. S. Ornstein. On the root problem in ergodic theory.
{ Proc. 6th Berkeley Symposium on Math. Stats. Prob.} University of
California Press. Berkeley. (1970). 348-356.

\bibitem{Par}
Park K.   GL(2,Z) action on a two torus. Proc. Amer. Math. Soc.
        {\bf 144}(1992). No 4.  955--963.

\bibitem{Pr} Prikhodko A.A.
Ergodic joinings of $GL(n,\Z)$--actions on n--torus.
    J. Dynamical and Control Systems.
    (1999).   {\bf 5}. No 3. 385--395.

\bibitem{PR} Prikhodko A.A. Ryzhikov V.V.
     Disjointness of convolutions for Chacon's automorphism.
     Colloquium Mathematicum. {\bf 84/85}(2000).
     67--74

\bibitem{Ra}
  Ratner M. Horocycle flows, joinings and rigidity of products.
   Annals of Mathematics. \bf 118 \rm (1983).  277--313.
\bibitem{Rat}
            Ratner M. On Raghunathan's measure conjecture.
            Ann. Math. {\bf 134}(1991). 545--607.


\bibitem{R85}
Rudolph D. k--fold mixing lifts to weakly mixing isometric extension.
           Ergod. Th. Dynam. Sys. \bf 5 \rm (1985). 445--447.

\bibitem{Rud}
Rudolph D. An example of a measure--preserving map with minimal
           self--joinings, and applications.
           J.d'Analyse Math. \bf35 \rm (1979). 97--122.

\bibitem{Rudolph2}  Rudolph D.J.  {  Fundamentals of Measurable
Dynamics. } Oxford University Press. Oxford. 1990.

\bibitem{Sch}
Schmidt K. Mixing automorphisms of compact group and a theorem by
Kurt Mahler. Pacific J.Math. \bf 137 \rm (1989). 371--385.


\bibitem{St5}  Stepin A.M. Les spectres des systemes dynamique.
{  Actes, Congres Intern. Math. } {\bf 2}(1970), 941--946.

\bibitem{T75} Thouvenot J.--P.
       Une classe de systemes pour lesquels la conjecture de
       Pinsker est vraie. Israel J. Math. {\bf 21}(1975).  208--214.

\bibitem{Tho}
Thouvenot J.--P. Some properties and applications of joinings in ergodic
theory.
Ergodic Theory and Its Connections with Harmonic Analysis: Proceedings
of the Alexandria 1993 Conference, K.E. Petersen and I.A. Salama, eds., LMS
Lecture Note Series \bf 205\rm  (1995).  207 -- 235.

\bibitem{T00}
Thouvenot J.--P.  { Les systems simples sont disjoints de ceux
                qui sont infiniment divisibles  et
                plongeabls  dans un flot.  }
                Colloq. Math. Vol 84/85 (2000). part 2. 481--483.

\bibitem{Veech}  Veech W.A. A criterion for a process to be prime.
{  Monatshefte Math. } {\bf 94} (1982). 335--341.
\medskip
\begin{center}
 {\bf Работы автора по теме диссертации}
\end{center}
\bibitem{preprint}
Рыжиков В.В.  Перемешивание, ранг и минимальное самоприсоединение
   сохраняющих меру преобразований. Препринт ВИНИТИ. 1991. 1--68.
\bibitem{UMN91}
Рыжиков В.В.  Джойнинги динамических систем. Аппроксимации и перемешивание.
   УМН. \bf 46\rm (1991). No 5.   177--178.

\bibitem{MZ92}
Рыжиков В.В.   Стохастические сплетения и джойнинги динамических систем.
Матем. заметки.  \bf 52\rm (1992). No 3.     130--140.
\bibitem{F1}
Рыжиков В.В.   Джойнинги и кратное перемешивание действий конечного ранга.
    Функц. анализ и его прил.  {\bf 27}(1993). No 2.   63--78.
\bibitem{Izv}
Рыжиков В.В. Джойнинги, сплетения, факторы и перемешивающие свойства
динамических систем.
Изв. АН СССР. сер. матем.  {\bf 57}(1993). No 1. 102--128.
\bibitem{UMN94}
Рыжиков В.В.  Косые произведения и кратное перемешивание динамических систем.
   УМН. \bf 49\rm (1994). No 2.   163--164.
\bibitem{F2}
Рыжиков В.В.  Кратное перемешивание и локальный ранг динамических систем.
    Функц. анализ и его прил.  \bf 29\rm (1995). No 2. 88--91.
\bibitem{UMN95}
Рыжиков В.В.  Функциональный взгляд на теорему Семереди. Замечание о кратном
перемешивании потоков.
   УМН. \bf 50\rm (1995). No 6.  213-214.

\bibitem{1}
Ryzhikov V.V.  Stochastic intertwinings and multiple mixing of
              dynamical systems.
   J. Dynamical and Control Systems. {\bf 2}(1996). No 1. 1--19.
\bibitem{MZ96a}
Рыжиков В.В.   Типичность изоморфизма преобразований при изоморфизме
их декартовых степеней.
Матем. заметки.  \bf 59\rm (1996). No 4.  630-632.

\bibitem{MZ96}
Рыжиков В.В.   Четная и нечетная простота
 динамических систем с инвариантной мерой.
    Матем. заметки. {\bf 60}(1996).  No 3. 470--473.
\bibitem{2}
Ryzhikov V.V. Around simple systems. Induced joinings and multiple mixing.
J. Dynamical and Control Systems. {\bf 3}(1997). No 1. 111--127.
\bibitem{MS}
Рыжиков В.В.  Сплетения тензорных произведений и стохастический централизатор
динамических систем. Матем. сборник. {\bf 188}(1997). No 2. 67--94.
\bibitem{F3}
Рыжиков В.В.  Полиморфизмы, джойнинги и тензорная простота динамических систем.
    Функц. анализ и его прил. {\bf 31}(1997). No 2. 45--57.
\bibitem{TM}
Рыжиков В.В. Об асимметрии каскадов.  Труды МИРАН.  \bf 216\rm (1997). 154--157.
\bibitem{3}
Ryzhikov V.V.  Transformations having homogeneous spectra.
    J. Dynamical and Control Systems.
      {\bf 5}(1999). No 1.  145--148.
\bibitem{SRM}
Ryzhikov V.V.      Homogeneous  spectrum, disjointness of convolutions,
         and mixing properties of  dynamical systems.
Selected Russian Mathematics.   {\bf  1}(1999).   13--24.
\bibitem{F4}
Рыжиков В.В.  Проблема Рохлина о кратном перемешивании в классе действий
положительного локального ранга.
    Функц. анализ и его прил.  \bf 34\rm(2000). No 1.   90--93.
\bibitem{F5}
Рыжиков В.В.    О рангах эргодического автоморфизма $T\times T$.
    Функц. анализ и его прил.  \bf 35\rm(2001). No 2.   84--87.
\end{thebibliography}
\end{document}